\documentclass[twoside,a4paper]{article}
\usepackage[top=4.5cm, bottom=4.5cm, inner=3.5cm,outer=3.5cm]{geometry}
\usepackage[utf8]{inputenc}
\usepackage{amsthm}
\usepackage{mathtools}
\usepackage{csquotes}
\usepackage{amsmath}
\usepackage{graphicx}
\usepackage{dsfont}
\usepackage{soul}
\usepackage{amssymb}
\usepackage{multicol}
\usepackage{courier}
\usepackage{listings}
\usepackage{eurosym}
\usepackage{subfiles}
\usepackage{multicol}
\usepackage{pdfpages}
\usepackage{footnote}
\usepackage{color}
\definecolor{myblue}{rgb}{.25, .25, .9}
\definecolor{myred}{rgb}{.6, .4, .4}
\definecolor{myred2}{rgb}{.9, .1, .1}
\definecolor{mygreen}{rgb}{.25, .6, .5}
\usepackage{adjustbox}
\usepackage[colorlinks=true,
            linkcolor=myred2,
            urlcolor=myred,
            citecolor=myred]{hyperref}
\usepackage{pdfpages}
\usepackage[english]{babel}
\usepackage{array}
\usepackage{enumerate}
\usepackage{gensymb}
\usepackage[nottoc]{tocbibind}
\usepackage{lettrine}
\usepackage[normalem]{ulem}
\usepackage[font=small,labelfont=bf]{caption}
\usepackage{cancel}
\usepackage{nomencl}
\usepackage{mathrsfs}
\usepackage{changepage}
\usepackage{amsmath,stmaryrd,graphicx,float}
\numberwithin{equation}{section}

\usepackage{tikz-cd}
\usetikzlibrary{arrows}

\usepackage{subcaption} 

\usepackage{titlesec}

\newtheoremstyle{mystyle}
  {}
  {}
  {\itshape}
  {}
  {\bfseries}
  {.}
  { }
  {}

\theoremstyle{mystyle}

\newtheorem{theorem}{Theorem}[section]

\newtheorem{definition}[theorem]{Definition}
\newtheorem{lemma}[theorem]{Lemma}
\newtheorem{proposition}[theorem]{Proposition}
\newtheorem{corollary}[theorem]{Corollary}

\newtheorem{example}[theorem]{Example}
 
\newtheorem{algorithm}[theorem]{Algorithm} 
 
\newtheorem{assumption}[theorem]{Assumption}

\usepackage{etoolbox}

\DeclareMathAlphabet{\mathbbold}{U}{bbold}{m}{n}

\titleformat{\section}
{\large\filcenter\bfseries}{\llap{\thesection\hskip 9pt}}{0pt}{}
\titlespacing*{\section}
{0pt}{3.25ex plus 1ex minus .2ex}{1.5ex plus .2ex}

\titleformat{\subsection}[runin]
{\bfseries}{\llap{\thesubsection\hskip 9pt}}{0pt}{}
\titlespacing*{\subsection}
{0pt}{3.25ex plus 1ex minus .2ex}{1.5ex plus .2ex}

\titleformat{\subsubsection}[runin]
{\bfseries}{\llap{\thesubsubsection\hskip 9pt}}{0pt}{}
\titlespacing*{\subsubsection}
{0pt}{3.25ex plus 1ex minus .2ex}{1.5ex plus .2ex}

\titleformat{\paragraph}[runin]
{\bfseries}{\llap{\theparagraph\hskip 9pt}}{0pt}{}
\titlespacing*{\paragraph}
{0pt}{3.25ex plus 1ex minus .2ex}{1.5ex plus .2ex}

\usepackage[titles]{tocloft}

\usepackage{fancyhdr}
\usepackage{lastpage}
\fancypagestyle{plain}{
\fancyhf{}
\fancyhead[LE]{\thepage\quad\textbf{\nouppercase\leftmark}}
\fancyhead[RO]{\textbf{\nouppercase\rightmark}\quad\thepage}
}

\fancypagestyle{basicstyle}{%
\fancyhf{}
\fancyhead[LE]{\thepage}
\fancyhead[RO]{\thepage}
}

\DeclareMathOperator*{\subjectto}{\text{subject to}}
\DeclareMathOperator*{\minimize}{\text{minimize}}

\DeclareMathOperator*{\argmin}{argmin}

\usepackage{algorithm}
\usepackage{algorithmic}[1]

\usepackage{chngcntr}
\counterwithin{figure}{section}

\usepackage[backend=bibtex,
maxbibnames=99,
style=alphabetic,
bibencoding=ascii,
giveninits=true, 
]{biblatex}
\renewbibmacro{in:}{}

\usepackage{rotating}
\usepackage{mathtools}

\title{{\textbf{\Large Small errors in random zeroth-order optimization\\ are imaginary}}}
\author{{Wouter Jongeneel}$^{\dagger}$\quad Man-Chung Yue$^{\ddagger}$\quad Daniel Kuhn$^{\dagger}$}
\date{\small{
    $^{\dagger}$Risk Analytics and Optimization Chair, \'Ecole Polytechnique F\'ed\'erale de Lausanne, \{wouter.jongeneel,daniel.kuhn\}@epfl.ch \\
    $^{\ddagger}$Musketeers Foundation Institute of Data Science and Department of Industrial and Manufacturing Systems Engineering, The University of Hong Kong, mcyue@hku.hk}
  \\~\\ \today}

\begin{document}
\maketitle
\thispagestyle{empty}

\begin{abstract}
Most zeroth-order optimization algorithms mimic a first-order algorithm but replace the gradient of the objective function with some gradient estimator that can be computed from a small number of function evaluations. This estimator is constructed randomly, and its expectation matches the gradient of a smooth approximation of the objective function whose quality improves as the underlying smoothing parameter~$\delta$ is reduced. Gradient estimators requiring a smaller number of function evaluations are preferable from a computational point of view. While estimators based on a single function evaluation can be obtained by use of the divergence theorem from vector calculus, their variance explodes as~$\delta$ tends to~$0$. Estimators based on multiple function evaluations, on the other hand, suffer from numerical cancellation when~$\delta$ tends to~$0$. To combat both effects simultaneously, we extend the objective function to the complex domain and construct a gradient estimator that evaluates the objective at a complex point whose coordinates have small imaginary parts of the order~$\delta$. As this estimator requires only one function evaluation, it is immune to cancellation. In addition, its variance remains bounded as~$\delta$ tends to~$0$. We prove that zeroth-order algorithms that use our estimator offer the same theoretical convergence guarantees as the state-of-the-art methods. Numerical experiments suggest, however, that they often converge faster in practice. 
\end{abstract}

{\footnotesize{
}}


\section{Introduction}
We study optimization problems of the form 
\begin{equation}
\label{equ:opt:main}
    \minimize_{x\in \mathcal{X}}\quad f(x),
\end{equation}
where $f:\mathcal{D}\to \mathbb{R}$ is a {real analytic and thus smooth} objective function defined on an open set~$\mathcal{D}\subseteq \mathbb{R}^n$, and~$\mathcal{X}\subseteq \mathcal{D}$ is a non-empty closed feasible set. 
Throughout the paper we assume that problem~\eqref{equ:opt:main} admits a global minimizer~$x^{\star}$ and that the objective function~$f$ can only be accessed through a deterministic zeroth-order oracle, which outputs function evaluations at prescribed test points. Under this premise, we aim to develop optimization algorithms that generate a (potentially randomized) sequence of iterates~$x_1,x_2,\ldots, x_K\in\mathcal X$ approximating~$x^\star$. As they only have access to a zeroth-order oracle, these algorithms fall under the umbrella of \textit{zeroth-order optimization}, \textit{derivative-free optimization} or, more broadly, \textit{black-box optimization}, see, \textit{e.g.,}  \cite{ref:audet2017derivative}. {As we will explain below and in contrast to all prior work on zeroth-order optimization, we will assume that our zeroth-order oracle also accepts {\em complex} inputs beyond~$\mathcal D$.}

Zeroth-order optimization algorithms are needed when problem~\eqref{equ:opt:main} cannot be addressed with first- or higher-order methods. This is the case when there is no simple closed-form expression for~$f$ and its partial derivatives or when evaluating the gradient of~$f$ is expensive. In simulation-based optimization, for example, the function~$f$ can be evaluated via offline or online simulation methods, but its gradient is commonly inaccessible. Zeroth-order optimization algorithms can also be used for addressing minimax, bandit or reinforcement learning problems, and they lend themselves for hyperparameter tuning in supervised learning~\cite{ref:spall2005introduction,ref:Conn,ref:nesterov2017random}. As they can only access function values, zeroth-order optimization methods are inevitably somewhat crude. This simplicity is both a curse and a blessing. On the one hand, it has a detrimental impact on the algorithms' ability to converge to local minima, on the other hand---and this requires further formalization~\cite{ref:scheinberg2022finite}, it may enable zeroth-order methods to escape from saddle points and thus makes them attractive for non-convex optimization.

Zeroth-order optimization algorithms can be categorized into direct search methods, model-based methods and random search methods~\cite{ref:larson2019derivative}. Direct search methods evaluate the objective function at a set of trial points without the goal of approximating the gradient. A representative example of a direct search method is the popular Nelder–Mead algorithm~\cite{ref:nelder1965simplex}. 
Model-based methods use zeroth-order information acquired in previous iterations to calibrate a $C^{r}$-smooth model for some~$r\in\mathbb{Z}_{\geq 0}$ that approximates the black-box function~$f$ locally around the current iterate and then construct the next iterate via $r^{\mathrm{th}}$-order optimization methods. These approaches typically attain a higher accuracy than the direct and random search methods, and they have the additional advantage that function evaluations can be re-used. In general, however, they require at least ${O}(n)$ function evaluations in each iteration to construct a well-defined local model~\cite{ref:berahas2021theoretical}. 
Examples of commonly used models include polynomial models, interpolation models and regression models~\cite{ref:larson2019derivative}. In contrast to model-based methods, random search methods estimate the gradients of~$f$ at the iterates directly from finitely many function evaluations and use the resulting estimators as surrogates for the actual gradients in a first-order algorithm. More precisely, random search methods typically approximate~$f$ by a smooth function~$f_{\delta}$ that is close to~$f$ for small~$\delta$ and construct an unbiased estimator $g_{\delta}(x)$ for $\nabla f_{\delta}(x)$ by sampling~$f$ at test points in the vicinity of~$x$ \cite{ref:Flaxman,ref:nesterov2017random}. For many popular approximations~$f_\delta$ there exists~$p\geq 1$ such that~$\|\nabla f_{\delta}(x)-\nabla f(x)\|\leq O(\delta^{p})$. In analogy to the model-based methods, $g_{\delta}(x)$ can thus be used as a surrogate for the actual gradient in a first-order algorithm. A striking advantage of these random search methods over model-based methods is that the computation of~$g_{\delta}(x)$ requires only~${O}(1)$ function evaluations, yet at the expense of weaker approximation guarantees~\cite{ref:IEEEsurvey,ref:berahas2021theoretical,ref:scheinberg2022finite}.  In principle, the approximation quality of the surrogate gradients (and therefore also the convergence rate of the first-order method at hand) can be improved by reducing the smoothing parameter~$\delta$. As $g_{\delta}(x)$ is often reminiscent of a difference quotient with increment~$\delta$, however, its evaluation is plagued by numerical cancellation. This means that if~$\delta$ drops below a certain threshold, innocent round-off errors in the evaluations of~$f$ have a dramatic impact on the evaluations of~$g_{\delta}$. Hence, the actual numerical performance of a random search zeroth-order algorithm may fall significantly short of its theoretical performance~\cite{ref:shi2021numerical}, however, the awareness for this phenomenon seems to be lacking.


Inspired by techniques for numerically differentiating analytic functions, we propose here a new smoothed approximation~$f_{\delta}$ as well as a corresponding stochastic gradient estimator~$g_\delta$ that can be evaluated rapidly and faithfully for arbitrarily small values of~$\delta$ without suffering from cancellation effects. 
Integrating the new estimator into the gradient-descent-type algorithm
\begin{equation}
\label{equ:GD:alg}
     x_{k+1} \leftarrow x_k - \mu_k \cdot g_{\delta_k}(x_k)
\end{equation}
with adaptive stepsize~$\mu_k$ and smoothing parameter~$\delta_k$ gives rise to new randomized zeroth-order algorithms. The performance of such algorithms is measured by the  decay rate of the regret $R_K=\textstyle\frac{1}{K} \sum^{K}_{k=1}\mathbb{E}\left[f(x_k)-f(x^{\star}) \right]$ as $K$ grows. Here, {$x^\star$ is a global minimizer, and} the expectation $\mathbb{E}[\cdot]$ is taken with respect to the randomness introduced by the algorithm. Note that if~$f$ is convex, then Jensen's inequality ensures that the expected suboptimality gap (or expected optimization error) of the \textit{averaged} iterate $\Bar{x}_{K}=\frac{1}{K}\sum^{K}_{k=1}x_k$ satisfies $\mathbb{E}\left[f(\Bar{x}_{K})-f(x^{\star}) \right]\leq R_K$. The main goal of this paper is to understand how~$R_K$ scales with the total number~$K$ of iterations and with critical problem parameters such as the dimension of~$x$ or Lipschitz moduli of~$f$. Whenever possible (\textit{e.g.,} when $f$ is strongly convex\textbf{}), we also analyze the expected suboptimality gap $\mathbb{E}[f(x_K)-f(x^{\star})]$ of the last iterate~$x_K$. The scaling behavior of $R_K$ with respect to~$K$ reflects the algorithm's \textit{convergence~rate}. We will show that algorithms of the form~\eqref{equ:GD:alg} equipped with the new gradient estimator offer provable convergence rates, are numerically stable, and empirically outperform algorithms that exploit existing smoothed approximations both in terms of accuracy and runtime.

\paragraph*{Notation}
We reserve the symbol~$i=\sqrt{-1}$ for the imaginary unit. The real and imaginary parts of a complex number $z=a+ib$ for $a,b\in \mathbb{R}$ are denoted by $\Re(z)=a$ and $\Im(z)=b$. In addition, $V_n$ stands for the volume of the unit ball~$\mathbb{B}^n=\{x\in \mathbb{R}^n : \|x\|_2\leq 1\}$, and~$S_{n-1}$ stands for the surface area of the unit sphere~$\mathbb{S}^{n-1}=\{x\in \mathbb{R}^n:\|x\|_2=1\}$. The family of all $r$ times continuously differentiable real-valued functions on an open set~$\mathcal D\subseteq\mathbb R^n$ is denoted by $C^r(\mathcal{D})$, and the family of all real analytic functions on~$\mathcal D$ is denoted by~$C^\omega(\mathcal D)$. 


\subsection{Related work}
Given a deterministic zeroth-order oracle, one could address problem~\eqref{equ:opt:main} with a gradient-descent algorithm that approximates the gradient of~$f$ with a vector of coordinate-wise finite differences \cite{ref:kiefer1952stochastic,ref:kushner2003stochastic,ref:spall2005introduction,ref:berahas2021theoretical}. The corresponding finite-difference methods for zeroth-order optimization are reminiscent of inexact gradient methods~\cite{ref:Aspremont:2008smooth,ref:devolder2014first}. 
Maybe surprisingly, there is merit in using {\em stochastic} gradient estimates even if a {\em deterministic} zeroth-order oracle is available~\cite{ref:nesterov2017random}. The randomness not only helps to penetrate previously unexplored parts of the feasible set but also simplifies the convergence analysis. Specifically, if $f$ is convex, then it is often easy to show that $f(x_k)$ converges {\em in expectation} to the global minimum~$f(x^\star)$~\cite{ref:nesterov2017random}. 

Zeroth-order optimization algorithms that mimic gradient descent algorithms can be categorized by the number of oracle calls needed for a single evaluation of the gradient estimator. The most efficient algorithms of this kind make do with one single oracle call. Arguably the first treatise on zeroth-order optimization with a random single-point gradient estimator appeared in~\cite[\S~9.3]{nemirovsky1983problem}, where the objective function~$f(x)$ is approximated by the smoothed function $f_{\delta}(x) = \textstyle{V_n^{-1}}\int_{\mathbb{B}^n}f(x+\delta y)\, \mathrm{d}y$, and the degree of smoothing is controlled by the parameter~$\delta>0$. By leveraging the dominated convergence theorem and the classical divergence theorem, one can then derive the following integral representation for the gradient of $f_{\delta}(x)$,
\begin{equation*}
    \nabla f_{\delta}(x) = \textstyle\frac{n}{\delta} \textstyle\int_{ \mathbb{S}^{n-1}}f(x+ \delta y) y\, \sigma(\mathrm{d} y),
\end{equation*}
where $\sigma$ represents the uniform distribution on the unit sphere~$\mathbb{S}^{n-1}$  (see also the proof of Proposition~\ref{prop:CR:grad}). Hence, the gradient of the smoothed function~$f_\delta$ admits the unbiased stochastic estimator 
\begin{equation}
    \label{equ:g:rho:N}
    g_{\delta}(x) = \textstyle\frac{n}{\delta }f(x+\delta y)y \quad \text{with}\quad y\sim \sigma,
\end{equation}
which can be accessed with merely a single function evaluation. Stochastic gradient estimators of this kind have been used as surrogate gradients in gradient descent algorithms, for example, in the context of bandit problems~\cite{ref:Flaxman}. However, as already pointed out in~\cite{nemirovsky1983problem}, the variance of the gradient estimator~\eqref{equ:g:rho:N} is of the order~$O(n^2/\delta^2)$ for small~$\delta$ even if the function~$f$ is constant. This is inconvenient because a smaller $\delta$ reduces the bias of~$f_\delta$ {\em vis-\`a-vis}~$f$. To improve this bias-variance trade-off, it has been proposed to subtract from~$g_\delta(x)$ the control variate $n\delta^{-1}f(x)y$, which has a vanishing mean but is strongly correlated with~$g_\delta(x)$ and therefore leads to a variance reduction~\cite{ref:agarwal2010optimal,ref:nesterov2017random}. The resulting unbiased stochastic gradient is representable as
\begin{equation}
    \label{equ:g:rho:N:nest}
    g'_{\delta}(x) = \textstyle\frac{n}{\delta }\left(f(x+\delta y)-f(x)\right)y\quad \text{with}\quad y\sim \sigma,
\end{equation}
which is reminiscent of a directional derivative and can be accessed via two function evaluations. Now, under mild conditions on $f$, the variance of~$g'_{\delta}(x)$ remains bounded as~$\delta$ tends to~$0$. If we aim to solve problem~\eqref{equ:opt:main} to an arbitrary precision, however, the smoothing parameter~$\delta$ needs to be made arbitrarily small, in which case $f(x+\delta y)$ and $f(x)$ become numerically indistinguishable. Subtractive cancellation therefore makes it impossible to evaluate estimators of the form~\eqref{equ:g:rho:N:nest} to an arbitrarily high precision. This phenomenon is exacerbated when the function evaluations are noisy, which commonly happens in simulation-based optimization~\cite{ref:Lian2016}. Generalized stochastic gradient estimators requiring multiple function evaluations are discussed in~\cite{ref:hazan2014bandit}, and in~\cite{duchi2015optimal,ref:lam2021minimax} various optimality properties of zeroth-order schemes with multi-point gradient estimators are discussed. 
 
Stochastic gradient estimators akin to~\eqref{equ:g:rho:N:nest} with~$u$ following a Gaussian instead of a uniform distribution are studied in~\cite{ref:nesterov2017random}. The corresponding stochastic gradient descent algorithms may converge as fast as~$O(n/K)$ if~$f$ is convex and has a Lipschitz continuous gradient, but they are typically~$O(n)$ times slower than their deterministic counterparts. Convergence can be accelerated by leveraging central finite-difference schemes or by adding random perturbations to the gradient estimators~\cite{duchi2015optimal,ref:gasnikov2017stochastic,ref:shamir2017optimal}. Local convergence results for nonconvex optimization problems are investigated in~\cite{ref:ghadimi2013stochastic}, and second-order algorithms similar to~\eqref{equ:GD:alg}, which use a Stein identity to estimate the Hessian matrix, are envisioned in~\cite{ref:balasubramanian2018zeroth}. Lower bounds on the convergence rate of algorithm~\eqref{equ:GD:alg} are established in~\cite{ref:agarwal2009information,ref:JamiesonNR12,shamir2013complexity}.

Another stream of related research investigates zeroth-order optimization methods that have only access to a {\em stochastic} zeroth-order oracle, which returns function evaluations contaminated by noise. The performance of these methods critically depends on the smoothness properties of~$f$. Indeed, the higher its degree of smoothness, the more terms in the Taylor series of~$f$ can be effectively averaged out~\cite{polyak1990optimal}. Improved convergence results for zeroth-order optimization methods under convexity assumptions are derived in~\cite{ref:bach:smooth,ref:akhavan2020exploiting,ref:novitskii2021improved}. When function evaluations are noisy, the smoothing parameter~$\delta$ controls a bias-variance tradeoff. Indeed, reducing~$\delta$ reduces the bias introduced by smoothing~$f$, while increasing~$\delta$ reduces the variance of the gradient estimator induced by the noisy oracle, which scales as~$1/\delta$ for small~$\delta$. The variance can be further reduced by mini-batching~\cite{ref:ji2019improved}. The impact of exact line search methods and adaptive stepsize selection schemes is discussed in~\cite{ref:stich2013optimization,ref:berahas2019linesearch}. Better stepsize rules are available if~$f$ displays a latent low-dimensional structure~\cite{ref:golovin2019gradientless}.

Generalized zeroth-order methods for optimizing functions defined on Riemannian manifolds are proposed in~\cite{ref:li2020stochastic}, and algorithms that have only access to a {\em comparison oracle}, which is less informative than a zeroth-order oracle, are investigated in~\cite{ref:cai2021onebit}.


For comprehensive surveys of zeroth-order optimization and derivative-free optimization we refer to~\cite{ref:larson2019derivative,ref:IEEEsurvey}. Abstract zeroth-order methods for convex optimization are studied in~\cite{ref:hu2016bandit}. The minimax regret bounds derived in this work reveal the importance of having control over the randomness of the zeroth-order oracle. Accordingly, most existing methods rely on the assumption that the noise distribution is light-tailed. In contrast, if the zeroth-order oracle is affected by adversarial noise, then optimization is easily obstructed~\cite[Thm~3.1]{ref:singer2015information}.

\subsection{Contributions}
Most existing zeroth-order schemes approximate the gradient of~$f$ in a way that makes them susceptible to numerical instability. For example, if~$f\in C^1(\mathbb{R})$ is Lipschitz continuous with Lipschitz constant $L$, then, in theory, the finite-difference approximation $(f(x+\delta)-f(x))/\delta$ converges to $\partial_x f(x)$ as $\delta>0$ tends to zero. In practice, however, $f$ can only be evaluated to within machine precision, which means that~$f(x+\delta)$ and $f(x)$ become indistinguishable for sufficiently small $\delta$. More precisely, as $f$ is Lipschitz continuous, we have $|f(x+\delta)-f(x)|\leq L\cdot |\delta| $, and thus cancellation errors are prone to occur when $L\cdot |\delta|$ approaches machine precision \cite[\S~11]{ref:overton2001numerical}. Other gradient estimators that are based on multiple function evaluations or that involve interpolation schemes suffer from similar cancellation errors. Nevertheless, the convergence guarantees of the corresponding zeroth-order methods require that the smoothing parameter~$\delta$ must be driven to zero. For example, \cite[Thm.~3.1]{ref:akhavan2020exploiting} establishes regret bounds under the assumption that the smoothing parameter of a multi-point estimator scales as~$\delta_k = O(1/\sqrt{k})$. 

The randomized gradient estimator~\eqref{equ:g:rho:N} avoids cancellation errors because it requires only one single function evaluation---an attractive feature that has, to the best of our knowledge, gone largely unnoticed to date. However, as pointed out earlier, the variance of this estimator diverges as~$\delta$ decays, which leads to suboptimal convergence rates. In this paper we propose a numerically stable gradient estimator that enables competitive convergence rates and is immune to cancellation errors. More precisely, we will use complex arithmetic to construct a one-point estimator akin to~\eqref{equ:g:rho:N} that offers similar approximation and convergence guarantees as state-of-the-art two-point estimators. Maybe surprisingly, we will see that computing this new estimator is not significantly more expensive than evaluating~\eqref{equ:g:rho:N}. Our results critically rely on the assumption that the objective function~$f$ is real analytic on its domain~$\mathcal D$. Recall that~$f$ is real analytic if it locally coincides with its multivariate Taylor series. We emphasize that real analyticity does not imply $\beta^{\mathrm{th}}$-order smoothness for some~$\beta\in\mathbb{Z}_{>0}$ in the sense of~\cite[\S~1.1]{ref:bach:smooth}, which means that~$f$ is almost surely $\beta-1$ times differentiable and that the $(\beta-1)^{\mathrm{th}}$-order term of its Taylor series is globally Lipschitz continuous. We will recall that~$f$ can be extended to a complex analytic function~$f:\Omega\rightarrow\mathbb C$ defined on some open set~$\Omega \subseteq \mathbb{C}^n$ that covers~$\mathcal D\subseteq\mathbb R^n$. By slight abuse of notation, this extension is also denoted by~$f$. Given an oracle that evaluates~$f$ at any query point in~$\Omega$, we will devise new zeroth-order methods that combine the superior convergence rates and low variances of multi-point schemes reported in~\cite{duchi2015optimal,ref:lam2021minimax} with the numerical robustness of single-point approaches.

We now use~$R=\|x_1-x^{\star}\|_2$ to denote the distance from the initial iterate~$x_1$ to a minimizer~$x^\star$ and~$F=f(x_1)-f(x^{\star})$ to denote the suboptimality of~$x_1$. Assuming that the objective function~$f$ is real analytic and has an $L$-Lipschitz continuous gradient, we will devise zeroth-order methods that offer the following convergence guarantees. {If~\eqref{equ:opt:main} represents a (constrained or unconstrained) convex optimization problem with $x^{\star}\in \mathrm{int}(\mathcal{X})$}, then our algorithm's regret decays as $O(nLR^2/K)$ with the iteration counter~$K$. If, in addition, $f$ is $\tau$-strongly convex for some~$\tau>0$, then the expected suboptimality decays at the linear rate~$O((1-\tau/(4nL))^K L R^2)$. {If~\eqref{equ:opt:main} represents a non-convex optimization problem, finally, we establish local convergence to a stationary point and prove that~$\min_{k\in [K]} \mathbb{E}[\|\nabla f(x_k)\|_2^2]\leq O(nLF/K)$}. All of these convergence rates are qualitatively equivalent to the respective rates reported in~\cite[Thm.~8]{ref:nesterov2017random}, and they are sharper than the rates provided in~\cite[\S~3]{ref:akhavan2020exploiting} in the noise-free limit. {The latter rely on higher-order smoothness properties of~$f$ but do not require $f$ to be analytic.} The key difference to all existing methods is that we can drive the smoothing parameter to~$0$,~\textit{e.g.}, as~$\delta_k=\delta/k$, without risking numerical instability. 

As highlighted in the recent survey article~\cite{ref:larson2019derivative}, an important open question in zeroth-order optimization is whether single-point estimators enable equally fast convergence rates as multi-point estimators. The desire to reap the benefits of multi-point estimators at the computational cost of using single-point estimators has inspired multi-point estimators with memory, which only require a single new function evaluation per call~\cite{ref:resfeedback2022}. However, this endeavor has not yet led to algorithms that improve upon the theoretical and empirical performance of the state-of-the-art methods in~\cite{ref:nesterov2017random}. Filtering techniques inspired by ideas from extremum seeking control can be leveraged to improve the convergence rates obtained in~\cite{ref:resfeedback2022} to~$O(n/K^{2/3})$ \cite{ref:chen2021improve}. However, this rate is still inferior to the ones reported in~\cite{ref:nesterov2017random}. To our best knowledge, we propose here the first single-point zeroth-order algorithm that enjoys the same convergence rates as the multi-point methods in~\cite{ref:nesterov2017random} but often outperforms these methods in experiments. {The price we pay for these benefits is the assumption that there exists a zeroth-order oracle accepting complex queries. This assumption is restrictive as it rules out oracles that depend on performing a physical experiment or timing a computational run etc. 
Nevertheless, as we will see in Section~\ref{sec:num}, the approach can excel in the context of simulation-based optimization.}

Numerical experiments built around standard test problems as well as a model predictive control (MPC) problem corroborate our theoretical results and demonstrate the practical efficiency of the proposed algorithms. Although cancellation effects are caused by rounding to machine precision, which is nowadays of the order $10^{-16}$, our single-point gradient estimator improves both the accuracy as well as the speed of zeroth-order algorithms already when $\varepsilon$-optimal solutions with $\varepsilon\gg 10^{-16}$ are sought. 


\paragraph*{Structure}
Section~\ref{sec:prelim} reviews basic tools from multivariate complex analysis and introduces the complex-step method from numerical differentiation. Section~\ref{sec:oracle} then combines smoothing techniques with complex arithmetic to construct a new single-point gradient estimator, and Sections~\ref{sec:convex}-\ref{sec:nonconvex} analyze the favorable convergence rates of zeroth-order optimization methods equipped with the new gradient estimator in the context of convex, strongly convex and non-convex optimization, respectively. Section~\ref{sec:num} reports on numerical experiments, and Section~\ref{sec:conclusions} concludes. 

\section{Preliminaries}
\label{sec:prelim}
Before presenting our main results, we review some tools that may not usually belong to the standard repertoire of researchers in optimization. Specifically, 
Section~\ref{sec:complex_analysis} reviews the relevant basics of multivariate complex analysis, Section~\ref{sec:Im:trick} introduces the complex-step approach, which uses complex arithmetic for computing highly precise and numerically stable approximations of derivatives based on a single function evaluation, and Section~\ref{sec:lipschitz} provides a survey of inequalities that will be needed for the analysis of the algorithms proposed in this paper.

\subsection{Multivariate complex analysis}
\label{sec:complex_analysis}
For any multi-index $\alpha \in \mathbb{Z}_{\geq 0}^n$ and vector~$x\in\mathbb R^n$, we use~$x^\alpha$ as a shorthand for the monomial~$x_1^{\alpha_1}\cdots x_n^{\alpha_n}$, and we denote the degree of~$x^\alpha$ by $|\alpha|=\sum_{i=1}^n \alpha_i$. The factorial of~$\alpha$ is defined as~$\alpha!=\prod^n_{i=1}\alpha_i !$, and~$\partial_x^{\alpha}$ stands for the higher-order partial derivative~$\partial_{x_1}^{\alpha_1} \cdots \partial_{x_n}^{\alpha_n}$. Multi-index notation facilitates a formal definition of real analytic functions. 

\begin{definition}[Real analytic function]
\label{def:real-analytic}
The function $f:\mathcal D\to\mathbb R$ is real analytic on~$\mathcal D\subseteq \mathbb R^n$, denoted $f\in C^{\omega}(\mathcal{D})$, if for every~$x'\in\mathcal D$ there exist~$f_\alpha\in\mathbb R$,  $\alpha \in \mathbb{Z}_{\geq 0}^n$, and an open set $U\subseteq \mathcal{D}$ containing $x'$ such that
\begin{equation}
	\label{eq:power-series-real}
    	f(x) = \textstyle\sum_{\alpha\in \mathbb{Z}_{\geq 0}^n} f_\alpha \cdot(x-x')^{\alpha}\quad \forall x\in U. 
\end{equation}
\end{definition}

Whenever we write that a series has a finite value, we mean that it converges absolutely, that is, it converges when the summands of the series are replaced by their absolute values. In this case any ordering of the summands results in the same value.

One can show that any real analytic function is {infinitely  differentiable} and that the coefficients of its power series are given by~$f_\alpha=\frac{1}{\alpha!}\, \partial_x^{\alpha}f(x')$ for every~$\alpha \in \mathbb{Z}_{\geq 0}^n$. This implies that the power series is unique and coincides with the multivariate Taylor series of~$f$ around~$x'$ \cite[\S~2.2]{ref:Real_Analytic_02}. We will now recall that every real analytic function admits a complex analytic extension. 

\begin{definition}[Complex analytic function]
{The function $f:\Omega\to\mathbb C$ is complex analytic on~$\Omega\subseteq \mathbb C^n$, denoted $f\in H(\Omega)$, if for every~$z'\in\Omega$ there exist~$f_\alpha\in\mathbb C$,  $\alpha \in \mathbb{Z}_{\geq 0}^n$, and an open set $U\subseteq \Omega$ containing $z'$ such that}
\begin{equation}
	\label{eq:power-series}
    	f(z) = \textstyle\sum_{\alpha\in \mathbb{Z}_{\geq 0}^n} f_\alpha \cdot (z-z')^{\alpha}\quad \forall z\in U. 
\end{equation}
\end{definition}

Complex analytic functions are intimately related to holomorphic functions.

\begin{definition}[Holomorphic function]
The function $f:\Omega\to\mathbb C$ is holomorphic on an open set~$\Omega\subseteq \mathbb C^n$ if the complex partial derivatives~$\partial_{z_j}f$, $j=1,\ldots,n$, exist and are finite at every~$z\in\Omega$.
\end{definition}
{The requirement that $\Omega$ be open is essential, and $f$ may fail to be holomorphic on a neighborhood of a point~$z$ even if it is complex differentiable at~$z$. For example, the Cauchy-Riemann equations reviewed below imply that $f(z)=|z|^3$ is complex differentiable at $z=0$ but fails to be complex differentiable on any neighborhood of~$0$.} Holomorphic functions are in fact infinitely often differentiable \cite[Prop.~1.1.3]{ref:lebl2019tasty}. Moreover, a function is holomorphic if and only if it is complex analytic \cite[Thm.~1.2.1]{ref:lebl2019tasty}. 

It is common to identify any complex vector~$z\in\mathbb C^n$ with two real vectors~$x,y\in\mathbb R^n$ through~$z=x+iy$. Similarly, we may identify any complex function~$f:\mathbb C^n\to\mathbb C$ with two real functions~$u:\mathbb R^n\to\mathbb R$ and $v:\mathbb R^n\to\mathbb R$ through the relation~$f(x+iy)=u(x,y)+i v(x,y)$. Clearly, $u$ and~$v$ inherit the differentiability properties of~$f$ and vice versa. In particular, one can show that if~$f$ is holomorphic, then the partial derivatives of~$u$ and~$v$ exist and satisfy the multivariate Cauchy-Riemann equations.

\begin{theorem}[Multivariate Cauchy-Riemann equations]
\label{thm:multivar:CR}
\!If $f(x+iy)=u(x,y)+iv(x,y)$ is a holomorphic function on an open set~$\Omega \subseteq \mathbb C^n$, then the multivariate Cauchy-Riemann equations
\begin{equation}
    \label{equ:multivar:Cr}
    \partial_{x_j} u(x,y) = \partial_{y_j} v(x,y) \quad \text{and}\quad -\partial_{x_j} v(x,y) = \partial_{y_j} u(x,y) \quad \forall j=1,\dots,n
\end{equation}
 hold for all~$x,y\in\mathbb R^n$ with~$x+iy\in \Omega$.
\end{theorem}

Theorem~\ref{thm:multivar:CR} is a standard result in complex analysis; see, {\em e.g.}, \cite[Thm.~11.2]{ref:rudin:CR} or~\cite{ref:Krantz:Complex}. 
{Nevertheless, we provide here a short proof to keep this paper self-contained.
\begin{proof}[Proof of Theorem~\ref{thm:multivar:CR}]
We use $e_j$ to denote the $j^{\mathrm{th}}$ standard basis vector in~$\mathbb{R}^n$. By the definition of the complex partial derivative, for any~$z\in\Omega$ we have
\begin{equation*}
    \partial_{z_j}f (z)= \lim_{{\delta \in \mathbb C,\, \delta \to 0}}\tfrac{1}{\delta}({f(z+\delta e_j)-f(z)}) ,
\end{equation*}
where the limit exists and is independent of how~$\delta\in\mathbb C$ converges to~$0$ because~$f$ is holomorphic on~$\Omega$. In particular, $\delta$ may converge to~$0$ along the real or the imaginary axis without affecting the result. Using our conventions that~$z=x+iy\in \Omega$ and~$f(x+iy)=u(x,y)+i v(x,y)$, we thus have
\begin{align*}
    \partial_{x_j}(u(x,y)+i v(x,y)) & =  \lim_{{\delta \in \mathbb R, \, \delta \to 0}}\frac{f((x+\delta e_j) + iy)-f(x+i y)}{\delta} \\ & =  \lim_{{\delta \in \mathbb R, \, \delta \to 0}}\frac{f(x+i (y+\delta e_j))-f(x+i y)}{i\delta}\\
    &= \tfrac{1}{i}\partial_{y_j}(u(x,y)+i v(x,y))
\end{align*}
for all $x,y\in\mathbb R^n$, where the second equality holds because both limits are equal to~$\partial_{z_j}f(z)$. Matching the real and imaginary parts of the above equations yields~\eqref{equ:multivar:Cr}.
\end{proof}}
Under additional assumptions one can further show that the Cauchy-Riemann equations imply that~$f$ is holomorphic~\cite{ref:conv:Holomorphic}. However, this reverse implication will not be needed in this paper.
The following lemma based on~\cite[\S~2.3]{ref:Krantz:Complex} establishes that any real analytic function defined on an open set~$\mathcal D\subseteq \mathbb R^n$ admits a complex analytic extension defined on an open set~$\Omega\subseteq \mathbb C^n$ that covers~$\mathcal D$. 

\begin{lemma}[Complex analytic extensions]\label{lem:analytic_extension}
If~$f\in C^{\omega}(\mathcal{D})$, then there exists an open set $\Omega \subseteq \mathbb{C}^n$ and a complex analytic function~$g\in H(\Omega)$ such that $\mathcal{D} \subseteq \Omega$ and $f(x)=g(x)$ for every~$x\in\mathcal D$, with $\mathcal{D}$ understood as embedded in $\mathbb{C}^n$. 
\end{lemma}

{
\begin{proof}
Select any~$x'\in\mathcal D$. As $f\in C^{\omega}(\mathcal{D})$, there exists a neighborhood~$U\subseteq \mathcal D$ of~$x'$ such that~$f$ admits a power series representation of the form~\eqref{eq:power-series} on~$U$. Also, as~$U$ is open, there exists~$x\in U$ with~$r_j=|x_j-x'_j|>0$ for every~$j=1,\ldots,n$. By Abel's lemma~\cite[Prop.~2.3.4]{ref:Krantz:Complex}, the power series~\eqref{eq:power-series} extended to~$\mathbb C^n$ is thus guaranteed to converge on the open polydisc $\Delta(x')=\{z\in\mathbb C^n:  |z_j-x'_j|< r_j \,\forall j=1,\ldots,n\}$.
This reasoning implies that~$f$ extends locally around~$x'$ to a complex analytic function, which we henceforth denote as~$g_{x'}$. It remains to be shown that the local extensions corresponding to different reference points~$x'\in\mathcal D$ are consistent. To this end, select any~$x',x''\in\mathcal D$ such that the polydiscs~$\Delta(x')$ and~$\Delta(x'')$ overlap. We need to prove that~$g_{x'}$ and~$g_{x''}$ coincide on the open convex set~$\Delta=\Delta(x')\cap \Delta(x'')$, which has a non-empty intersection with~$\mathbb R^n$. For ease of exposition, we will equivalently prove that the holomorphic function~$h=g_{x'}-g_{x''}$ vanishes on~$\Delta$. We first notice that~$h$ vanishes on~$\Delta\cap\mathbb R^n$ because~$g_{x'}$ and~$g_{x''}$ are constructed to coincide with~$f$ on~$\Delta\cap\mathbb R^n$. This implies that~$\partial_{z_j}h=\partial_{x_j} h=0$ on~$\Delta\cap\mathbb R^n$, where the first equality follows from standard arguments familiar from the proof of Theorem~\ref{thm:multivar:CR}. As any partial derivative of a holomorphic function is also holomorphic, one can use induction to show that all higher-order partial derivatives of~$h$ must vanish on~$\Delta\cap\mathbb R^n$. Hence, the Taylor series of~$h$ around any reference point in~$\Delta\cap\mathbb R^n$ vanishes, too. We have thus shown that~$h=0$ on an open subset of~$\Delta$. By standard results in complex analysis, this implies that~$h$ vanishes throughout~$\Delta$; see, {\em e.g.}, \cite[Thm.~1.2.2]{ref:lebl2019tasty}. In summary, this reasoning confirms that all local complex analytic extensions~$g_{x'}$, $x'\in\mathcal D$, of $f$ are consistent and thus coincide with a complex analytic function~$g$ defined on the open set~$\Omega=\cup_{x'\in\mathcal D} \Delta(x')$. This observation completes the proof.
\end{proof}}

{Lemma~\ref{lem:analytic_extension} implies that the complex extension of a real analytic function is unique. For example, $f(x)=\log(x)$ is real analytic on the positive real line. Representing $z=re^{i\theta}\in\mathbb C$ in polar form with $r\geq 0$ and $\theta\in(-\pi,\pi]$, the complex logarithm has countably many branches, that is, $\log(z)$ can be defined as $g_k(z)=\log(r)+i(\theta+2\pi k)$ for any $k\in \mathbb{Z}$. However, only the branch $g_0$ corresponding to $k=0$ matches $f$ on the positive reals.} 
We will henceforth use the same symbol~$f$ to denote both the given real analytic function as well as its {unique} complex analytic extension~$g$. {We now explicitly derive the complex analytic extensions of a few simple univariate functions.
\begin{example}[Complex analytic extensions]
    The unique complex analytic extension of $f(x)= e^x$ is the entire function $g(z)=g(x+iy)=e^x(\cos(y)+i\sin( y))$. Similarly, the unique complex analytic extension of the even polynomial $f(x)=x^{2p}$ with $p\in\mathbb{Z}_{\geq 0}$ is the entire function
    \[
    \textstyle g(z)=g(x+iy)=\sum_{k=0}^p (-1)^k \binom{2p}{2k} y^{2k}x^{2(p-k)}+i\sum_{k=0}^{p-1} (-1)^k \binom{2p}{2k+1}y^{2k+1}x^{2(p-k)-1}.
    \]
    Finally, the unique solution $f(x)$ to the Lyapunov equation $f(x)=x^2f(x)+1$ parametrized by~$x\in\mathbb R$ is real analytic on $\mathbb R\setminus\{1\}$. It admits the extension 
    \[
         g(z)=g(x+iy)=\frac{1-x^2+y^2-2i xy}{(1-x^2+y^2)^2+4x^2y^2},
    \]
    which is analytic throughout~$\mathbb C\setminus\{(1,0)\}$.
\end{example}} 

{The next example shows that the domain~$\Omega$ of the complex analytic extension is not always representable as~$\mathbb{R}^n+ i\cdot(-\bar{\delta}, \bar{\delta})^n$ for some~$\bar{\delta}>0$ even if~$\mathcal D=\mathbb R^n$. 

\begin{example}[Non-trivial extension]
Consider the function $f(x) =\sum^{\infty}_{k=1} 2^{-k} (1+k^2(x-k)^2)^{-1}\in C^{\omega}(\mathbb{R})$, which admits a unique complex analytic extension with domain~$\Omega=\mathbb C\backslash \{k+ik^{-1}:k\in \mathbb{Z}_{>0}\}$. In addition, $f$ can be extended to a meromorphic function on~$\mathbb C$ with countably many poles~$k+ik^{-1}$, $k\in \mathbb{Z}_{>0}$. As these poles approach~$\mathbb R$ arbitrarily closely, however, $\Omega$ cannot contain any strip of the form $\mathbb{R}\times i\cdot (-\bar{\delta}, \bar{\delta})$.
\end{example}}
To avoid technical discussions of limited practical impact, we will from now on restrict attention to functions~$f\in C^\omega(\mathcal D)$ that admit a complex analytic extension to~$\Omega=\mathcal{D}\times i\cdot (-\bar{\delta}, \bar{\delta})^n$ for some~$\bar{\delta}>0$. One can show that such an extension always exists if~$f\in C^{\omega}(\mathbb{R}^n)$ and~$\mathcal{D}$ is bounded or if~$f$ is \textit{entire}, that is, if~$f$ has a globally convergent power series representation. The latter condition is restrictive, however, because it rules out simple functions such as~$f(x)=1/(1+x^2)$. Provided there is no risk of confusion, we will sometimes call a real analytic function~$f\in C^\omega$ and its complex analytic extension simply an {\em analytic} function.

\subsection{Complex-step approximation}
\label{sec:Im:trick}
The \textit{finite-difference} method~\cite[Ch.~3]{ref:vuik2007numerical} is arguably the most straightforward approach to numerical differentiation. {It simply approximates the derivative of any sufficiently smooth function $f\in C^2(\mathbb{R})$ by a difference quotient. For example, the forward-difference method uses the approximation
\begin{equation}
\label{equ:finite:difference:forward}
    \partial_x f(x) = \tfrac{1}{\delta}({f(x+\delta)-f(x)})+O(\delta).
\end{equation}
The continuity of the second derivative of~$f$ allows for a precise formula for the $O(\delta)$ remainder term. However, as explained earlier, the finite difference method suffers from cancellation errors when~$\delta$ becomes small.} The \textit{complex-step} approximation proposed in~\cite{Ref:Lyness:Moler, ref:SquireTrapp} and further refined in~\cite{ref:Martins,ref:Abreu:geo,ref:Abreu:wave:2018} leverages ideas from complex analysis to approximate the derivative of any real analytic function~$f\in C^{\omega}(\mathbb{R})$ on the basis of one single function evaluation only, thereby offering an elegant remedy for numerical cancellation. Denoting by~$u$ and~$v$ as usual the real and imaginary parts of the unique complex analytic extension of~$f$, which exists thanks to Lemma~\ref{lem:analytic_extension}, we observe that $\partial_xf(x)$ equals
\begin{equation*}
\label{equ:cs:1}
    \partial_x u(x,0)=\partial_y v(x,0)  = \lim_{\delta\downarrow 0} \tfrac{1}{\delta}({v(x,\delta)-v(x,0)}) = \lim_{\delta\downarrow 0} \tfrac{1}{\delta}{v(x,\delta)} = \lim_{\delta\downarrow 0}\tfrac{1}{\delta}{\Im\big(f(x+i\delta) \big)},
\end{equation*}
where the first and the fourth equalities hold because~$f(x)$ must be a real number, which implies that~$v(x,0)=0$, while the second equality follows from the Cauchy-Riemann equations. The derivative~$\partial_x f(x)$ can thus be approximated
by the fraction~$\Im(f(x+i\delta))/\delta$, which requires merely a single function evaluation. To estimate the approximation error, we consider the Taylor expansion
\begin{equation}
\label{equ:im:taylor}
    f(x+i\delta) = f(x) + \partial_x f(x) i\delta - \tfrac{1}{2}\partial_x^2 f(x) \delta^2 - \tfrac{1}{6}\partial_x^3 f(x) i\delta^3 + O(\delta^4)
\end{equation}
{of the unique complex analytic extension of~$f$, which exists thanks to Lemma~\ref{lem:analytic_extension}.} 
Separating the real and imaginary parts of~\eqref{equ:im:taylor} then yields
\begin{equation*}
\label{equ:imag:trick}
    f(x)=\Re(f(x+i\delta))+O(\delta^2)\quad \text{and} \quad \partial_x f(x) = \tfrac{1}{\delta}{\Im\big(f(x+i\delta) \big)}+O(\delta^2).
\end{equation*}
This reasoning shows that a single {\em complex} function evaluation~$f(x+i\delta)$ is sufficient to approximate both~$f(x)$ as well as~$\partial_x f(x)$ without the risk of running into numerical instability caused by cancellation effects. In addition, the respective approximation errors scale quadratically with~$\delta$ and are thus one order of magnitude smaller than the error incurred by~\eqref{equ:finite:difference:forward}. Note also that the complex-step approximation recovers the derivatives of quadratic functions {\em exactly} irrespective of the choice of~$\delta$. For example, if~$f(x)=x^2$, then~$\Im(f(x+i \delta))/\delta =2x= \partial_x f(x)$. This insight suggests that the approximation is numerically robust for locally quadratic functions. 


The error of the complex-step approximation can be further reduced to~$O(\delta^4)$ by enriching it with a finite-difference method~\cite{ref:Abreu:geo,ref:Hare2020}. However, the resulting scheme requires multiple function evaluations and is thus again prone to cancellation errors. Unless time is expensive, the standard complex-step approximation therefore remains preferable. The complex-step approximation can also be generalized to handle matrix functions~\cite{ref:Higham} or to approximate higher-order derivatives~\cite{ref:lantoine2012}. {Its ramifications for automatic differentiation (AD) are discussed in~\cite{ref:Martins}. We return to AD below.}

Being immune to cancellation effects, the complex-step approach offers approximations of almost arbitrary precision. For example, software by the UK's National Physical Laboratory is reported to use smoothing parameters as small as~$\delta=10^{-100}$~\cite[p.~44]{ref:cox2004software}. {The complex-step approach also emerges in various other domains. For example,} it is successfully used in airfoil design~\cite{ref:airfoil}. However, its potential for applications in optimization has not yet been fully exploited. Coordinate-wise complex-step approximations with noisy function evaluations show promising performance in line search experiments~\cite{ref:complexwithnoise} but come without a rigorous convergence analysis. In addition, the complex-step approach is used to approximate the gradients and Hessians in deterministic Newton algorithms for blackbox optimization models~\cite{ref:Hare2020}. The potential of leveraging complex arithmetic in mathematical optimization is also mentioned in~\cite{ref:snyman2018practical,ref:Berahas}. In this paper we use the complex-step method to construct an estimator akin to~\eqref{equ:g:rho:N} and provide a full regret analysis. Our approach is most closely related to the recent works~\cite{ref:wang2021improved,ref:wang2021model}, which integrate the complex-step and simultaneous perturbation stochastic approximations~\cite{ref:spall1992multivariate} into a gradient-descent algorithm and offer a rigorous asymptotic convergence theory. In contrast, we will derive convergence \textit{rates} for a variety of zeroth-order optimization problems.

In optimization, the ability to certify that the gradient of an objective function is sufficiently small ({\em i.e.}, smaller than a prescribed tolerance) is crucial to detect local optima. The following example shows that, with the exception of the complex-step approach, standard numerical schemes to approximate gradients fail to offer such certificates---at least when a high precision is required.

\begin{figure*}[t!]
    \centering
    \begin{subfigure}[b]{0.3\textwidth}
        \includegraphics[width=\textwidth]{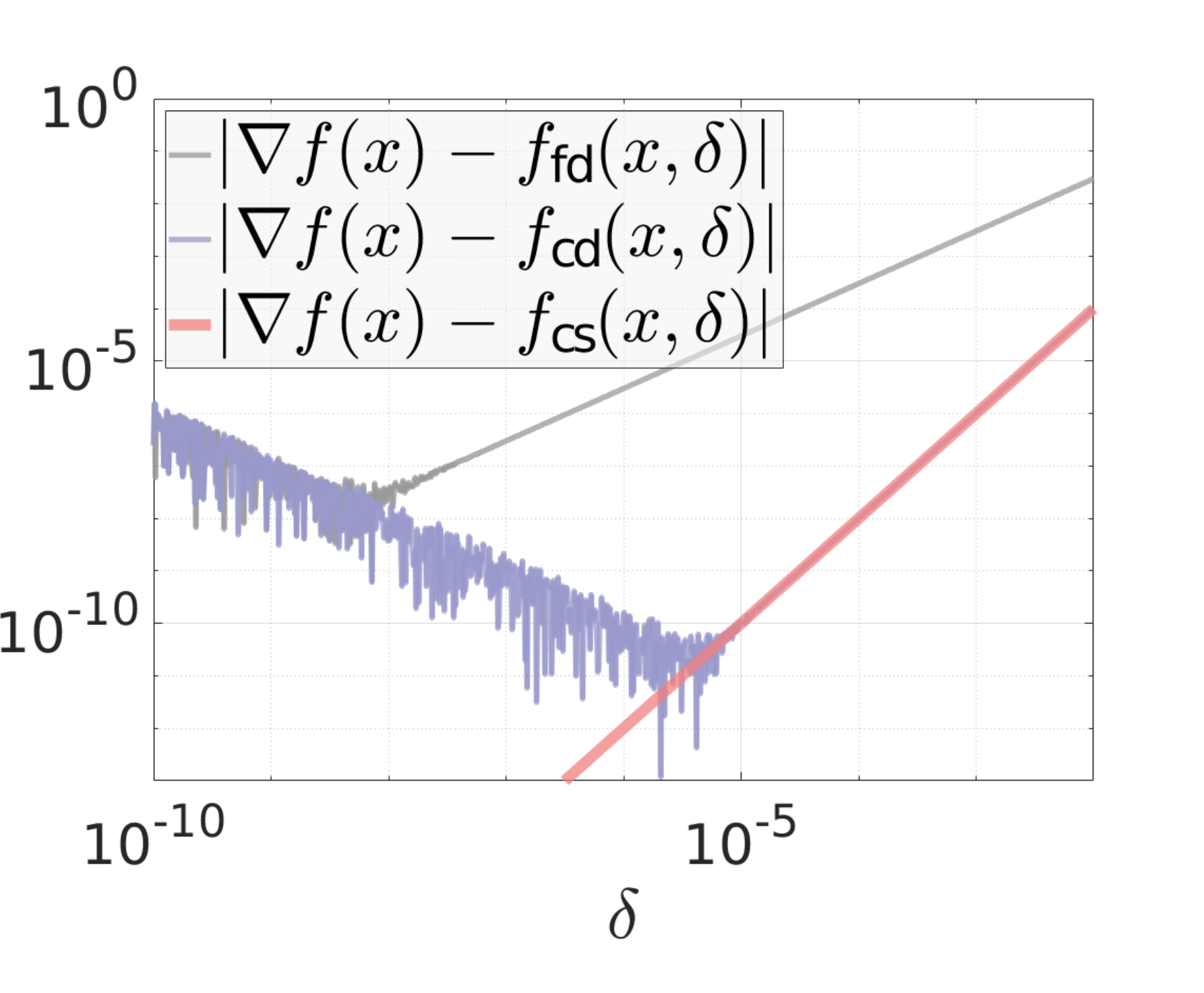}
        \caption{$x=-1$}
        \label{fig:x0-1}
    \end{subfigure}\quad
    \begin{subfigure}[b]{0.3\textwidth}
        \includegraphics[width=\textwidth]{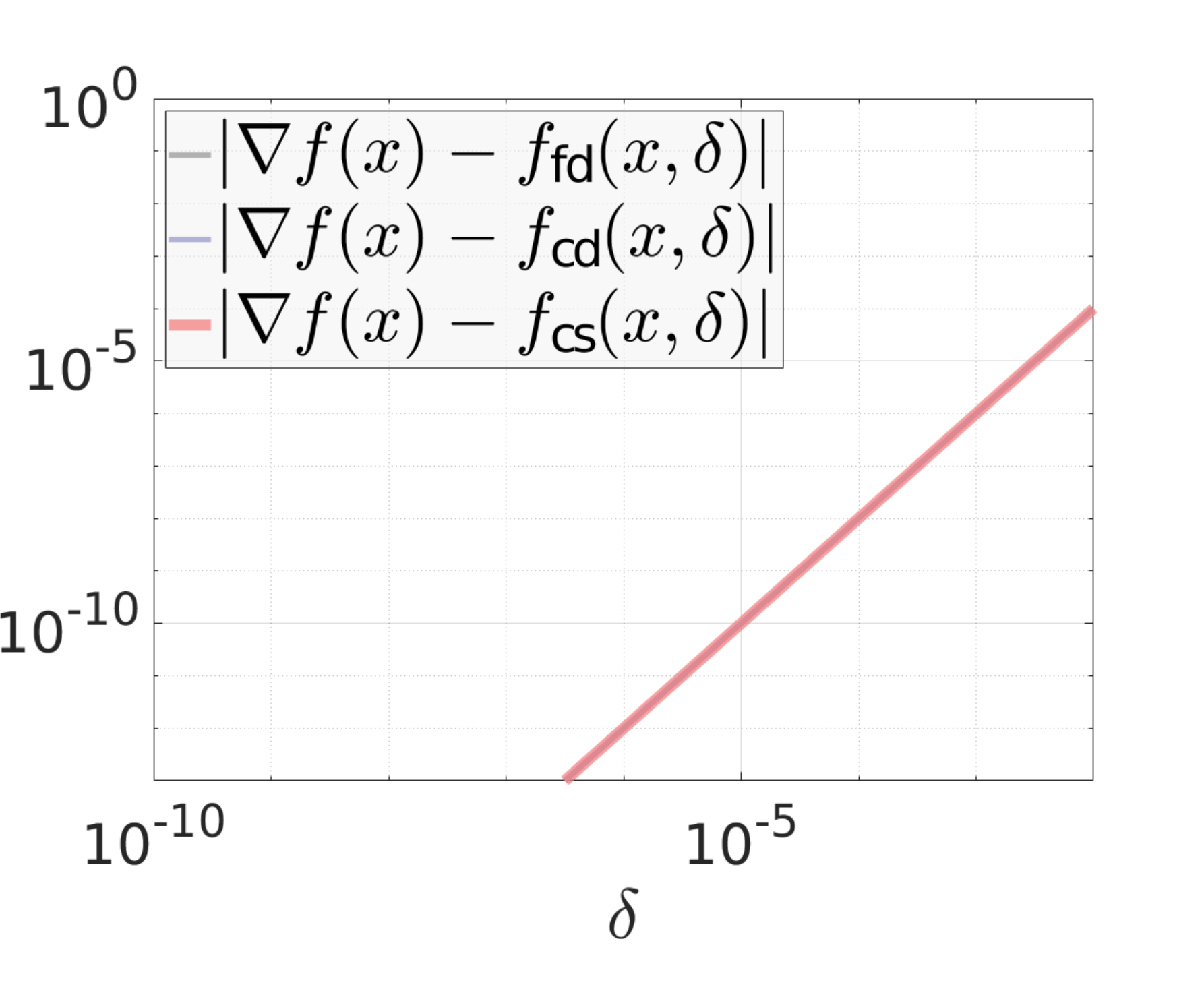}
        \caption{$x=0$}
        \label{fig:x00}
    \end{subfigure}\quad 
    \begin{subfigure}[b]{0.3\textwidth}
        \includegraphics[width=\textwidth]{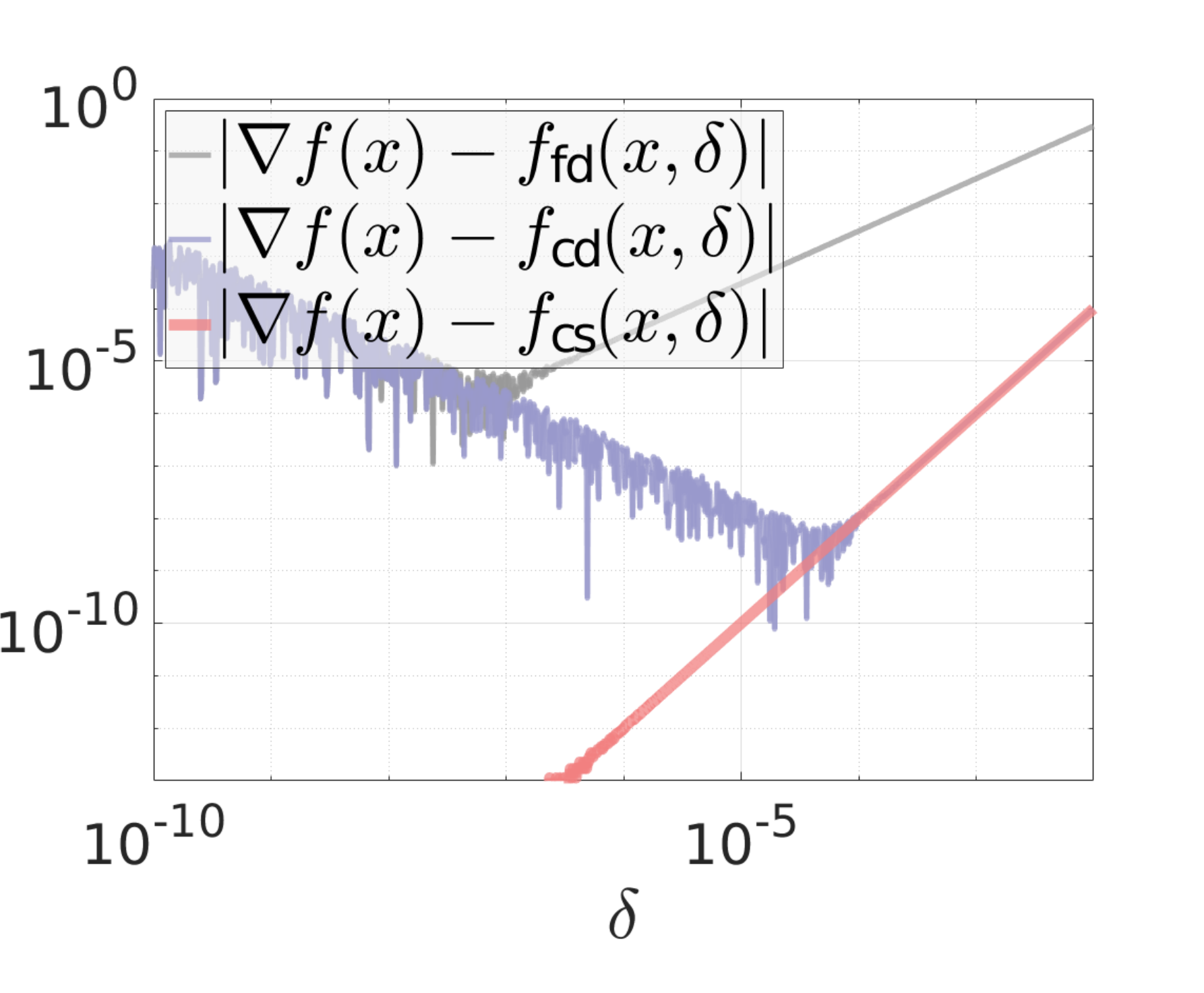}
        \caption{$x=10$}
        \label{fig:x010}
    \end{subfigure}
    \caption[]{Comparison of the gradient estimators of Example~\ref{ex:num_or} at different test points.}
    \label{fig:oracles}
\end{figure*}

\begin{example}[Numerical stability of gradient estimators]
\label{ex:num_or}
To showcase the power of the complex-step method and to expose the numerical difficulties encountered by finite-difference methods, we approximate the derivative of $f(x)=x^3$ at $x\in\{-1,0,10\}$ via a forward-difference~($\mathsf{fd}$), central-difference~($\mathsf{cd}$) and complex-step~($\mathsf{cs}$) method, that is, for small values of~$\delta$ we compare $f_{\mathsf{fd}}(x,\delta) = \tfrac{1}{\delta}({f(x+\delta)-f(x)})$, $f_{\mathsf{cd}}(x,\delta) = \tfrac{1}{2\delta}({f(x+\delta)-f(x-\delta)})$ and $f_{\mathsf{cs}}(x,\delta) = \tfrac{1}{\delta}{\Im( f(x+i\delta))}$. Figure~\ref{fig:oracles} visualizes the absolute approximation errors as a function of $\delta$. We observe that $f_{\mathsf{cd}}$ and $f_{\mathsf{cs}}$ offer the same approximation quality and incur an error of~$O(\delta^2)$ for all sufficiently large values of $\delta$. However, only the complex-step approximation reaches machine precision ($\approx 10^{-16}$), whereas both finite-difference methods deteriorate below $\delta\approx 10^{-6}$ due to subtractive cancellation errors. Note that for $x=0$ all errors are equal to~$\delta^2$ because~$f(0)=0$. As most existing zeroth-order optimization methods use finite-difference-based gradient estimators, we conclude that there is room for numerical improvements by leveraging complex arithmetic.
\end{example}

{
\subsection{Automatic differentiation}

    The complex-step approach is closely related to \textit{automatic differentiation} (AD)~\cite{ref:griewank2008evaluating,elliott2009beautiful}. AD decomposes the evaluation of $f$ into a partially ordered set of elementary operations and evaluates its derivative recursively using the rules of differentiation such as the chain and product rules etc. 
    Moreover, while the complex step approach evaluates $f$ at {\em complex} numbers of the form $a+ib$ with $a,b\in\mathbb R$ and an abstract imaginary unit $i$ satisfying $i^2=-1$, {a version of} forward-mode AD evaluates $f$ at {\em dual} numbers of the form $a+b\varepsilon$ with $a,b\in \mathbb{R}$ and $\varepsilon\neq 0$ an abstract number satisfying $\varepsilon^2=0$. The arithmetics of dual numbers imply that $f(x+\varepsilon)=f(x)+\varepsilon \partial_x f(x)$ whenever~$f$ is real analytic, and hence one can compute both $f(x)$ as well as $\partial_x f(x)$ in one forward pass. Intuitively, $\varepsilon$ should thus be interpreted as a nilpotent infinitesimal unit. While the set of complex numbers forms a {\em field}, the set of dual numbers only forms a {\em ring} (in fact, it forms the quotient ring $\mathbb{R}[\varepsilon]/\varepsilon^2$, which fails to be a field because multiplicative inverses and hence terms such as $\varepsilon^2/\varepsilon$ and $\sqrt{\varepsilon^2}$ are not defined).
    %
    The assumptions that $\varepsilon\neq 0$ and $\varepsilon^2=0$ require us to give up the law of the excluded middle \cite{ref:bell2008primer} and thus also the axiom of choice~\cite{ref:bauer2017five}.


The complex-step approach is computationally cheap, can {\em approximate} the derivative of an analytic function~$f$ at extremely high accuracy levels (see Figure~\ref{fig:oracles}) and even remains applicable when function evaluations are noisy~\cite{ref:more2014you}. AD is computationally more expensive than the complex-step approach~\cite{ref:martins2001connection}, {in settings where the latter applies}, yet, AD usually finds the {\em exact} derivative of~$f$. Whether or not AD will succeed, however, depends on the representation of~$f$; see, {\em e.g.}, \cite{ref:huckelheim2023understanding} for a discussion of possible pitfalls. That is, AD must be able to evaluate~$f$ at all dual numbers of the form~$x+\varepsilon$. The following example inspired by~\cite{ref:berz1992automatic} illustrates why this is restrictive. Consider the function $f\in C^{\omega}(\mathbb{R})$ defined through $f(x)=\mathrm{sinc}(x)=\sin(x)/x$ for all $x\in \mathbb{R}\setminus\{0\}$ with $f(0)=1$. This function is entire and has a unique global maximum of $1$ at $x=0$. Nevertheless, AD breaks down at $x=0$ because $1/\varepsilon$ is not defined. For instance, the deep learning toolbox in MATLAB as well as the state-of-the-art AD tools in Julia~\cite{Julia} (\textit{e.g.}, \texttt{ForwardDiff.jl}~\cite{ref:revels2016forward}, \texttt{Zygote.jl}~\cite{ref:innes2018don} and \texttt{Enzyme.jl}~\cite{ref:NEURIPS2020_9332c513}) or in Python (\textit{e.g.}, \texttt{JAX}~\cite{jax2018github}) evaluate $\partial_x f(0)$ to \texttt{NaN}, whereas the complex-step method provides a close approximation of the correct value $\partial_x f(0)=0$. We remark that the derivative of $f(x)=\mathrm{sinc}(x)$ is hard-coded in Julia.\footnote{\url{https://github.com/JuliaDiff/DiffRules.jl/blob/9030629bbea6b25851789af5f236f35c9009b1f6/src/rules.jl}} Hence, for AD to succeed the representation of~$f$ is critical, that is, $f$ must be defined as $f(x)=\mathrm{sinc}(x)$ instead of $f(x)=\sin(x)/x$. A simpler, yet contrived, example is $f(x)=x/x$, for which {the version of AD under consideration} evaluates $\partial_xf(0)$ to \texttt{NaN}, too. Note also that AD fails to compute the derivative of $f(x)=\mathrm{exp}((\sqrt{x^2})^2)$ at $0$ even though there is no division by~$0$. All of these problems emerge because the dual numbers form only a ring instead of a field. As pointed out in~\cite{ref:berz1992automatic}, these theoretical deficiencies of AD could be remedied by working with the \textit{Levi-Civita field}, whose members generalize the dual numbers and are representable as $\sum_{q\in \mathbb{Q}}a_q \varepsilon^q$ with $a_q\in \mathbb{R}$ for all $q\in \mathbb{Q}$. Unfortunately, the members of the Levi-Civita field do not admit a finite representation in general and are therefore difficult to handle computationally.
}

\subsection{Lipschitz inequalities}
\label{sec:lipschitz}
In order to be able to design reasonable zeroth-order optimization algorithms, we need to impose some regularity on the objective function~$f$. This is usually done by requiring~$f$ to display certain Lipschitz continuity properties. Following~\cite{nesterov2003introductory}, for any integers~$p,k\ge 0$ with~$p\leq k$, we thus use~$C^{k,p}_L(\mathcal{D})$ to denote the family of all $k$ times continuously differentiable functions on~$\mathcal D$ whose $p^{\mathrm{th}}$ derivative is Lipschitz continuous with Lipschitz constant~$L\geq 0$. Similarly, we use~$C^{\omega,p}_L(\mathcal{D})$ to denote the family of all analytic functions in~$C^{p,p}_L(\mathcal{D})$.

For example, if $f\in C^{1,1}_{{L_1}}(\mathcal{D})$, then $f$ has a Lipschitz continuous gradient, that is,
\begin{equation}
    \label{equ:grad:Lipschitz}
    \|\nabla f(x) - \nabla f(y) \|_2 \leq {L_1} \|x-y\|_2 \quad \forall x,y\in \mathcal{D}.
\end{equation}
By~\cite[Eq.~(6)]{ref:nesterov2017random}, this condition is equivalent to the inequality
\begin{equation}
    \label{equ:L1:upper}
    |f(y)-f(x) - \langle \nabla f(x),y-x \rangle | \leq \tfrac{1}{2}{L_1} \|x-y\|_2^2 \quad \forall x,y\in \mathcal{D}.  
\end{equation}
If $f\in C^{1,1}_{{L_1}}(\mathcal{D})$ is also convex then, the Lipschitz condition~\eqref{equ:grad:Lipschitz} is also equivalent to
\begin{equation}
    \label{equ:grad:Lipschitz:bound2}
    f(y)\geq f(x) + \langle \nabla f(x) , y-x  \rangle +\tfrac{1}{2L_1}\|\nabla f(y) - \nabla f(x) \|_2^2\quad \forall x,y \in \mathcal{D} ,
\end{equation}
see, \textit{e.g.,}~\cite{nesterov2003introductory}. In particular, if~$x$ is a local minimizer of~$f$ with $\nabla f(x)=0$, then the estimate~\eqref{equ:grad:Lipschitz:bound2} simplifies to $2{L_1}\left(f(y)-f(x)\right)\geq \|\nabla f(y)\|_2^2$ for all~$y\in\mathcal D$. 

If $f\in C^{2,2}_{{L_2}}(\mathcal{D})$, then $f$ has a Lipschitz continuous Hessian,~\textit{i.e.},
\begin{equation}
    \label{equ:Hess:Lip}
    \|\nabla^2 f(x)- \nabla^2 f(y)\|_2 \leq {L_2}\|x-y\|_2 \quad \forall x,y \in \mathcal{D}.
\end{equation}
By~\cite[Lem.~1.2.4]{nesterov2003introductory}, this condition is equivalent to the inequality
\begin{equation}
\label{equ:L2:bound}
    |f(y)-f(x) - \langle \nabla f(x), y-x\rangle - \tfrac{1}{2}\langle \nabla^2 f(x)(y-x), y-x \rangle |\leq \tfrac{1}{6}{L_2} \|x-y\|_2^3 \quad \forall x,y \in \mathcal{D}.
\end{equation}

More generally, any $f\in C^{p,p}_{{L_p}}(\mathcal{D})$ has a Lipschitz continuous $p^{\rm th}$ derivative. Recalling the definitions of higher-order partial derivatives and multi-indices, this requirement can be expressed~as
\begin{equation*}
	|\textstyle\sum_{|\alpha|=p} \partial^\alpha_x f(x) \cdot u^\alpha  - \sum_{|\alpha|=p} \partial^\alpha_x f(y) \cdot u^\alpha | \leq {L_p}\|x-y\|_2 \quad \forall x,y \in \mathcal{D},\;u\in \mathbb{S}^{n-1}.
\end{equation*}
It is often referred to as a $(p+1)^{\rm th}$-order smoothness condition~\cite[\S~1.1]{ref:bach:smooth} as it implies that any~$f\in C^{p+1,p}_{{L_p}}(\mathcal{D})\subseteq C^{p,p}_{{L_p}}(\mathcal{D})$ has a bounded $(p+1)^{\rm th}$ derivative, that~is, 
\begin{equation}
    \label{equ:pth:deriv:bound}
    | \textstyle\sum_{|\alpha|=p+1} \partial^\alpha_x f(x) \cdot u^\alpha | = |\partial^{p+1}_t f(x+tu)|_{t=0}|
    \leq {L_p}\quad \forall x\in \mathcal{D}, \; u\in \mathbb{S}^{n-1}.
\end{equation}

\section{A smoothed complex-step approximation}
\label{sec:oracle}
We now use ideas from~\cite{nemirovsky1983problem,ref:nesterov2017random} to construct a new gradient estimator, which can be viewed as a complex-step generalization of the estimators proposed in~\cite{nemirovsky1983problem,ref:Flaxman}. Our construction is based on the following assumption, which we assume to hold throughout the rest of the paper.
\begin{assumption}[Analytic extension]
\label{ass:holomorphic}
The function $f:\mathcal{D}\to \mathbb{R}$ of problem~\eqref{equ:opt:main} admits an analytic extension to the strip~$\mathcal{D}\times i\cdot (-\bar{\delta},\bar{\delta})^n$ for some $\bar{\delta}\in (0,1)$. 
\end{assumption}
Recall from Lemma~\ref{lem:analytic_extension} that~$f$ admits an analytic extension to some open set~$\Omega\subseteq \mathbb C^n$ covering~$\mathcal D$ whenever~$f\in C^\omega(\mathcal D)$. However, unless~$f$ is entire or~$\mathcal D$ is bounded, $\Omega$ may not contain a strip of the form envisaged in Assumption~\ref{ass:holomorphic}. Hence, this assumption is {\em not} automatically satisfied for any real analytic function~$f\in C^\omega(\mathcal D)$. The requirement~$\bar{\delta}\in (0,1)$ is unrestrictive and has the convenient consequence that~$\delta^{p}\leq \delta^{p-1}$ for any~$\delta \in(0, \bar\delta)$ and~$p\in\mathbb{Z}_{\geq 0}$. 
All subsequent results are based on a smoothed complex-step approximation~$f_\delta$ of~$f$, which is defined through
\begin{equation}
\label{equ:f:delta}
	 f_{\delta}(x) = \textstyle{V_n^{-1}}\int_{\mathbb{B}^n}\Re\big(f(x+i\delta y)\big) \mathrm{d}y. 
\end{equation}
Here, the radius~$\delta \in (0,\bar\delta)$ of the ball used for averaging represents a tuneable smoothing parameter. Given prior structural knowledge about~$f$, one could replace~$ \mathbb{B}^n$ with a different compact set~\cite{ref:hazan2014bandit,ref:JongeneelIZO2021}. We emphasize that the integral in~\eqref{equ:f:delta} is well-defined whenever~$\delta \in (0,\bar\delta)$, which ensures that~$f$ has no singularities in the integration domain. Next, we address the approximation quality of~$f_\delta$.

\begin{proposition}[Approximation quality of $f_\delta$]
\label{prop:approx:function}
If $f\in C^{\omega,1}_{{L_1}}(\mathcal{D})$ satisfies Assumption~\ref{ass:holomorphic}, then for $f_{\delta}$ defined as in~\eqref{equ:f:delta} and for any fixed~$x\in \mathcal{D}$ and~$\kappa\in(0,1)$ there exists~$C_\kappa\ge 0$ with 
\begin{align*}
    \left|f_\delta (x) - f(x) \right| &\le\tfrac{1}{2} L_1 \delta^2 +C_\kappa \delta^4\quad \forall \delta \in(0,\kappa\bar\delta].
\end{align*}
\end{proposition}

\begin{proof}
By the definition of~$f_\delta$ in~\eqref{equ:f:delta}, we have 
$| f_\delta (x)-f(x)|  \leq  \textstyle{V_n^{-1}}\int_{\mathbb{B}^n} | \Re\big(f(x+i\delta y)\big) - f(x) | \mathrm{d}y$. The Taylor series of~$f(x+i\delta y)$ around~$x$ then yields
\begin{align*}
	\Re\left(f(x+i\delta y)\right)- f(x) &=  \textstyle\sum^{\infty}_{k=0} \frac{(-1)^{k} \delta^{2k}}{(2k)!} \sum_{|\alpha|=2k} \partial^\alpha_x f(x) y^\alpha - f(x) \\
	&=  -\tfrac{1}{2} \delta^{2}\langle \nabla^2 f(x) y, y\rangle + \delta^4R(y,\delta),
\end{align*}
where the real-valued remainder term~$R(y,\delta)$ is continuous in~$y\in\mathbb B^n$ and~$\delta\in [0,\bar{\delta})$. Substituting the last expression into the above estimate and using~\eqref{equ:pth:deriv:bound}, we obtain
\begin{align*}
    \left| f_\delta (x)-f(x)\right| & \leq \textstyle{V_n^{-1}}\int_{\mathbb{B}^n}\tfrac{1}{2}\delta^2  L_1 + \delta^4 |R(y,\delta)| \mathrm d y 
    \le \tfrac{1}{2} \delta^2 {L_1} + C_\kappa \delta^4 \quad \forall \delta\in(0,\kappa\bar \delta],
\end{align*}
where the non-negative constant
$	C_\kappa = \max_{y\in \mathbb{B}^{n}}\max_{\delta\in[0,\kappa \bar\delta]} |R(y,\delta)| $
is finite due to continuity of~$R(y,\delta)$ and compactness of~$\mathbb{B}^{n}$ and~$[0,\kappa \bar\delta]$. Hence, the claim follows.
\end{proof}

Note that if $f$ is affine, then~ $f_{\delta}=f$. Note also that $[0,\kappa\bar \delta]$ is a compact subset of the set~$[0,\bar \delta)$ on which~$R(y,\delta)$ is continuous in~$\delta$ and that~$R(y,\delta)$ may be unbounded on~$[0,\bar\delta)$.  
The following proposition provides an integral representation for the gradient of~$f_{\delta}$. It extends~\cite[\S~9.3]{nemirovsky1983problem} and~\cite[Lem.~1]{ref:Flaxman} to the realm of complex arithmetic.

\begin{proposition}[Gradient of the smoothed complex-step function]
\label{prop:CR:grad}
If $f\in C^{\omega}(\mathcal{D})$ satisfies Assumption~\ref{ass:holomorphic}, then~$f_{\delta}$ defined as in~\eqref{equ:f:delta} is differentiable, and we have
\begin{equation}
\label{equ:unif:grad}
    \nabla f_{\delta}(x) = \textstyle\frac{n}{\delta}\mathbb{E}_{y\sim \sigma}\left[\Im\left(f(x+i\delta y)\right)y \right]\quad \forall x\in \mathcal{D}, \; \delta\in (0,\bar{\delta}),
\end{equation}
where $\sigma$ denotes the uniform distribution on $\mathbb{S}^{n-1}$.
\end{proposition}
\begin{proof}
Any function~$g\in C^1(\mathbb R^n)$ and vector~$w\in\mathbb R^n$ define a vector field~$v(y)=g(y)\cdot w$. 
The divergence theorem~\cite[Thm.~16.32]{Lee2} then implies that
\begin{align*}
     \textstyle\int_{\mathbb{B}^n} \langle w, \nabla g(y)\rangle \, \mathrm d y = \int_{\mathbb{B}^n}\mathrm{div}(v(y))\, \mathrm d y &= S_{n-1} \textstyle\int_{\mathbb{S}^{n-1}}\langle v(y),y \rangle \,\sigma(\mathrm d y)\\
     &=S_{n-1} \textstyle\int_{\mathbb{S}^{n-1}} g(y) \langle w,y \rangle \,\sigma(\mathrm d y),
\end{align*}
where the scaling factor~$S_{n-1}$ accounts for the fact that the uniform distribution~$\sigma$ is normalized on~$\mathbb{S}^{n-1}$. Note also that the outward-pointing unit normal vector of~$\mathbb{S}^{n-1}$ at any point~$y\in \mathbb{S}^{n-1}$ is exactly~$y$ itself. As the above equation holds for all vectors~$w\in\mathbb R^n$ and as both the leftmost and rightmost expressions are linear in~$w$, their gradients must coincide. This reasoning implies that
\begin{equation}
\label{equ:div2}
    \textstyle\int_{\mathbb{B}^n}\nabla g(y)\, \mathrm dy = S_{n-1} \int_{\mathbb{S}^{n-1}}g(y) y\, \sigma(\mathrm dy). 
\end{equation}
We are now ready to prove~\eqref{equ:unif:grad} by generalizing tools developed in~\cite{nemirovsky1983problem,ref:Flaxman} to the complex domain. Specifically, by the definition of~$f_\delta$ in~\eqref{equ:f:delta} we have
\begin{align*}
     \nabla f_{\delta}(x) 
     =\textstyle  V_n^{-1}\int_{\mathbb{B}^n}\nabla_x\Re\left(f(x+i\delta y)\right)\mathrm{d} y&= \textstyle  (V_n \delta)^{-1}\int_{\mathbb{B}^n}\nabla_y\Im\left(f(x+i\delta  y)\right)\mathrm{d} y\\
     &= \textstyle{S_{n-1}}{(V_n\delta)^{-1}}\int_{\mathbb{S}^{n-1}}\Im\left(f(x+i\delta y)\right)y\, \sigma(\mathrm{d}y)\\ &=\textstyle{S_{n-1}}{(V_n\delta)^{-1}} \mathbb{E}_{y\sim \sigma}\left[\Im\left(f(x+i\delta y)\right)y \right], 
\end{align*}
where the interchange of the gradient and the integral in the first equality is permitted by the dominated convergence theorem, which applies because~$\mathbb B^n$ is compact and because any continuously differentiable function on a compact set is Lipschitz continuous. The second equality is a direct consequence of the Cauchy-Riemann equations, and the third equality, finally, holds thanks to the generalized Achimedean principle~\eqref{equ:div2} with pressure function~$g(y)=\Im\left(f(x+i\delta  y)\right)$. We finally observe that the volume of the unit ball and the surface of the unit sphere satisfy $\textstyle	V_n=\int_{\mathbb B^n} \mathrm{d}y = S_{n-1}\int_0^rr^{n-1}\mathrm{d}r= {S_{n-1}}/{n}$ $\implies$ ${S_{n-1}}/{V_n}= n$. Thus, the claim follows.
\end{proof}

Proposition~\ref{prop:CR:grad} reveals that~$\nabla f_\delta$ admits the unbiased single-point estimator 
\begin{equation}
    \label{equ:grad:est:g}
    g_{\delta}(x) = \tfrac{n}{\delta}\Im\left(f(x + i \delta y)\right)y \quad \text{with}\quad y\sim \sigma.
\end{equation}
Now we show that~$\nabla f_\delta(x)$ approximates~$\nabla f(x)$ arbitrarily well as~$\delta$ drops to~$0$. 

\begin{proposition}[Approximation quality of $\nabla f_\delta$]
\label{prop:grad:approx}
If $f\in  C^{\omega,2}_{{L_2}}(\mathcal{D})$ satisfies Assumption~\ref{ass:holomorphic}, then for $f_{\delta}$ defined as in~\eqref{equ:f:delta} and for any fixed~$x\in \mathcal{D}$ and~$\kappa\in(0,1)$ there exists~$C_\kappa\ge 0$ with 
\begin{align}
\label{equ:grad:delta:approx:error}
    \|\nabla f_{\delta}(x) - \nabla f(x) \|_2 &\leq \tfrac{1}{6} n L_2 \delta^2 + n C_\kappa\delta^4 \quad \forall \delta \in(0,\kappa\bar\delta].
\end{align}
\end{proposition}

\begin{proof}
If we denote as usual by~$I_n$ the identity matrix in~$\mathbb R^n$, then the covariance matrix of the uniform distribution~$\sigma$ on the unit sphere~$\mathbb S^{n-1}$ can be expressed as
\begin{equation}
\label{equ:int:sphere}
	\textstyle\int_{\mathbb S^{n-1}} yy^{\mathsf{T}} \sigma(\mathrm d y)  = \int_{\mathbb S^{n-1}} \|y\|_2^2\,\sigma(\mathrm d y) \cdot\tfrac{1}{n} I_n  = \tfrac{1}{n} I_n ,
\end{equation}
where the two equalities hold because the sought covariance matrix must be isotropic and because~$\|y\|_2=1$ for all~$y\in\mathbb S^{n-1}$, respectively. Thus, the gradient of~$f$ can be represented as $\nabla f(x) = n \int_{\mathbb{S}^{n-1}} \langle \nabla f(x), y \rangle y \,\sigma( \mathrm{d}y)$. Together with Proposition~\ref{prop:CR:grad}, this yields the estimate
\begin{align*}
    \|\nabla f_{\delta}(x) - \nabla f(x) \|_2 &= \textstyle\frac{n}{\delta} \textstyle\left\| \int_{\mathbb{S}^{n-1}} \Im\left(f(x+i\delta y)\right) y -\delta \langle \nabla f(x), y \rangle  y\,\sigma(\mathrm d y) \right\|_2 \\
    &\leq \textstyle \frac{n}{\delta}\textstyle\int_{\mathbb{S}^{n-1}} \left|\Im\left(f(x+i\delta y)\right)-\delta \langle \nabla f(x), y \rangle  \right| \|y\|_2\,\sigma(\mathrm d y).
\end{align*}
By using the Taylor series representation of~$f(x+i\delta y)$ around~$x$, we find
\begin{align*}
	\Im\left(f(x+i\delta y)\right)-\delta \langle \nabla f(x), y \rangle &=  \textstyle\sum^{\infty}_{k=0} \frac{(-1)^k \delta^{2k+1}}{(2k+1)!} \sum_{|\alpha|=2k+1} \partial^\alpha_x f(x) y^\alpha -\delta \langle \nabla f(x), y \rangle \\
	&=  \textstyle\sum^{\infty}_{k=1} \frac{(-1)^k \delta^{2k+1}}{(2k+1)!} \sum_{|\alpha|=2k+1} \partial^\alpha_x f(x) y^\alpha \\
	&=  -\tfrac{1}{6} \delta^{3}\textstyle\sum_{|\alpha|=3} \partial^\alpha_x f(x) y^\alpha + \delta^5 R(y,\delta),
\end{align*}
where the real-valued remainder term~$R(y,\delta)$ is continuous in~$y\in\mathbb B^n$ and~$\delta\in [0,\bar{\delta})$. Substituting the last expression into the above and using~\eqref{equ:pth:deriv:bound}, we obtain
\begin{align*}
    \|\nabla f_{\delta}(x) - \nabla f(x) \|_2 \leq & \tfrac{n}{\delta}\textstyle \int_{\mathbb{S}^{n-1}} \left(\tfrac{1}{6}{\delta^3 {L_2}} + \delta^5\,|R(y,\delta)| \right)\|y\|_2 \,\sigma(\mathrm{d}y)\\ 
    \leq& \tfrac{1}{6} \delta^2 n {L_2} + n C_\kappa \delta^4 \quad \forall \delta\in(0,\kappa\delta],
\end{align*}
where the non-negative constant $C_\kappa = \max_{y\in \mathbb{S}^{n-1}}\max_{\delta\in[0,\kappa \bar\delta]} |R(y,\delta)| $ is again finite due to the continuity of~$R(y,\delta)$ and the compactness of~$\mathbb{S}^{n-1}$ and~$[0,\kappa \bar\delta]$. \end{proof}

Proposition~\ref{prop:grad:approx} implies that the \textit{single-point} estimator~\eqref{equ:grad:est:g} incurs only errors of the order~$O(\delta^2)$ on average. Equally small errors were attained in~\cite{ref:nesterov2017random} for~$f\in C^{2,2}_{{L_2}}$ by using Gaussian smoothing and a \textit{multi-point} estimator. Unfortunately, the latter is susceptible to cancellation effects. Proposition~\ref{prop:grad:approx} also implies that $\lim_{\delta\downarrow 0}\nabla f_{\delta} (x) = \nabla f(x)$. In addition, one readily verifies that if~$f$ is quadratic (that is, if~$L_2=0$), then~$\nabla f_{\delta} (x) = \nabla f(x)$ for all~$x\in\mathcal D$ and~$\delta\in(0, \bar \delta)$. 
The single-point estimator~$g_{\delta}(x)$ introduced in~\eqref{equ:grad:est:g} is unbiased by construction. In addition, as for the multi-point estimator proposed in~\cite{ref:nesterov2017random}, the second moment of~$g_{\delta}(x)$ admits a convenient bound.

\begin{corollary}[Second moment of~$g_{\delta}(x)$]
\label{cor:im:lip}
If $f\in  C^{\omega,2}_{{L_2}}(\mathcal{D})$ satisfies Assumption~\ref{ass:holomorphic}, then for $g_{\delta}$ as in~\eqref{equ:grad:est:g} and for any fixed~$x\in \mathcal{D}$ and~$\kappa\in(0,1)$ we have 
\begin{equation}
    \label{equ:g:var}
    \begin{aligned}
    \mathbb{E}_{y\sim \sigma}\left[ \|g_{\delta}(x)\|_2^2\right] &\leq n^2(\tfrac{1}{6}L_2\delta^2 + C_{\kappa}\delta^4)^2+n\|\nabla f(x)\|_2^2\\ & \hspace{1.5cm} + 2n^2\left(\tfrac{1}{6}L_2\delta^2 + C_{\kappa}\delta^4\right) \|\nabla f(x)\|_2,
\end{aligned}
\end{equation}
where~$C_\kappa\geq 0$ is the same constant as in Proposition~\ref{prop:grad:approx}.
\end{corollary}

\begin{proof}
Using the definition of~$g_\delta$ and the fact that~$\|y\|_2=1$~$\forall\,y\in\mathbb{S}^{n-1}$, we find
\begin{equation}
	\label{equ:int:moment}
	\textstyle\mathbb{E}_{y\sim \sigma}\left[ \|g_{\delta}(x)\|_2^2\right] = \frac{n^2}{\delta^2}\mathbb{E}_{y\sim \sigma}[\left(\Im\left(f(x+i\delta y) \right)\right)^2 ].
\end{equation}
By essentially the same arguments as in the proof of Proposition~\ref{prop:grad:approx}, we further have
\begin{equation*}
    \begin{aligned}
    \left|\Im\left(f(x+i\delta y) \right)\right| &= \left|\Im\left(f(x+i\delta y) \right) - \langle \nabla f(x), \delta y\rangle + \langle \nabla f(x), \delta y \rangle\right|\\
    &\leq  \left|\tfrac{1}{6}\delta^3 {L_2} + \delta^5 C_{\kappa}\right| + \left|\langle \nabla f(x), \delta y\rangle\right|. 
\end{aligned}
\end{equation*}
Squaring the above and applying the Cauchy-Schwarz inequality yields
\begin{equation*}
    \begin{aligned}
    \left|\Im\left(f(x+i\delta y) \right)\right|^2 \leq  \left(\tfrac{1}{6}\delta^3 {L_2} + \delta^5 C_{\kappa}\right)^2 &+ \langle \nabla f(x), \delta y\rangle^2\\
    &+2\delta \left(\tfrac{1}{6}\delta^3L_2+\delta^5C_{\kappa}\right)\|\nabla f(x)\|_2\|y\|_2.
\end{aligned}
\end{equation*}
The claim then follows from substituting the above into~\eqref{equ:int:moment} and using~\eqref{equ:int:sphere}.
\end{proof}
In analogy to Proposition~\ref{prop:grad:approx}, one readily verifies that if~$f$ is quadratic (\textit{i.e.}, if~$L_2=0$), then the right hand side of~\eqref{equ:g:var} vanishes. Under a third-order smoothness condition, there exist multi-point estimators that satisfy a bound akin to~\eqref{equ:g:var}~\cite[Thm.~4.3]{ref:nesterov2017random}.


Unlike the smooth approximations proposed in~\cite{ref:nesterov2017random}, the smoothed complex-step approximation~$f_{\delta}$ does frequently {\em not} belong to the same function class as~$f$. For example, even though the Lorentzian function~$f(x)=1/(1+x^2)$ has a Lipschitz continuous gradient with ${L_1}=2$, the Lipschitz modulus of its approximation~$f_\delta$ strictly exceeds~$2$ for some values of~$\delta$ close to~$1$ because~$f$ has two poles at~$i$ and~$-i$. Similarly, $f_\delta$ does not necessarily inherit convexity from~$f$.
\begin{example}[Loss of convexity]
\label{ex:ncvx:f}
If $f\in C^\omega(\mathbb R)$ is entire, then it has a globally convergent power series representation with real coefficients. Consequently, $f$ satisfies
\begin{equation*}
    \Re(f(x+i\delta y)) = \textstyle\sum^{\infty}_{k=0}(-1)^k \frac{f^{(2k)}(x)}{(2k)!}(\delta y)^{2k}. 
\end{equation*}
In the special case when $f(x)=x^2$, the complex-step approximation $\Re(f(x+i\delta y))=x^2-(\delta y)^2$ inherits convexity from~$f$ regardless of the choice of~$\delta>0$ and $y\in \mathbb{R}$. Thus, $f_\delta$ is also convex because convexity is preserved by integration. However, if~$f(x) = x^4$, then we find~$\Re(f(x+i \delta y))=x^4-6x^2 (\delta y)^2 + (\delta y)^4$, which fails to be convex in~$x$ for any~$\delta>0$ and~$y\neq 0$. 
In this case,~$f_\delta$ remains non-convex despite the smoothing. {Finally, if~$f$ is strongly convex ({\em e.g.}, if $f(x)=x^2+x^4$), then one readily verifies that $\Re(f(x+i\delta y))$ is convex in~$x$ provided that~$\delta$ is sufficiently small.
}
\end{example}

If $f_{\delta}$ inherited convexity from~$f$, one could simply incorporate the estimator~\eqref{equ:grad:est:g} into the algorithms studied in~\cite[\S~5]{ref:nesterov2017random}, and the corresponding convergence analysis would carry over with minor modifications. As the smoothed complex-step approximation may destroy convexity, however, a different machinery is needed here.

\section{Convex optimization}
\label{sec:convex}
We now study the convergence properties of zeroth-order algorithms for solving problem~\eqref{equ:opt:main} under the assumption that~$f$ is a convex function on~$\mathcal D$ and~$\mathcal X$ is a non-empty closed convex subset of~$\mathcal D$. Our methods mimic existing algorithms developed in~\cite{ref:nesterov2017random} but use the single-point estimator~$g_{\delta}$ defined in~\eqref{equ:grad:est:g} instead of a multi-point estimator that may suffer from cancellation effects. Our method is described in Algorithm~\ref{alg:convex_unconstrained}, where~$\Pi_{\mathcal{X}}:\mathcal{D}\to \mathcal{X}$ denotes the Euclidean projection onto~$\mathcal X$. Note that~$\Pi_{\mathcal{X}}$ reduces to the identity operator if~$\mathcal X=\mathcal D$. 

In the remainder we will assume that the iterates~$\{x_k\}_{k\in\mathbb{Z}_{>0}}$ generated by Algorithm~\ref{alg:convex_unconstrained} as well as all samples~$\{y_k\}_{k\in\mathbb{Z}_{>0}}$ and the corresponding gradient estimators~$\{g_{\delta_k}(x_k)\}_{k\in\mathbb{Z}_{>0}}$ represent random objects on an abstract filtered probability space~$(\Omega,\mathcal F, \{\mathcal F_k\}_{k\in\mathbb{Z}_{>0}},\mathbb P)$, where~$\mathcal F_k$ denotes the $\sigma$-algebra generated by the independent and identically distributed samples~$y_1,\ldots,y_{k-1}$. Therefore, $x_k$ is $\mathcal F_k$-measurable. In the following, we use~$\mathbb E[\cdot]$ to denote the expectation operator with respect to~$\mathbb P$. 

\begin{algorithm}[t]
\begin{algorithmic}[1]
\STATE \textbf{Input:} initial iterate $x_1\in\mathcal{X}$, stepsizes $\{\mu_k \}_{k\in\mathbb{Z}_{\geq 0}}$, smoothing parameters $\{\delta_k \}_{k\in\mathbb{Z}_{\geq 0}}$
  \FOR{$k = 1,2,\ldots, K-1$} 
  \STATE sample $y_k \sim \sigma$
  
  \STATE set $g_{\delta_k}(x_k) = \frac{n}{\delta_k}\Im\left(f(x_k + i \delta_k y_k)\right)y_k$
  \STATE set $x_{k+1} = \Pi_{\mathcal{X}}\left( x_k - \mu_k\, g_{\delta_k}(x_k)\right)$
\ENDFOR
\STATE \textbf{Output:} last iterate $x_K$ and averaged iterate $\Bar{x}_{K}=\frac{1}{K}\sum^{K}_{k=1}x_k$
\end{algorithmic}
\caption{Imaginary zeroth-order optimization} 
\label{alg:convex_unconstrained}
\end{algorithm}

\begin{theorem}[Convergence rate of Algorithm~\ref{alg:convex_unconstrained} for convex optimization]
\label{thm:unconstrained}
Suppose that~$f$ is a convex, real analytic function satisfying Assumption~\ref{ass:holomorphic} as well as the Lipschitz conditions~\eqref{equ:grad:Lipschitz} and~\eqref{equ:Hess:Lip} with ${L_1}>0$ and ${L_2}\geq 0$. Also assume that $\mathcal{X}$ is non-empty, closed and convex and that there exists $x^{\star}\in \mathcal{X}$ with $\nabla f(x^{\star})=0$. Denote by $\{x_k\}_{k\in \mathbb{Z}_{>0}}$ the iterates generated by Algorithm~\ref{alg:convex_unconstrained} with constant stepsize $\mu_k=\mu=1/(2nL_1)$ and adaptive smoothing parameter $\delta_k\in (0,\kappa \bar{\delta}]$ for all $k\in \mathbb{Z}_{>0}$, where $\kappa \in (0,1)$, and define $R=\|x_1-x^{\star}\|_2$. Then, the following hold for all $K\in \mathbb{Z}_{>0}$. 
\begin{enumerate}[(i)]
    \item There is a constant $C_1\geq 0$ such that 
   \begin{align*}
    \mathbb{E} \left[ f(\bar{x}_{K}) - f(x^\star) \right] \leq&  \tfrac{1}{\mu K} R^2 + \tfrac{1}{K}{C_1 n }R \textstyle\sum_{k=1}^{K} \delta_k^2 + \tfrac{1}{K}{\mu C_1^2 n^2  } (\sum_{k=1}^{K} \delta_k^2)^2 \\
        & + \tfrac{1}{K}{\mu C_1C_2 n^2  } (\textstyle\sum_{k=1}^{K} \delta_k^2) (\textstyle\sum_{k=1}^{K} \delta_k^4)^{\frac{1}{2}} + \tfrac{1}{K}{\mu C_2^2 n^2 } \textstyle\sum_{k=1}^{K}  \delta_k^4.
    \end{align*}    
\item \noindent If $\delta_k=\delta $ for all $k\in \mathbb{Z}_{>0}$, then we have
\begin{equation*}
    \mathbb{E} \left[ f(\bar{x}_{K}) - f(x^\star) \right] \le \textstyle\frac{1}{K}{2n{L_1}} R^2 + {C_1 n} R\delta^2 + \tfrac{1}{{L_1}} (1+\sqrt{K})^2 C_1^2n\delta^4. 
\end{equation*}
\item If $\delta_k = \delta/k$ for all $k\in \mathbb{Z}_{>0}$, then there is a constant $C_2\geq 0$ such that
\begin{equation*}
   \mathbb{E} \left[ f(\bar{x}_{K}) - f(x^\star) \right] \le \textstyle\frac{n}{K} \left( \sqrt{2 {L_1}} R + C_2 \delta^2 \right)^2.
\end{equation*}
\end{enumerate}
\end{theorem}

Under the assumptions of Theorem~\ref{thm:unconstrained}, problem~\eqref{equ:opt:main} is convex and $x^{\star}$ represents a global minimizer. Note, however, that~$\mathcal{X}$ may not contain any~$x^{\star}$ with~$\nabla f(x^\star)=0$ even if~$\mathcal X$ is compact. This is usually the case if the global minimum of~\eqref{equ:opt:main} is attained at the boundary of~$\mathcal X$. 
{If~$x^\star$ is not unique, one should set $R=\|x_1-P^{\star}(x_1)\|_2$ for the bounds not to be trivial, with $P^{\star}(x_1)=\argmin_{x^{\star}}\|x^{\star}-x_1\|_2^2$, which is well-defined since $f$ is convex, real analytic.} 
Explicit formulas for~$C_1$ and $C_2$ in terms of $\kappa$, $L_2$ etc.\ are derived in the proof of Theorem~\ref{thm:unconstrained}. 

\begin{proof}[Proof of Theorem~\ref{thm:unconstrained}]
For ease of notation, we define~$r_{k} =\|x_k-x^{\star}\|_2$ for all~$k \in\mathbb{Z}_{> 0}$. We prove the theorem first under the simplifying assumption that~$\mathcal{X}=\mathcal{D}$, 
which implies the projection onto~$\mathcal X$ becomes obsolete, that is, $x_{k+1}=x_k - \mu_k\cdot g_{\delta_k}(x_k)$. Thus, we have
\begin{align*}
     \mathbb{E} \left[ \left. r_{k+1}^2\,\right| \mathcal F_k \right] & =  \mathbb{E} \left[ \left. r_k^2 -2\mu_k \langle g_{\delta_k}(x_k),x_k-x^{\star}\rangle + \mu_k^{2} \, \|g_{\delta_k}(x_k)\|_2^2\, \right| \mathcal F_k \right] \\
    & =   r_k^2 - 2\mu_k \langle \nabla f_{\delta_k} (x_k),x_k-x^{\star}\rangle + \mu_k^{2} \, \mathbb{E} \left[ \left. \|g_{\delta_k}(x_k)\|_2^2 \,\right| \mathcal F_k \right],
\end{align*}
where the second equality follows from~\eqref{equ:unif:grad}, the definition of~$g_{\delta_k}(x_k)$ and the $\mathcal F_k$-measurability of~$x_k$ and~$r_k$. The Cauchy-Schwartz inequality then implies that
\begin{align*}
    &\mathbb{E} \left[ \left. r_{k+1}^2\,\right| \mathcal F_k \right]\\ & \le  r_k^2 - 2\mu_k \langle \nabla f (x_k),x_k-x^{\star}\rangle + 2\mu_k \| \nabla f_{\delta_k} (x_k) - \nabla f (x_k) \|_2\, r_k \\ & \qquad
     + \mu_k^{2} \, \mathbb{E} \left[ \left. \|g_{\delta_k}(x_k)\|_2^2\,\right| \mathcal F_k \right] \\
    & \le  r_k^2 - 2\mu_k \left( f(x_k) - f(x^\star) \right) + 
    2\mu_k\left(\tfrac{1}{6} n L_2 \delta_k^2 + n C_\kappa \delta_k^4 \right) \,r_k \\
    &  \qquad + \mu_k^{2}n^2 \left(\left(\tfrac{1}{6}L_2\delta_k^2 + C_\kappa \delta_k^4\right)^2 + \tfrac{1}{n}\|\nabla f(x_k)\|_2^2+2 \left(\tfrac{1}{6}L_2\delta_k^2 + C_\kappa\delta_k^4\right)\|\nabla f(x_k)\|_2 \right)  \\
    & \le r_k^2 - 2\mu_k \left( f(x_k) - f(x^\star) \right)  + 
    2n\mu_k\delta_k^2 \left(\tfrac{1}{6} L_2 + C_\kappa \delta_k^2 +nL_1\mu_k\left(\tfrac{1}{6}L_2 + C_\kappa\delta_k^2\right)\right) r_k \\
    &\qquad + \mu_k^{2}n^2 \left( \delta_k^4 \left(\tfrac{1}{6}L_2 + C_\kappa\delta_k^2\right)^2 + \tfrac{1}{n}{2 L_1}(f(x_k)-f(x^{\star})) \right), 
\end{align*}
where the second inequality exploits the convexity of~$f$ as well as Proposition~\ref{prop:grad:approx} and Corollary~\ref{cor:im:lip}, while the third inequality follows from the estimates~\eqref{equ:grad:Lipschitz} and~\eqref{equ:grad:Lipschitz:bound2}, which imply that $\|\nabla f(x_k)\|_2\leq L_1\|x_k-x^\star\|_2$ and $2L_1(f(x_k)-f(x^\star))\geq \|\nabla f(x_k)\|_2^2$, respectively. To simplify notation, we now introduce the constant~$C_1=\tfrac{1}{2}L_2 + 3C_\kappa$, which upper bounds~$\tfrac{1}{2}L_2 + 3C_\kappa\delta_k^2$ and~$\tfrac{1}{6}L_2 + C_\kappa\delta_k^2$ for any~$k\in\mathbb{Z}_{> 0}$ because all smoothing parameters belong to the interval~$[-1,1]$. Recalling that the stepsize is constant and equal to~$\mu=1/(2nL_1)$, the above display equation thus simplifies to
\begin{equation}
\label{eq:r2-estimate}
    \mathbb{E} \left[ \left. r_{k+1}^2\,\right| \mathcal F_k \right]  \le r_k^2 - \mu \left( f(x_k) - f(x^\star) \right) + n\mu \delta_k^2 C_1 r_k  + \mu^{2}n^2 C_1^2 \delta_k^4 .
\end{equation}
Taking unconditional expectations and rearranging terms then yields
\begin{align*}
    \mathbb{E} \left[ f(x_k) - f(x^\star) \right] & \le \tfrac{1}{\mu} \left(\mathbb{E} \left[  r_k^2 \right] - \mathbb{E} \left[ r_{k+1}^2  \right] \right) + n{C_1\delta_k^2  } \mathbb{E}\left[ r_k \right] + {\mu n^2 C_1^2  \delta_k^4} \\
    &\le \tfrac{1}{\mu} \left(\mathbb{E} \left[  r_k^2 \right] - \mathbb{E} \left[ r_{k+1}^2  \right] \right) + n{C_1 \delta_k^2  } \textstyle\sqrt{\mathbb{E}\left[ r_k^2 \right] }+ \mu n^2 C_1^2  \delta_k^4.
\end{align*}
Next, choose any~$k'\in\mathbb{Z}_{> 0}$ and sum the above inequalities over all $k\leq k'-1$ to obtain
\begin{equation}
   \label{ineq:f_k_suboptimality_sum}
   \textstyle \sum_{k = 1}^{k'-1} \mathbb{E} \left[f(x_k) - f(x^\star)  \right] \leq  \tfrac{1}{\mu} \left( r_1^2 - \mathbb{E} \left[ r_{k'}^2  \right] \right) + {C_1 n } \sum_{k=1}^{k'-1} \delta_k^2 \sqrt{\mathbb{E}\left[ r_k^2 \right] }+ {\mu C_1^2 n^2 } \sum_{k=1}^{k'-1}  \delta_k^4 .
\end{equation}
Clearly, the inequality~\eqref{ineq:f_k_suboptimality_sum} remains valid if we lower bound its left hand side by~$0$ and upper bound its right hand side by increasing the upper limits of the two sums to~$k'$. We then obtain $\textstyle \mathbb{E}[ r_{k'}^2 ] \le r_1^2 +  {\mu C_1 n } \sum_{k=1}^{k'} \delta_k^2 \sqrt{\mathbb{E}\left[ r_k^2 \right] } +  {\mu^2 C_1^2 n^2 } \sum_{k=1}^{k'}  \delta_k^4 $. Setting~$t_k=\sqrt{\mathbb{E}\left[ r_k^2 \right]}$ and $\nu_k= \mu C_1 n \delta_k^2$ for all~$k\in\mathbb{Z}_{> 0}$ and defining~$T_{k'}=r_1^2 +\mu^2 C_1^2 n^2 \sum_{k=1}^{k'}  \delta_k^4$ for all~$k'\in\mathbb{Z}_{> 0}$, we may use Lemma~\ref{lem:inequality_recursion} to conclude that
\begin{align*}
    \textstyle\sqrt{\mathbb{E}\left[ r_{k'}^2 \right]} & \le \tfrac{1}{2}\mu C_1 n \textstyle\sum_{k=1}^{k'} \delta_k^2 + \left( r_1^2 + {\mu^2 C_1^2 n^2 } \textstyle\sum_{k=1}^{k'} \delta_k^4 + ( \tfrac{1}{2}{\mu C_1 n }\textstyle\sum_{k=1}^{k'} \delta_k^2 )^2 \right)^{\frac{1}{2}} \\
    &\le {\mu C_1 n }\textstyle\sum_{k=1}^{K} \delta_k^2 + r_1 + ({\mu^2 C_1^2 n^2 } \textstyle\sum_{k=1}^{K} \delta_k^4)^{\frac{1}{2}}\quad \forall k\leq K,
\end{align*}
where the second inequality holds because $\sqrt{a+b+c}\leq \sqrt{a}+\sqrt{b}+\sqrt{c}$ for all~$a,b,c\geq 0$ and because the sums increase when we increase their upper limits from~$k'$ to~$K$. Next, consider the estimate~\eqref{ineq:f_k_suboptimality_sum} for $k'=K+1$, replace $\mathbb E[r_{K+1}^2]$ with its trivial lower bound~0 and replace~$ \sqrt{\mathbb{E}[ r_{k}^2]}$ with the above upper bound for every $k\leq K$. Noting that~$r_1=R$ and dividing by~$K$ then yields
\begin{align*}
    \tfrac{1}{K}& \textstyle\sum_{k=1}^{K} \mathbb{E} \left[ f(x_k) - f(x^\star) \right] \\
    \le &\, \tfrac{1}{\mu K} R^2 + \tfrac{1}{K}{\mu C_1^2 n^2 } \textstyle\sum_{k=1}^{K}  \delta_k^4\\
    &+ \tfrac{1}{K}{C_1 n } \textstyle\sum_{k=1}^{K} \delta_k^2 \left({\mu C_1 n }\textstyle\sum_{k=1}^{K} \delta_k^2 + R + ({\mu^2 C_1^2 n^2 } \sum_{k=1}^{K} \delta_k^4)^{\frac{1}{2}} \right) \\
    = &\, \tfrac{1}{\mu K} R^2 + \tfrac{1}{K}{C_1 n }R \textstyle\sum_{k=1}^{K} \delta_k^2 + \tfrac{1}{K}{\mu C_1^2 n^2  } (\sum_{k=1}^{K} \delta_k^2)^2 \\
    & + \tfrac{1}{K}{\mu C_1^2 n^2  } (\textstyle\sum_{k=1}^{K} \delta_k^2) (\textstyle\sum_{k=1}^{K} \delta_k^4)^{\frac{1}{2}} + \tfrac{1}{K}{\mu C_1^2 n^2 } \textstyle\sum_{k=1}^{K}  \delta_k^4.
\end{align*}
As~$\mathbb{E}[ f(\bar{x}_{K}) - f(x^\star)] \le \tfrac{1}{K} \sum_{k=1}^{K}\mathbb{E} [ f(x_k) - f(x^\star) ]$ by Jensen's inequality, assertion~{\em (i)} thus follows. If~$\delta_k = \delta \in (0,\kappa\bar{\delta}]$ for all~$k\in\mathbb{Z}_{> 0}$, then assertion~{\em (i)} implies that
\begin{align*}
    \mathbb{E} \left[ f(\bar{x}_{K}) - f(x^\star) \right]
    \leq  &\, \tfrac{1}{\mu K} R^2 + C_1 n R  \delta^2  + C_1^2 K \mu n^2  \delta^4  + C_1^2\textstyle\sqrt{K} \mu  n^2  \delta^4 + C_1^2 \mu n^2  \delta^4 \\
    \le &\, \tfrac{1}{\mu K} R^2 + {C_1 n R  \delta^2} +(C_1\textstyle\sqrt{K}+C_1)^2  \mu n^2 \delta^ 4\\
    \leq &\, \tfrac{1}{K}{2n{L_1}} R^2 + {C_1 n R  \delta^2} +(C_1\textstyle\sqrt{K}+C_1)^2  \tfrac{1}{L_1}n \delta^4\\
    \leq &\, \tfrac{1}{K}{2nL_1} R^2 + C_1 n R  \delta^2 +C_1^2 (1+ \textstyle\sqrt{K})^2  \tfrac{1}{L_1}n \delta^4,
\end{align*}
where the last two inequalities exploit the assumption~$\mu=1/(2nL_1)$. Thus, assertion~{\em (ii)} follows. Next, assume that $\delta_k =  \delta/k$ for all $k \in\mathbb{Z}_{> 0}$. In analogy to the proof of assertion~{\em (ii)}, we combine assertion~{\em (i)} with the standard zeta function inequalities~\eqref{equ:zeta} to conclude that 
\begin{align*}
    \mathbb{E} \left[ f(\bar{x}_{K}) - f(x^\star) \right]
     \leq & \tfrac{1}{\mu K} R^2 + \tfrac{1}{6} \pi^2C_1 \tfrac{1}{K} nR  \delta^2+  \tfrac{1}{90}\pi^4C_1^2 \tfrac{1}{K}\mu n^2 \delta^4\\
    & +  \tfrac{1}{36} \pi^ 4 C_1^2 \tfrac{1}{K} \mu n^2\delta^4 +   \tfrac{1}{6\sqrt{90}} \pi^4 C_1^2  \tfrac{1}{K}  \mu n^2 \delta^4 \\
    \leq &\, \tfrac{n}{K} 2 {L_1} R^2 + \tfrac{n}{K} R \tfrac{1}{6} \pi^2C_1  \delta^2 + \tfrac{n}{K} R \pi^4(\tfrac{1}{6}C_1+\tfrac{1}{\sqrt{90}}C_1)^2 \tfrac{1}{2L_1} \delta^4 \\
    \le &\, \tfrac{n}{K} ( \textstyle\sqrt{2 {L_1}} R + C_2 \delta^2 )^2, 
\end{align*}
where~$C_2= \pi^2 (C_1/3+C_1/\sqrt{90})/\sqrt{2 L_1}$. The third inequality holds because~$\mu=1/(2nL_1)$. Thus, assertion~\emph{(iii)} follows. This completes the proof for~$\mathcal X=\mathcal D$.

In the last part of the proof we show that the three assertions remain valid when~$\mathcal{X}$ is a non-empty closed convex subset of~$\mathcal{D}$. Indeed, as the projection~$\Pi_{\mathcal X}$ onto~$\mathcal X$ is contractive, we have
\begin{align*}
    r_{k+1}^2&=  \|x_{k+1}-x^\star\|^2_2 = \|\Pi_{\mathcal X} (x_k - \mu_k g_{\delta_k}(x_k)) - \Pi_{\mathcal X} (x^{\star})\|_2^2 \leq \|x_k - \mu_k g_{\delta_k}(x_k) - x^{\star}\|_2^2.
\end{align*}
Thus, all arguments used above carry over trivially to situations where~$\mathcal X\neq \mathcal D$.
\end{proof}

Theorem~\ref{thm:unconstrained}~{\em{(iii)}} shows that if $\delta_k$ decays as $O(1/k)$, then one needs $O\left({n{L_1}R^2}/{\epsilon} \right)$ iterations to guarantee that $\mathbb{E}\left[f(\Bar{x}_{K})-f(x^{\star}) \right]\leq \epsilon$. This is the first-order complexity scaled by $n$~\cite[\S~2.1.5]{nesterov2003introductory}. 
{
Theorem~\ref{thm:unconstrained} can be extended to a larger class of convex optimization problems by relaxing the assumption of constant stepsizes  \cite{ref:JongeneelIZO2021}. In particular, it can be extended to constrained optimization problems whose constraints are binding at optimality, in which case $\nabla f(x^{\star})\neq 0$; see also Example~\ref{ex:MPC} below. 
}


\section{Strongly convex optimization}
\label{sec:strongly:convex}

We now extend the results from Section~\ref{sec:convex} to analytic objective functions $f$ that are ${\tau}$-strongly convex over their domain~$\mathcal{D}$ for some ${\tau}>0$,~\textit{i.e.}, we assume that $f(y)\geq f(x) + \langle \nabla f(x), y-x \rangle + \tfrac{1}{2}{\tau}\|y-x\|_2^2$ $\forall x,y\in \mathcal{D}$. If $y$ is a stationary point with $\nabla f(y)=0$, then $\tau$-strong convexity ensures that 
\begin{equation}
\label{equ:strong:ineq}
    f(y)-f(x) \geq \tfrac{1}{2}{\tau}\|y-x\|_2^2\quad \forall x\in \mathcal{D},
\end{equation}
which in turn implies via the Polyak-\L ojasiewicz inequality $\|\nabla f(x)\|_2^2\geq 2 \tau(f(x)-f(y))$ for $\tau$-strongly convex functions \cite[Eq.~2.1.19]{nesterov2003introductory} that
\begin{equation}
    \label{equ:Lip:strongly:convex}
    \|\nabla f(x)\|_2\geq {\tau}\|y-x\|_2. 
\end{equation}

\begin{theorem}[Convergence rate of Algorithm~\ref{alg:convex_unconstrained} for strongly convex optimization]
\label{thm:strong:convex}
Suppose that all assumptions of Theorem~\ref{thm:unconstrained}~(iii) are satisfied and that~$f$ is $\tau$-strongly convex for some~$\tau>0$. Then, there is a constant $C\geq 0$ such that the following inequality holds for all $K\in \mathbb{Z}_{>0}$.
\begin{equation*}
    \mathbb{E}[f(x_{K})-f(x^{\star})] \leq \tfrac{1}{2}{L_1} \left(\delta^2 C+  ( 1 - \tfrac{{\tau}}{4n {L_1}})^{K-1}\left(R^2-\delta^2 C\right)\right)
\end{equation*}
\end{theorem}

An explicit formula for~$C$ in terms of $n$, $L_1$, $L_2$ and~$\tau$ is derived in the proof.

\begin{proof}[Proof of Theorem~\ref{thm:strong:convex}]
As in the proof of Theorem~\ref{thm:unconstrained}, we set~$C_1=3(\tfrac{1}{6} L_2  + C_\kappa)$ and~$r_{k}=\|x_k-x^{\star}\|_2$  for all~$k \in \mathbb{Z}_{>0}$, and we initially assume that $\mathcal{X}=\mathcal{D}$. Combining the estimate~\eqref{eq:r2-estimate} from the proof of Theorem~\ref{thm:unconstrained} with the strong convexity condition~\eqref{equ:strong:ineq} yields $\mathbb{E} \left[ r_{k+1}^2 | \mathcal{F}_k \right]
    \le \left(1 - \tfrac{\mu \tau}{2}\right) r_k^2+ \mu{C_1 n \delta_k^2 }r_k+ \mu^{2} C_1^2 n^2\delta_k^4$. By taking unconditional expectations and applying Jensen's inequality, we then find
\begin{subequations}
\begin{align}
    \label{equ:rho:tau:iter}
    \mathbb{E} [ r_{k+1}^2 ]& \leq  \left( 1 - \tfrac{\mu \tau}{2}\right) \mathbb{E} [ r_{k}^2 ] + \mu C_1 n\delta_k^2 \textstyle\sqrt{\mathbb{E} [ r_{k}^2 ]} + \mu^{2} C_1^2 n^2\delta_k^4\\
    & \leq \mathbb{E} [ r_{k}^2 ] + \mu C_1n\delta_k^2 \textstyle\sqrt{\mathbb{E} [ r_{k}^2 ]}+ \mu^{2} C_1^2 n^2\delta_k^4.
\end{align}
\end{subequations}
Next, choose any~$k'\in\mathbb{Z}_{>0}$ and sum the above inequalities over all $k\leq k'-1$ to obtain
\begin{align*}
    \mathbb{E} [ r_{k'}^2 ] & \leq r^2_1 + \mu C_1n\textstyle\sum^{k'-1}_{k=1} \delta_k^2 \sqrt{\mathbb{E} [ r_{k}^2 ]}  + \mu^{2} C_1^2 n^2\textstyle\sum^{k'-1}_{k=1}\delta_k^4 \\
    & \leq r^2_1 + \mu C_1 n\textstyle\sum^{k'}_{k=1}\delta_k^2 \sqrt{\mathbb{E} [ r_{k}^2 ]} + \mu^{2} C_1^2 n^2\textstyle\sum^{k'}_{k=1}\delta_k^4. 
\end{align*}
By using the same reasoning as in the proof of Theorem~\ref{thm:unconstrained}, the last bound implies
\begin{align*}
    \textstyle\sqrt{\mathbb{E}\left[ r_{k'}^2 \right]}  \le {\mu C_1n }\textstyle\sum_{k=1}^{k'} \delta_k^2 + r_1 + ({\mu^2 C_1^2 n^2 } \textstyle\sum_{k=1}^{k'} \delta_k^4)^{\frac{1}{2}}.
\end{align*}
Substituting this inequality into~\eqref{equ:rho:tau:iter} for $k=k'$ and noting that $r_1=R$ yields
\begin{align*}
    \mathbb{E} [ r_{k'+1}^2 ] \leq  \left( 1 - \tfrac{\mu \tau}{2}\right) \mathbb{E} [ r_{k'}^2 ]  &+ \mu^2C_1 ^2n^2 \delta_{k'}^4\\ 
    &+ \mu C_1 n \delta_{k'}^2 \left( {\mu C_1 n }\textstyle\sum_{k=1}^{k'} \delta_k^2 + R + ({\mu^2 C_1^2 n^2 } \textstyle\sum_{k=1}^{k'} \delta_k^4)^{\frac{1}{2}} \right).
\end{align*}
As $\delta_k = \delta/k$ for all $k\in \mathbb{Z}_{>0}$ and as the constant stepsize satisfies $\mu=1/(2nL_1)$,  we may then use the standard zeta function inequalities~\eqref{equ:zeta} to obtain
\begin{align*}
    \mathbb{E} [ r_{k'+1}^2 ]\! \leq &  ( 1 - \tfrac{\tau}{4nL_1}) \mathbb{E} [ r_{k'}^2 ]   + C_1^2 \tfrac{\delta^4}{4L_1^2(k')^4} + C_1^2 \tfrac{\pi^2 \delta^4}{24 L_1^2 (k')^2} +C_1R \tfrac{\delta^2}{2L_1 (k')^2} + C_1^2 \tfrac{\pi^2\delta^4}{4\sqrt{90} L_1^2(k')^2}\\
    \leq & ( 1 - \tfrac{{\tau}}{4n {L_1}}) \mathbb{E} [ r_{k'}^2 ]  + C_1 R \tfrac{\delta^2}{L_1} + 3C_1^2 \tfrac{\delta^4}{L_1^2},
\end{align*}
where the last inequality follows from the elementary bounds $\tfrac{1}{2(k')^2}<1$, $\tfrac{1}{4(k')^4}<1$, $\tfrac{\pi^2}{24(k')^2}<1$ and $\pi^2/(4\sqrt{90}(k')^2)<1$. As~$|\delta|<1$, we may set~$C=\frac{4n}{\tau}( C_1R + 3C_1^2/L_1)$ to obtain
\begin{align*}
    \mathbb{E} [ r_{k'+1}^2 ] \leq  ( 1 - \tfrac{{\tau}}{4n {L_1}}) \mathbb{E} [ r_{k'}^2 ]  + \tfrac{\tau}{4n L_1} \delta^2 C.
\end{align*}
Taken together, the Lipschitz inequality~\eqref{equ:grad:Lipschitz} and the strong convexity inequality~\eqref{equ:Lip:strongly:convex} imply that~$\tau\leq L_1$, which in turn ensures that~$\tau/(4nL_1)<1$. Hence, the above inequality implies $\mathbb{E} [ r_{k'+1}^2 ] -\delta^2C  \leq ( 1 - \tfrac{{\tau}}{4n {L_1}}) \left( \mathbb{E} ([ r_{k'}^2 ]  -\delta^2C\right)$. As this estimate holds for all~$k'<K$, we may finally conclude that
\begin{equation*}
    \mathbb{E} [ r_{K}^2 ] - \delta^2 C  \leq  ( 1 - \tfrac{{\tau}}{4n {L_1}}) \left(\mathbb{E} [ r_{K-1}^2 ]- \delta^2 C \right) \leq  \cdots \leq ( 1 - \tfrac{{\tau}}{4n {L_1}})^{K-1} (R-\delta^2 C). 
\end{equation*}
The claim then follows by combining this inequality with the estimate~$\mathbb{E}[f(x_K)-f(x^{\star})]\leq \frac{1}{2}{L_1} \mathbb{E} [ r_{K}^2 ]$, which follows from the Lipschitz condition~\eqref{equ:L1:upper}. This completes the proof for~$\mathcal X=\mathcal D$. To show that the claim remains valid when~$\mathcal{X}$ is a non-empty closed convex subset of~$\mathcal{D}$, we may proceed as in the proof of Theorem~\ref{thm:unconstrained}. Details are omitted for brevity.
\end{proof}

By Theorem~\ref{thm:strong:convex} and the construction of~$C$, we can enforce $\mathbb{E}[f({x}_{K})-f(x^{\star}) ]\leq \epsilon$ for a given tolerance~$\epsilon>0$ by selecting a sufficiently small smoothing parameter~$\delta\leq O(\sqrt{\epsilon \tau/(n L_1^2)})$ and by running Algorithm~\ref{alg:convex_unconstrained} over $O({n {L_1}}/{{\tau}}\log({{L_1}R^2}/{\epsilon} ) )$ iterations. 


\section{Non-convex optimization}
\label{sec:nonconvex}

We now extend the convergence guarantees for Algorithm~\ref{alg:convex_unconstrained} to unconstrained non-convex optimization problems. Our proof strategy differs from the one in~\cite{ref:nesterov2017random} as the smoothed objective function $f_{\delta}$ does not necessarily admit a Lipschitz continuous gradient. In this setting, convergence can still be guaranteed if the initial iterate~$x_1$ is sufficiently close to some global minimizer $x^{\star}$.

\begin{theorem}[Convergence rate of Algorithm~\ref{alg:convex_unconstrained} for nonconvex optimization]
\label{thm:ncvx}
Suppose that all assumptions of Theorem~\ref{thm:unconstrained}~(iii) hold, but assume that~$f$ may be non-convex, $\mathcal{X}=\mathcal{D}$ and $\mu_k=\mu=1/(nL_1)$ for all~$k\in\mathbb{Z}_{> 0}$. Define $F=f(x_1)-f(x^{\star})$, where~$x^\star$ is a global minimizer of problem~\eqref{equ:opt:main}. 
If~$\|\nabla f(x_1)\|_2^2\leq 2nL_1F$, then there is a constant $C\geq 0$ such that for all $K\in \mathbb{Z}_{>0} $ we have
\begin{equation*}
     \tfrac{1}{K}\textstyle\sum^{K}_{k=1}\mathbb{E}\left[\|\nabla  f(x_k)\|_2^2 \right]\leq \tfrac{n}{K}\left( 2 {L_1}F +\delta^2 C \right).
\end{equation*}
\end{theorem}
The dependence of~$C$ on $n$, $L_1$, $L_2$ and~$F$ can be derived from the proof of Theorem~\ref{thm:ncvx}. 
\begin{proof}[Proof of Theorem~\ref{thm:ncvx}]
As~$\mathcal{X}=\mathcal{D}$, the iterates of Algorithm~\ref{alg:convex_unconstrained} satisfy $x_{k+1}=x_k - \mu_k\cdot g_{\delta_k}(x_k)$. In addition, as $f$ has a Lipschitz continuous gradient, the Lipschitz inequality~\eqref{equ:L1:upper} implies that  
\begin{align*}
    f(x_{k+1}) \leq & f(x_k) - \mu_k \langle \nabla f(x_k), g_{\delta_k}(x_k) \rangle + \tfrac{1}{2}\mu_k^2 {L_1} \|g_{\delta_k}(x_k)\|_2^2\\
    = f(&x_k) - \mu_k\|\nabla f(x_k)\|_2^2 - \mu_k \langle \nabla f(x_k),g_{\delta_k}(x_k)-\nabla f(x_k)\rangle + \tfrac{1}{2}\mu_k^2 {L_1} \|g_{\delta_k}(x_k)\|_2^2.
\end{align*}
Taking conditional expectations on both sides of this expression, recalling that~$g_{\delta_k}$ is an unbiased estimator for~$\nabla f_{\delta_k}$ conditional on~$\mathcal F_k$ and applying the Cauchy-Schwarz inequality then yields
\begin{align*}
    \mathbb{E}\left[f(x_{k+1})|\mathcal{F}_k \right] \leq & f(x_k) - \mu_k \|\nabla f(x_k)\|_2^2 \\
    & + \mu_k \|\nabla f(x_k)\|_2 \| \nabla f_{\delta_k}(x_k) - \nabla f(x_k)\|_2+ \tfrac{1}{2}\mu_k^2 {L_1} \mathbb E  [\|g_{\delta_k}(x_k)\|_2^2|\mathcal{F}_k ].
\end{align*}
Defining~$C_0=\tfrac{1}{6} L_2  + C_\kappa$, we may use the estimates~\eqref{equ:grad:delta:approx:error} and~\eqref{equ:g:var} to obtain 
\begin{align*}
    \mathbb{E}\left[f(x_{k+1})|\mathcal{F}_k \right] \leq & f(x_k) - \mu_k \|\nabla f(x_k)\|_2^2 + \mu_k {C_0 n\delta_k^2 }\|\nabla f(x_k)\|_2 \\
    & + \tfrac{1}{2}\mu_k^2 {L_1}\left( n\|\nabla f(x_k)\|_2^2 + {C_0^2 n^2  \delta_k^4 }+2C_0n^2\delta_k^2\|\nabla f(x)\|_2\right)\\
    =& f(x_k) -\tfrac{1}{2nL_1} \|\nabla f(x_k)\|_2^2 +\tfrac{1}{2{L_1}}C_0^2 \delta_k^4+ \tfrac{2}{L_1}C_0\delta_k^2\|\nabla f(x_k)\|_2,
\end{align*}
where the equality holds because the stepsize is constant and equal to~$\mu_k = 1/(n {L_1})$. By taking unconditional expectations, applying Jensen's inequality and rearranging terms, we then find
\begin{equation}
\label{equ:ncvx:grad:bound}
\begin{aligned}
   \mathbb{E}[\|\nabla f(x_k)\|_2]^2&\leq \mathbb{E}[\|\nabla f(x_k)\|^2_2]\\ 
   &\leq 2nL_1 \mathbb{E}[f(x_k)-f(x_{k+1})]+4nC_0 \delta_k^2\mathbb{E}\left[\|\nabla f(x_k)\|_2\right] +nC_0^2 \delta_k^4. 
\end{aligned}
\end{equation}
Next, choose any~$k'\in\mathbb{Z}_{>0}$ and sum the left- and rightmost terms in~\eqref{equ:ncvx:grad:bound} over all $k\leq k'$ to obtain 
\begin{align*}
    &\mathbb{E}[\|\nabla f(x_{k'})\|_2]^2\\ & \leq \textstyle\sum^{k'}_{k=1} \mathbb{E}[\|\nabla f(x_k)\|_2]^2 \\
    & \leq 2n{L_1}\mathbb E[f(x_1)-f(x_{k'+1})] + 4nC_0\textstyle\sum^{k'}_{k=1}\delta_k^2 \mathbb{E}[\|\nabla f(x_k)\|_2] + nC_0^2\textstyle\sum^{k'}_{k=1}\delta_k^4\\
    & \leq 2n{L_1} F+ 4nC_0\textstyle\sum^{k'}_{k=1}\delta_k^2 \mathbb{E}[\|\nabla f(x_k)\|_2] + nC_0^2\textstyle\sum^{k'}_{k=1}\delta_k^4,
\end{align*}
where the third inequality holds because $x^\star$ is a global minimizer of problem~\eqref{equ:opt:main}, which implies $\mathbb E[f(x_1)-f(x_{k'+1})] = \mathbb E[f(x_1)-f(x^\star)] + \mathbb E[f(x^\star)-f(x_{k'+1})] \leq F$. Setting~$t_k=\mathbb{E}[\|\nabla f(x_k)\|_2]$ and~$\nu_k=4nC_0\delta_k^2$ for all~$k\in\mathbb{Z}_{> 0}$, and defining $T_{k'}=2nL_1F+nC_0^2\sum^{k'}_{k=1}\delta_k^4$ for all~$k'\in\mathbb{Z}_{>0}$, we may then use Lemma~\ref{lem:inequality_recursion}, which applies because~$\|\nabla f(x_1)\|_2^2\leq 2n {L_1}F$, to find
\begin{align*}
    \mathbb{E}[\|\nabla f(x_{k'})\|_2] \leq 2 n C_0 \textstyle\sum^{k'}_{k=1}\delta_k^2 + \left(2nL_1F + {nC_0^2}\textstyle\sum^{k'}_{k=1}\delta^4_k + (2n C_0 \textstyle\sum^{k'}_{k=1}\delta_k^2 )^2 \right)^{\frac{1}{2}}.
\end{align*}
As $\delta_k = \delta/k$ for all $k\in \mathbb{Z}_{>0}$,  the standard zeta function inequalities~\eqref{equ:zeta} imply that
\begin{align*}
    \mathbb{E}[\|\nabla f(x_{k'})\|_2] & \leq n C_0\delta^2 \tfrac{\pi^2}{3} + \left(2nL_1F + {nC_0^2}\delta^4 \tfrac{\pi^4}{90} + n^2 C^2_0 \delta^4\tfrac{\pi^4}{9} \right)^{\frac{1}{2}}\\
    & \leq \textstyle\sqrt{2nL_1F} + nC_0 \delta^2( \tfrac{2\pi^2}{3}  + \tfrac{\pi^2}{\sqrt{90}}) ,
\end{align*}
where the second inequality holds because $\sqrt{a+b+c}\leq \sqrt{a}+\sqrt{b}+\sqrt{c}$ for all~$a,b,c\geq 0$ and because~$\sqrt{n}\leq n$ for all~$n\in\mathbb{Z}_{\geq 0}$.  Averaging the second inequality in~\eqref{equ:ncvx:grad:bound} across all $k\leq K$ and using the above upper bound on $\mathbb{E}[\|\nabla f(x_{k})\|_2]$ for each~$k\leq K$ finally yields
\begin{align*}
    \tfrac{1}{K}\textstyle\sum^K_{k=1} \mathbb{E}\left[\|\nabla f(x_k)\|_2^2 \right] \leq & \tfrac{n}{K}\Big[  2L_1 F + 4C_0\textstyle\sum_{k=1}^K \delta_k^2\left( \sqrt{2nL_1 F} + nC_0 \delta^2( \tfrac{2\pi^2}{3}  + \tfrac{\pi^2}{\sqrt{90}})\right)\\& +C_0^2 \textstyle \sum_{k=1}^K \delta_k^4 \Big].
\end{align*}
Applying the zeta function inequalities~\eqref{equ:zeta} once again and recalling that~$\delta_k^2\leq 1$ for all~$k\in\mathbb{Z}_{> 0}$, it is then easy to construct a constant~$C\geq 0$ such that\\ $\frac{1}{K}\sum^K_{k=1} \mathbb{E}[\|\nabla f(x_k)\|_2^2] \le \frac{n}{K} ( 2 L_1F + \delta^2 C)$.
\end{proof}

By Theorem~\ref{thm:ncvx}, we can enforce $\tfrac{1}{K}\sum^{K}_{k=1}\mathbb{E}[\|\nabla f(x_k)\|_2^2]\leq \epsilon$ for a given $\epsilon>0$ by selecting a smoothing parameter~$\delta\leq O(\sqrt{K\epsilon/n})$  and by running Algorithm~\ref{alg:convex_unconstrained} over $O(nL_1 F /\epsilon )$ iterations.

\section{Numerical experiments}
\label{sec:num}

We will now assess the empirical performance of different variants of Algorithm~\ref{alg:convex_unconstrained} equipped with different gradient estimators on standard test problems. Specifically, we will compare the proposed complex-step estimator~$g_{\mathsf{cs}}$ defined in~\eqref{equ:grad:est:g} against the forward-difference estimator
\begin{align*}
    g_{\mathsf{fd}}(x,\delta) &= \tfrac{1}{\delta}({f(x+\delta y)-f(x)}) y\quad \text{with}\quad y\sim \mathcal{N}(0,I_n)
\end{align*}
and the central-difference estimator 
\begin{align*}
    g_{\mathsf{cd}}(x,\delta) &= \tfrac{1}{2\delta}({f(x+\delta y)-f(x-\delta y)}) y\quad \text{with}\quad y\sim \mathcal{N}(0,I_n),
\end{align*}
both of which rely on Gaussian smoothing~\cite[Eq.~(30)]{ref:nesterov2017random}; see also Example~\ref{ex:num_or}. As pointed out in the introduction, the single-point estimator~\eqref{equ:g:rho:N} displays a higher variance and thus leads to slow convergence, in general. Therefore, we exclude it from the numerical experiments. When using~$g_{\mathsf{fd}}$ or~$g_{\mathsf{cd}}$, we set the stepsize of Algorithm~\ref{alg:convex_unconstrained} to $\mu_k= 1/( 4(n+4){L_1})$ as recommended in~\cite[Eq.~(55)]{ref:nesterov2017random}. When using~$g_{\mathsf{cs}}$, on the other hand, we select the stepsize in view of the structural properties of the given objective function~$f$ in accordance with Theorems~\ref{thm:unconstrained}, \ref{thm:strong:convex} and~\ref{thm:ncvx}. {We recall that when $\delta_k$ is fixed, $\delta_k \approx 10^{-8}$ is optimal for $g_{\mathsf{fd}}$ whereas $\delta_k\approx 10^{-5}$ is optimal for $g_{\mathsf{cd}}$, \textit{cf.}~Figure~\ref{fig:x0-1}.} Also, the initial iterate $x_1$ is always set to $0$ unless stated otherwise. {All experiments are performed in MATLAB on a x86\_64 machine with a 4 GHz CPU and 16 GB RAM, using double precision, that is, machine precision is $2^{-52}\approx 2.2204\cdot 10^{-16}$.} 

From Sections~\ref{sec:convex}--\ref{sec:nonconvex} we know that Algorithm~\ref{alg:convex_unconstrained} with~$g_{\mathsf{cs}}$ is guaranteed to find stationary points of a wide range of convex and non-convex optimization problems provided that its stepsize is inversely proportional to the Lipschitz constant~$L_1$ of the gradient of~$f$. Implementing Algorithm~\ref{alg:convex_unconstrained} in practice thus requires knowledge of~$L_1$. Unfortunately, the Lipschitz modulus of~$f$ is typically unknown in the context of zeroth-order optimization, and the results of Sections~\ref{sec:convex}--\ref{sec:nonconvex} indicate that increasing~$L_1$ increases the number of iterations and decreases the smoothing parameter~$\delta$ needed to attain a desired suboptimality gap~$\epsilon$. These insights are consistent with classical results in zeroth-order optimization based on multi-point gradient estimators such as~$g_{\mathsf{fd}}$ or~$g_{\mathsf{cd}}$ (\textit{cf.}~\cite{ref:nesterov2017random}). As the complex-step method proposed in this paper remains numerically stable for almost arbitrarily small smoothing parameters~$\delta$, it may thus be preferable to classical methods when~$L_1$ is overestimated.

The experiments will show that if $\delta$ is sufficiently large for multi-point methods to be applicable, then our complex-step method converges equally fast or faster than the multi-point methods, which obey the theoretical convergence rates reported in~\cite{ref:nesterov2017random}. A theoretical explanation for the better empirical \textit{transient} convergence behavior of the complex-step method is left for future work.

\begin{figure*}[t!]
    \centering
    \begin{subfigure}[b]{0.24\textwidth}
        \includegraphics[width=\textwidth]{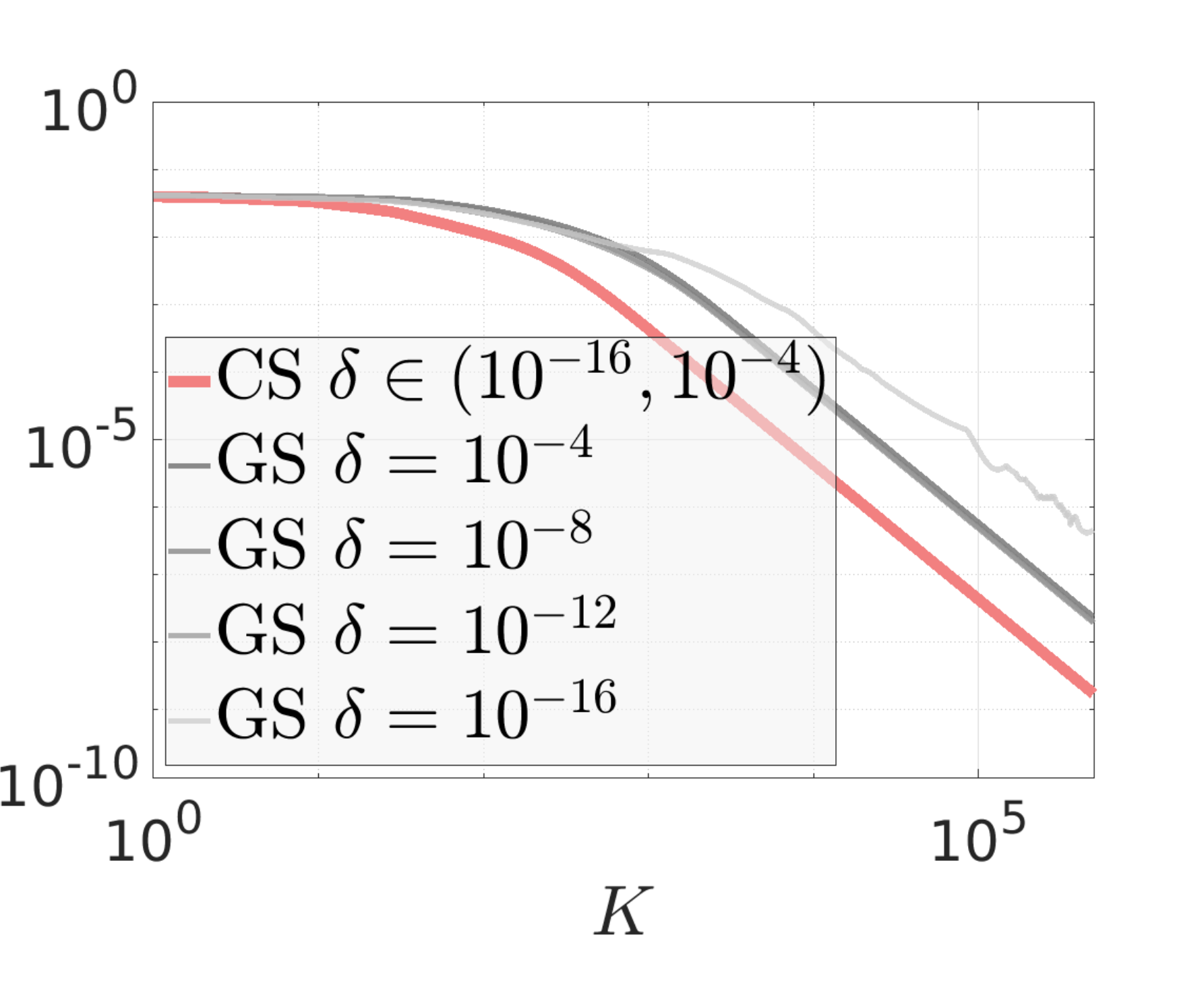}
        \caption{Suboptimality gap $f(\bar{x}_K)-f^{\star}$ for~\eqref{equ:nesterov:test}.}
        \label{fig:quad:test:xkbar}
    \end{subfigure}
         \begin{subfigure}[b]{0.24\textwidth}
        \includegraphics[width=\textwidth]{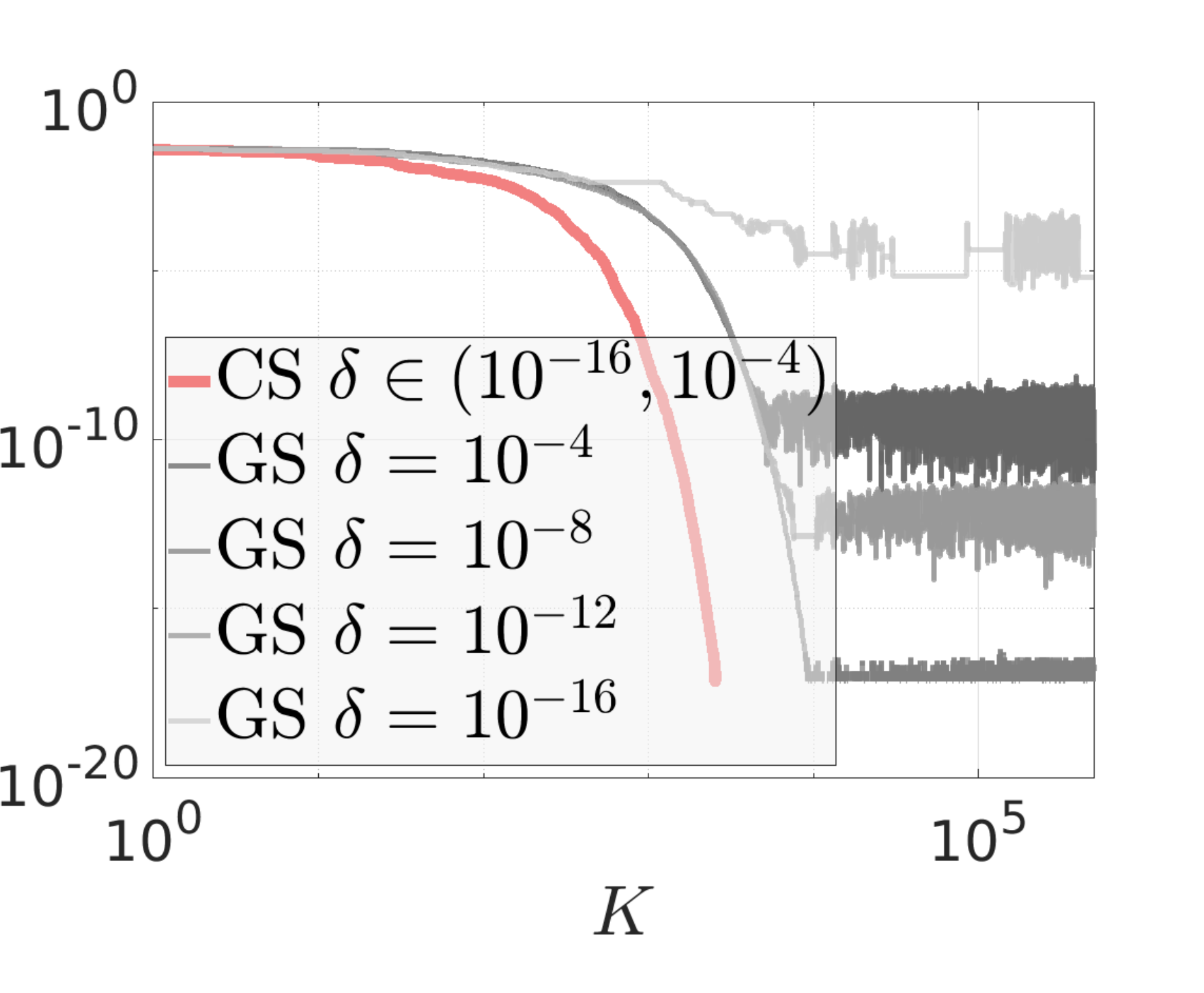}
        \caption{Suboptimality gap $f({x}_K)-f^{\star}$ for~\eqref{equ:nesterov:test}.}
        \label{fig:quad:test:xk}
    \end{subfigure}
    \begin{subfigure}[b]{0.24\textwidth}
        \includegraphics[width=\textwidth]{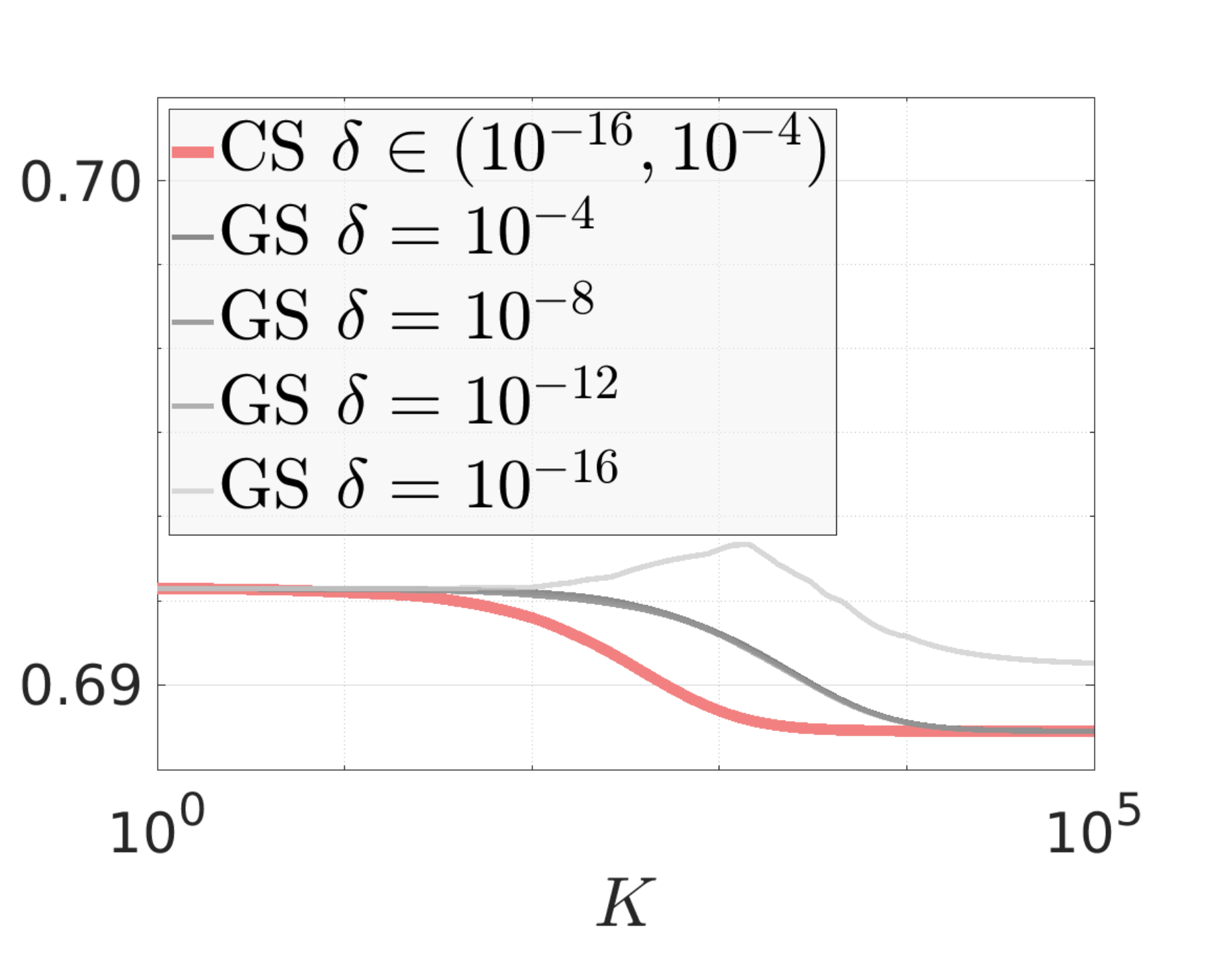}
        \caption{Cost $f(\bar{x}_K)$ for the log loss function~\eqref{equ:logistic}.}
        \label{fig:log:xkbar}
    \end{subfigure}
    \begin{subfigure}[b]{0.24\textwidth}
        \includegraphics[width=\textwidth]{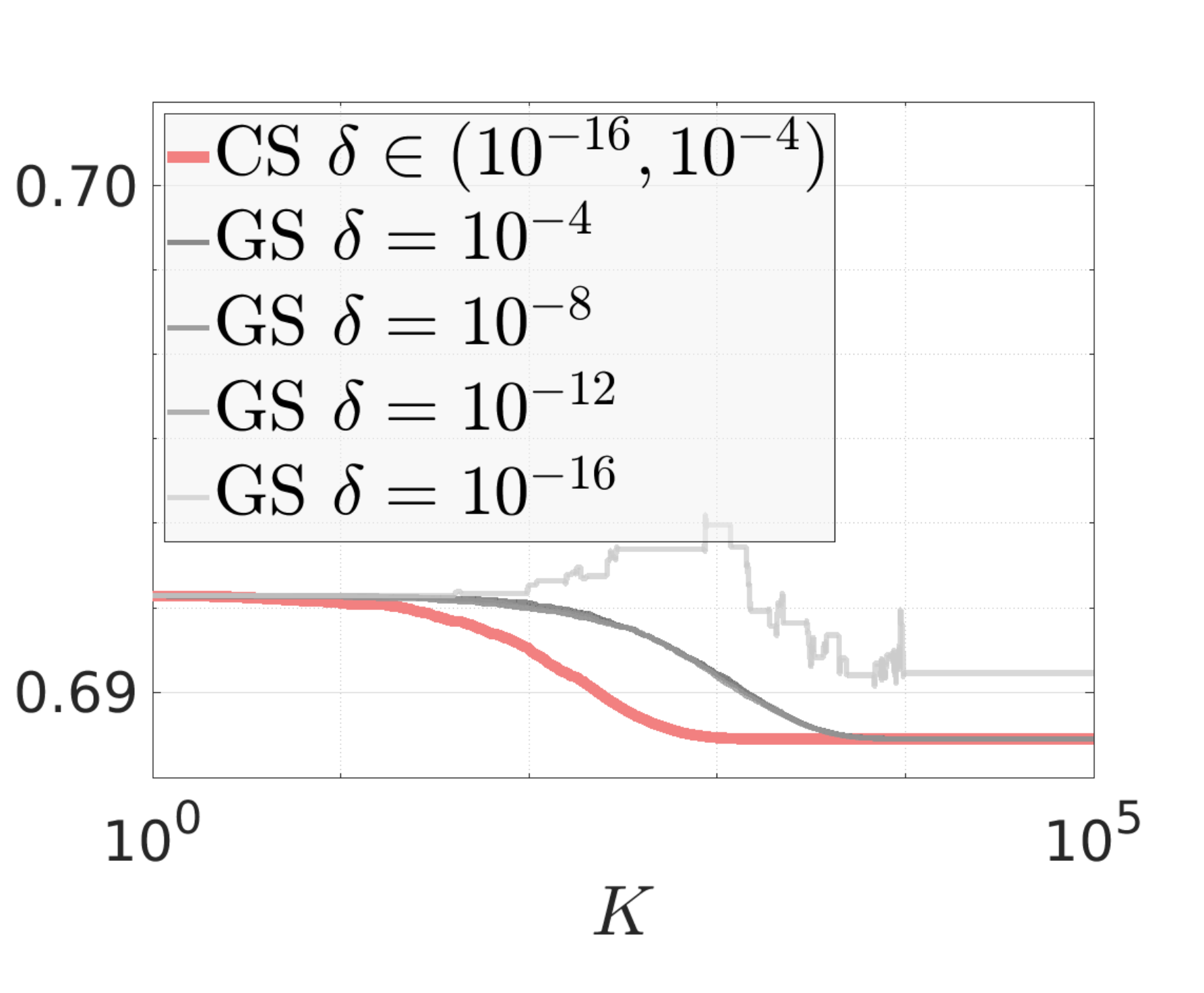}
        \caption{Cost $f({x}_K)$ for the log loss function~\eqref{equ:logistic}.}
        \label{fig:log:xk}
    \end{subfigure}
    \caption[]{Comparison of the single-point complex-step estimator $g_{\mathsf{cs}}$ ($\mathrm{CS}$) against the multi-point Gaussian smoothing estimator~$g_{\mathsf{fd}}$ ($\mathrm{GS}$) on two objective functions.}
    \label{fig:unconstrained}
\end{figure*}

\subsection{Unconstrained convex optimization}
\begin{example}[Quadratic test function]
\label{ex:quadratic_function}
\upshape{
Assume that~$\mathcal X=\mathbb{R}^n$ and~$f$ is an ill-conditioned version of what Nesterov calls the `\textit{worst function in the world}' \cite[\S~2.1.2]{nesterov2003introductory}, that is, assume that
\begin{equation}
    \label{equ:nesterov:test}
    f(x) = {L}\left(\tfrac{1}{2}\left[(x^{(1)})^2 + \textstyle\sum^{n-1}_{j=1}(x^{(j+1)}-x^{(j)})^2 + (x^{(n)})^2 \right]-x^{(1)}\right),
\end{equation}
where $n=5$, $L=10^{-8}$, and~$x^{(j)}$ denotes the $j^{\mathrm{th}}$ component of $x$ for any~$j\leq n$. One can show that $\nabla f$ has Lipschitz modulus~${L_1}=4 L$ and that the unique global minimizer~$x^\star$ of~$f$ has coordinates~$(x^{\star})^{(j)}=1-j/(n+1)$. In this case, the theoretical convergence guarantees of Algorithm~\ref{alg:convex_unconstrained} are independent of whether~$g_{\mathsf{cs}}$ or~$g_{\mathsf{fd}}$ is used. However, starting from~$x_1=0$ and gradually reducing the smoothing parameter~$\delta_k=\delta$ towards machine precision exposes the advantages of the one-point estimator~$g_{\mathsf{cs}}$ over the multi-point estimator~$g_{\mathsf{fd}}$. Figures~\ref{fig:quad:test:xkbar} and~\ref{fig:quad:test:xk} visualize the suboptimality gap of~$\bar x_K$ and~$x_K$ as a function of~$K$ along a single sample trajectory, respectively. Note that especially the performance of $x_K$ is significantly better when~$g_{\mathsf{cs}}$ is used. {One might argue that, {despite using the optimal $\delta\approx 10^{-8}$}, $g_{\mathsf{fd}}$ leads to a higher suboptimality gap than $g_{\mathsf{cs}}$ because of its inferior approximation quality; see, \textit{e.g.}, Figure~\ref{fig:oracles}. We shed more light on this conjecture in Example~\ref{ex:dim-dep}, where we compare $g_{\mathsf{cs}}$ against $g_{\mathsf{cd}}$.}  
}
\end{example}


\begin{example}[Logistic regression]
\upshape{
Assume that~$\mathcal X=\mathbb{R}^n$ and~$f$ is the log loss function used to quantify the prediction loss in logistic regression. Specifically, set
\begin{equation}
    \label{equ:logistic}
    f(x) = \tfrac{1}{m}\textstyle\sum^m_{i=1}\log\left(1+\exp\left(-v_i a_i^{\mathsf{T}}x \right) \right)
\end{equation}
for~$m=100$ and~$n=2$, and assume that the features~$a_i$ and the labels~$v_i$ are sampled independently from the standard normal distribution on~$\mathbb R^n$ and the uniform distribution on~$\{-1,1\}$, respectively. Denoting by~$A\in\mathbb R^{m\times n}$ the matrix with rows~$a_i^{\mathsf{T}}$ for all~$i\leq m$, one readily verifies that ${L_1}=\frac{1}{m}\|A\|_2$. We compare again the empirical convergence properties of Algorithm~\ref{alg:convex_unconstrained} equipped with~$g_{\mathsf{cs}}$ or~$g_{\mathsf{fd}}$. Figures~\ref{fig:log:xkbar} and~\ref{fig:log:xk} visualize the objective function values of~$\bar{x}_K$ and~$x_K$ as a function of~$K$. We observe that the cancellation effects in the cost of~$\bar{x}_K$ are mild even if~$g_{\mathsf{fd}}$ is used and~$\delta_k=\delta$ is small, whereas those in the cost of~$x_K$ are significantly more pronounced. 
}
\end{example}

\begin{figure*}[t!]
    \centering
    \begin{subfigure}[b]{0.24\textwidth}
        \includegraphics[width=\textwidth]{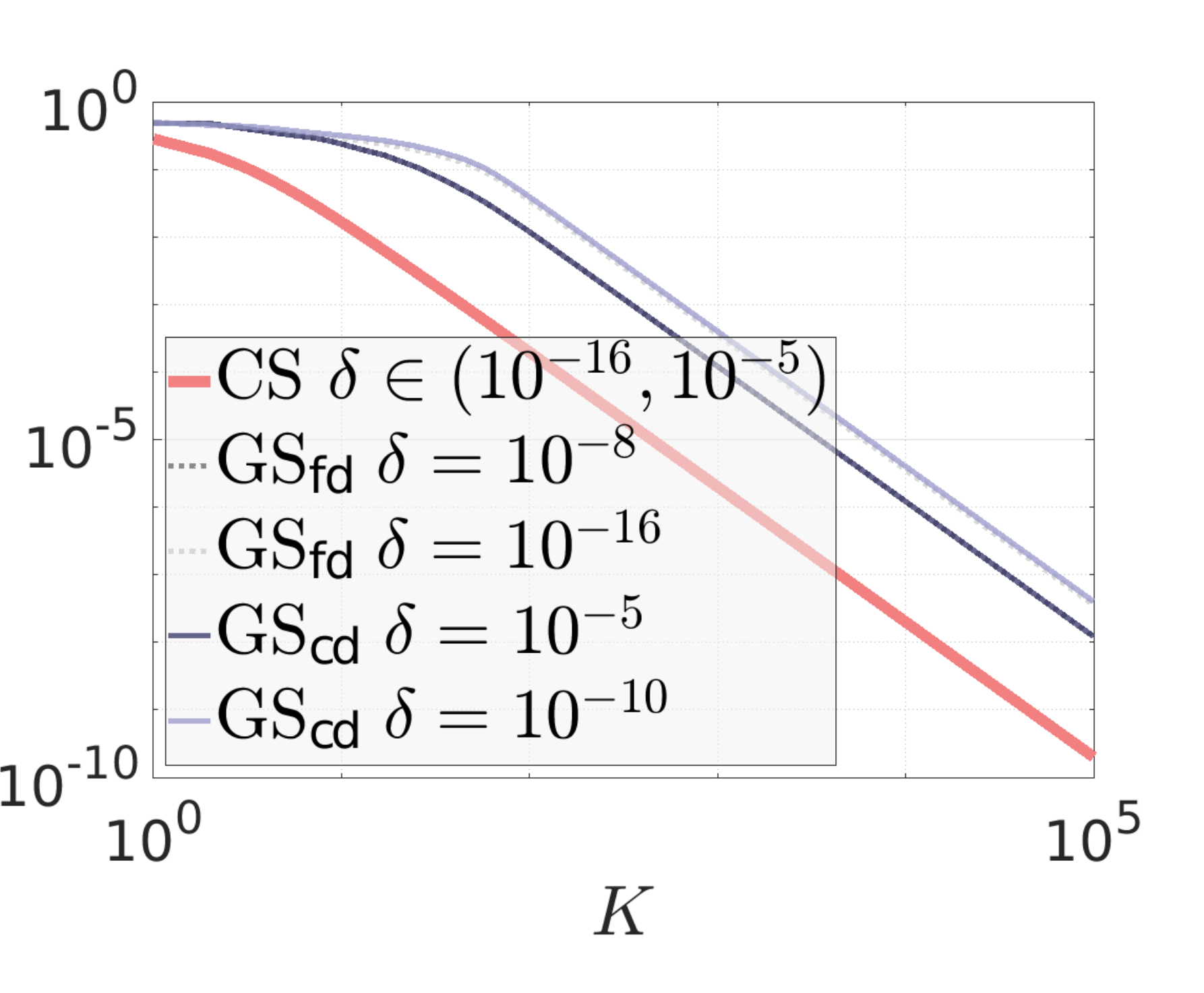}
        \caption{Suboptimality gap $f(\bar{x}_K)\!-\!f^{\star}$ for~$n=1$. }
        \label{fig:quadavgn1}
    \end{subfigure}\,
    \begin{subfigure}[b]{0.24\textwidth}
        \includegraphics[width=\textwidth]{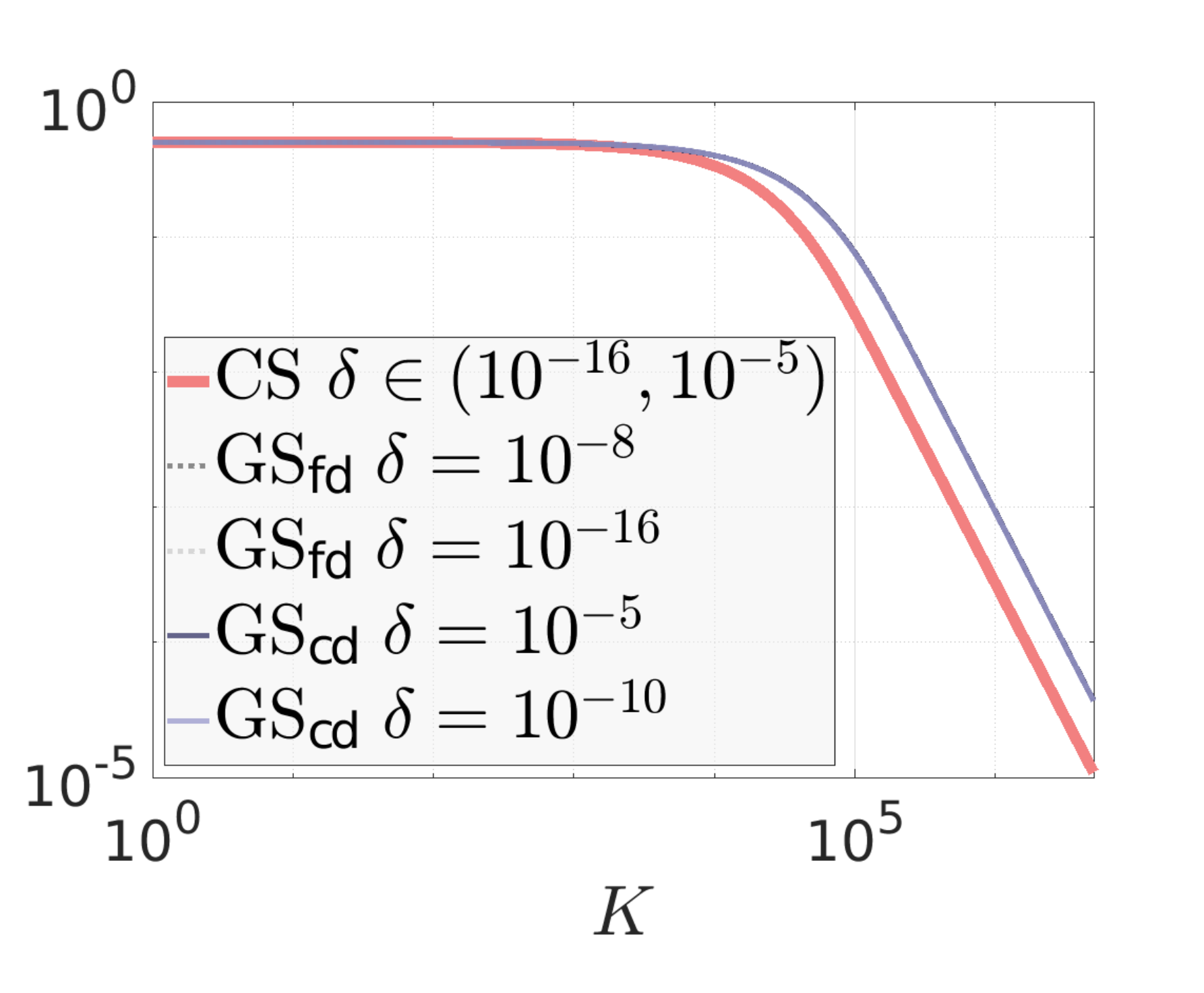}
        \caption{Suboptimality gap $f(\bar{x}_K)\!-\!f^{\star}$ for~$n=10^4$.}
        \label{fig:quadavgn100}
    \end{subfigure}\,
    \begin{subfigure}[b]{0.24\textwidth}
        \includegraphics[width=\textwidth]{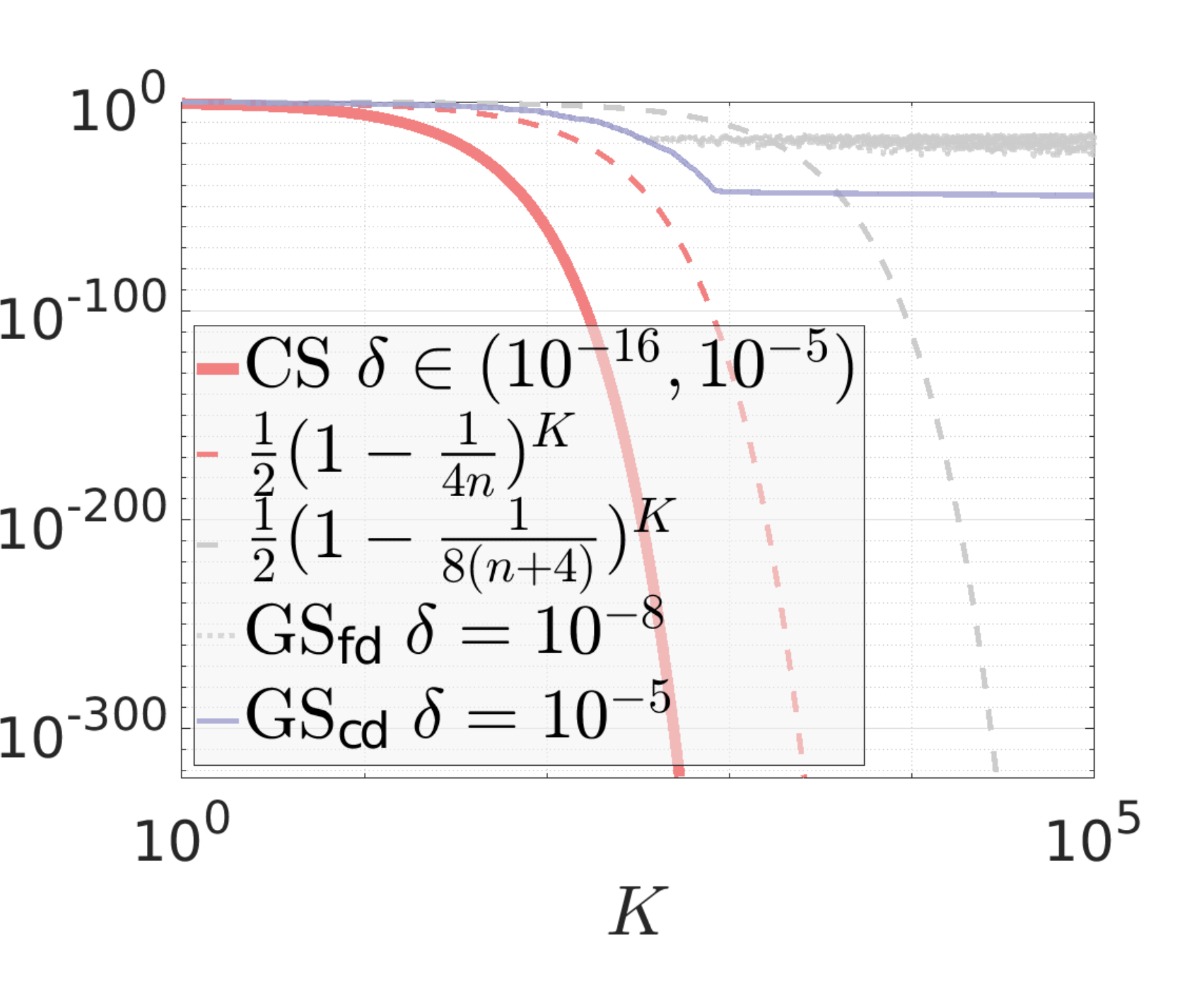}
        \caption{Suboptimality gap $f({x}_K)\!-\!f^{\star}$ for~$n=1$.}
        \label{fig:quadn1}
    \end{subfigure}\,
    \begin{subfigure}[b]{0.24\textwidth}
        \includegraphics[width=\textwidth]{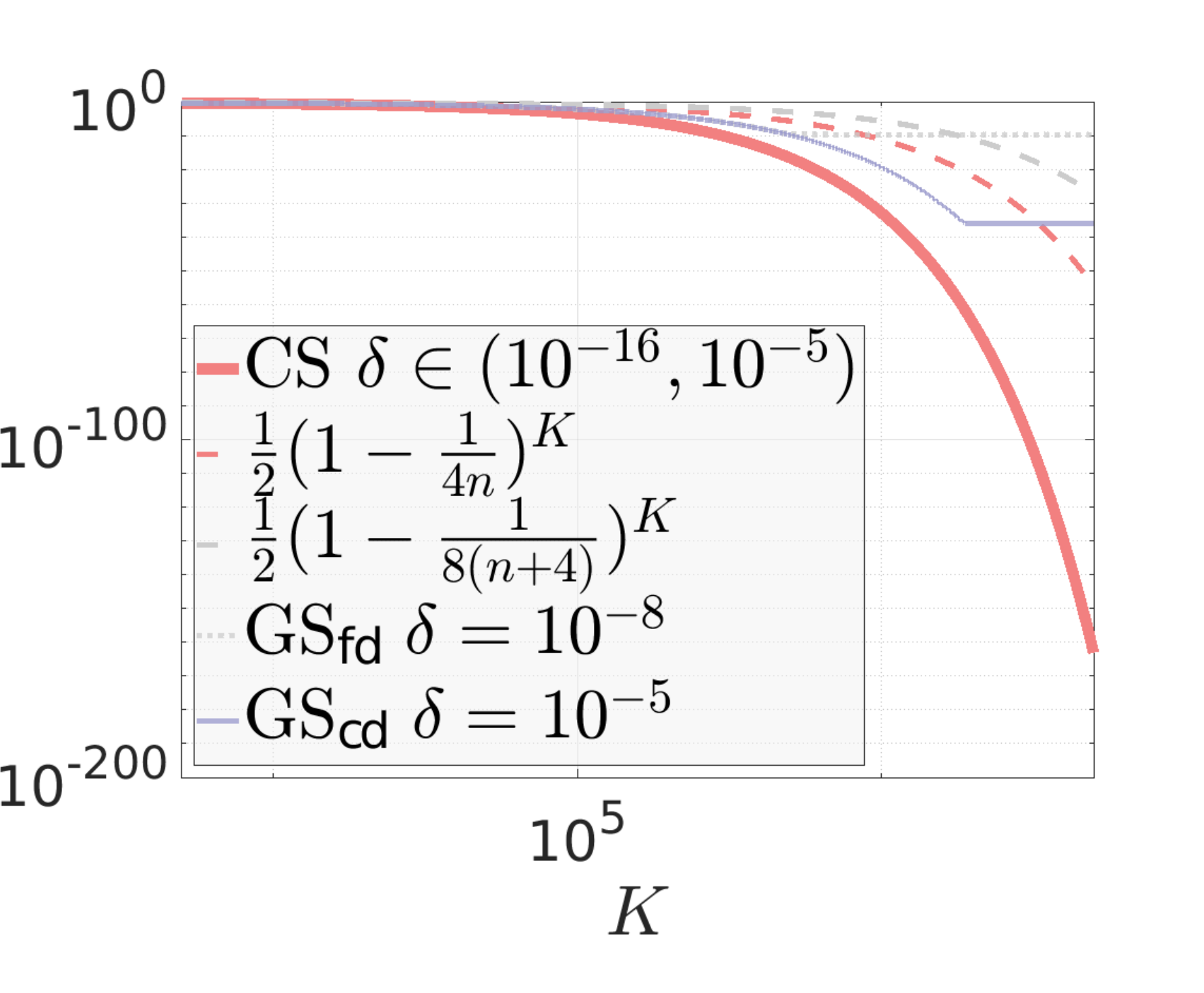}
        \caption{Suboptimality gap $f({x}_K)\!-\!f^{\star}$ for~$n=10^4$.}
        \label{fig:quadn100}
    \end{subfigure}
    \caption[]{Comparison of the single-point complex-step estimator~$g_{\mathsf{cs}}$ ($\mathrm{CS}$) against multi-point Gaussian smoothing estimators~$g_{\mathsf{fd}}$ ($\mathrm{GS}_{\mathsf{fd}}$) and~$g_{\mathsf{cd}}$ ($\mathrm{GS}_{\mathsf{cd}}$) on~$f(x)=\frac{1}{2}\|x\|_2^2$.}
    \label{fig:unconstrained:dim}
\end{figure*}

\begin{example}[Dimension-dependence and strong convexity]
\label{ex:dim-dep}
\upshape{
Figure~\ref{fig:unconstrained} not only confirms that the complex-step estimator~$g_{\mathsf{cs}}$ is less susceptible to cancellation effects than the forward-difference estimator~$g_{\mathsf{fd}}$, but it also suggests that Algorithm~\ref{alg:convex_unconstrained} converges faster if~$g_{\mathsf{cs}}$ is used instead of~$g_{\mathsf{fd}}$. In view of Example~\ref{ex:num_or}, we further expect that the central-difference estimator~$g_{\mathsf{cd}}$ should lead to faster convergence than~$g_{\mathsf{fd}}$. To verify these conjectures numerically, we now set~$\mathcal X=\mathbb R^n$ and~$f(x)=\frac{1}{2}\|x\|_2^2$. Figure~\ref{fig:unconstrained:dim} visualizes the suboptimality gap of~$\bar x_K$ and~$x_K$ as a function of~$K$ for the three gradient estimators~$g_{\mathsf{cs}}$, $g_{\mathsf{fd}}$ and~$g_{\mathsf{cd}}$ and for increasing dimensions~$n\in\{1, 10{,}000\}$, starting from $x_1=n^{-1/2}\mathbbold{1}$. Cancellation effects prevail even in this simple example, {for optimal smoothing parameters $\delta_k=\delta$}. We observe that~$\bar x_K$ is less susceptible to cancellation effects but converges significantly slower than~$x_K$. Even though the central-difference estimator~$g_{\mathsf{cd}}$ does indeed provide a speed-up compared to the finite difference estimator~$g_{\mathsf{fd}}$, it is still dominated by the complex-step estimator~$g_{\mathsf{cs}}$. As~$f$ has a Lipschitz continuous gradient with~${L_1}=1$ and is strictly convex with~${\tau}=1$, Theorem~\ref{thm:strong:convex} ensures that for a negligible smoothing parameter~$\delta$ and for an initial iterate with~$\|x_1\|_2=1$ the suboptimality gap of~$x_K$ decays at least as fast as~$\tfrac{1}{2}(1-\tfrac{1}{4n})^K$. Figure~\ref{fig:unconstrained:dim} also visualizes this theoretical convergence rate and contrasts it with the rate $\tfrac{1}{2}(1-\tfrac{1}{8(n+4)})^K$ for Algorithm~\ref{alg:convex_unconstrained} with~$g_{\mathsf{cd}}$ \cite[Eq.~(57)]{ref:nesterov2017random}. 
}
\end{example}

\subsection{Constrained convex optimization}
\begin{figure*}[t!]
    \centering
     \begin{subfigure}[b]{0.24\textwidth}
        \includegraphics[width=\textwidth]{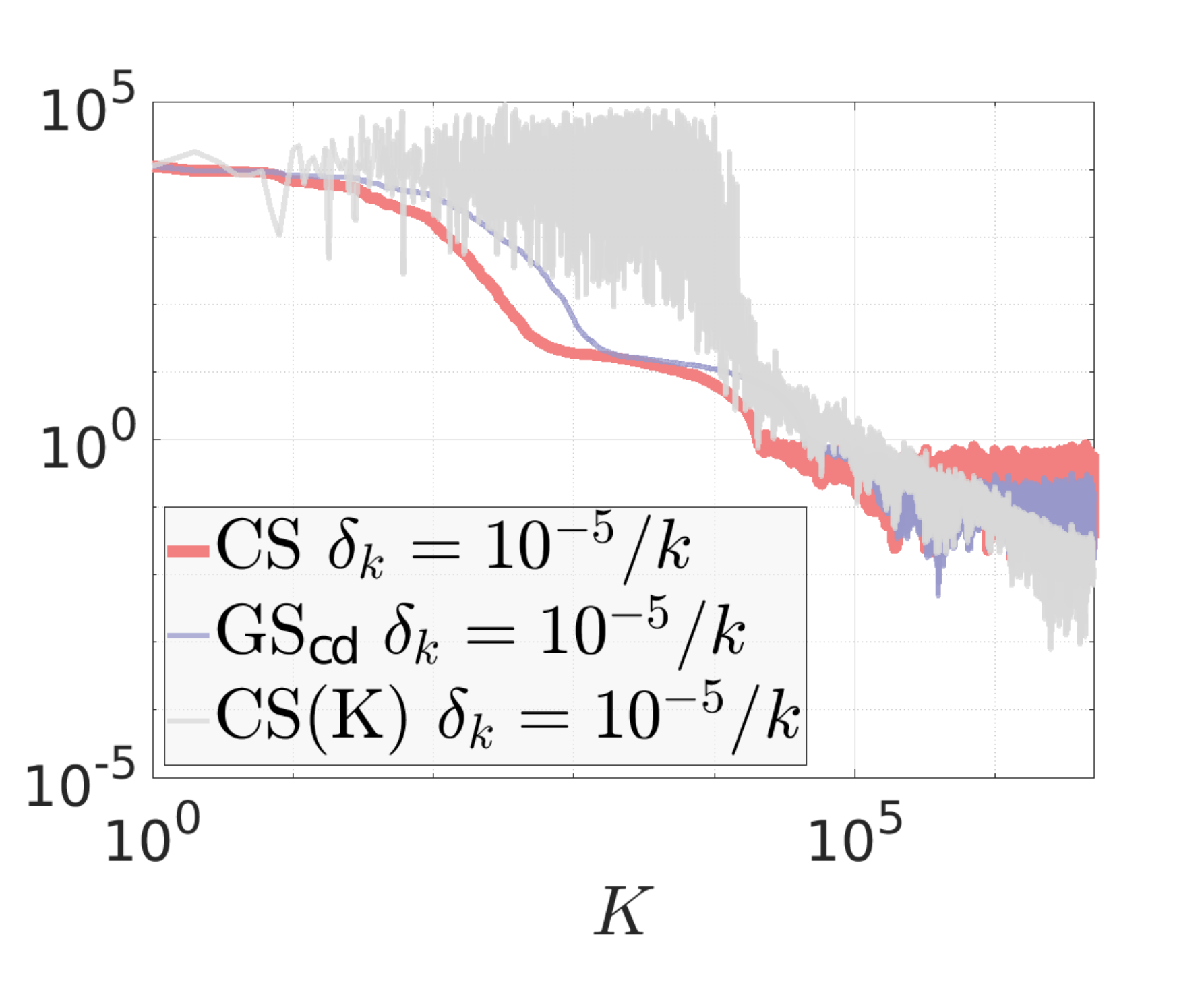}
        \caption{Suboptimality gap $f({x}_K)-f^{\star}$ for Example~\ref{ex:MPC}.}
        \label{fig:MPC}
    \end{subfigure}\,
    \begin{subfigure}[b]{0.24\textwidth}
        \includegraphics[width=\textwidth]{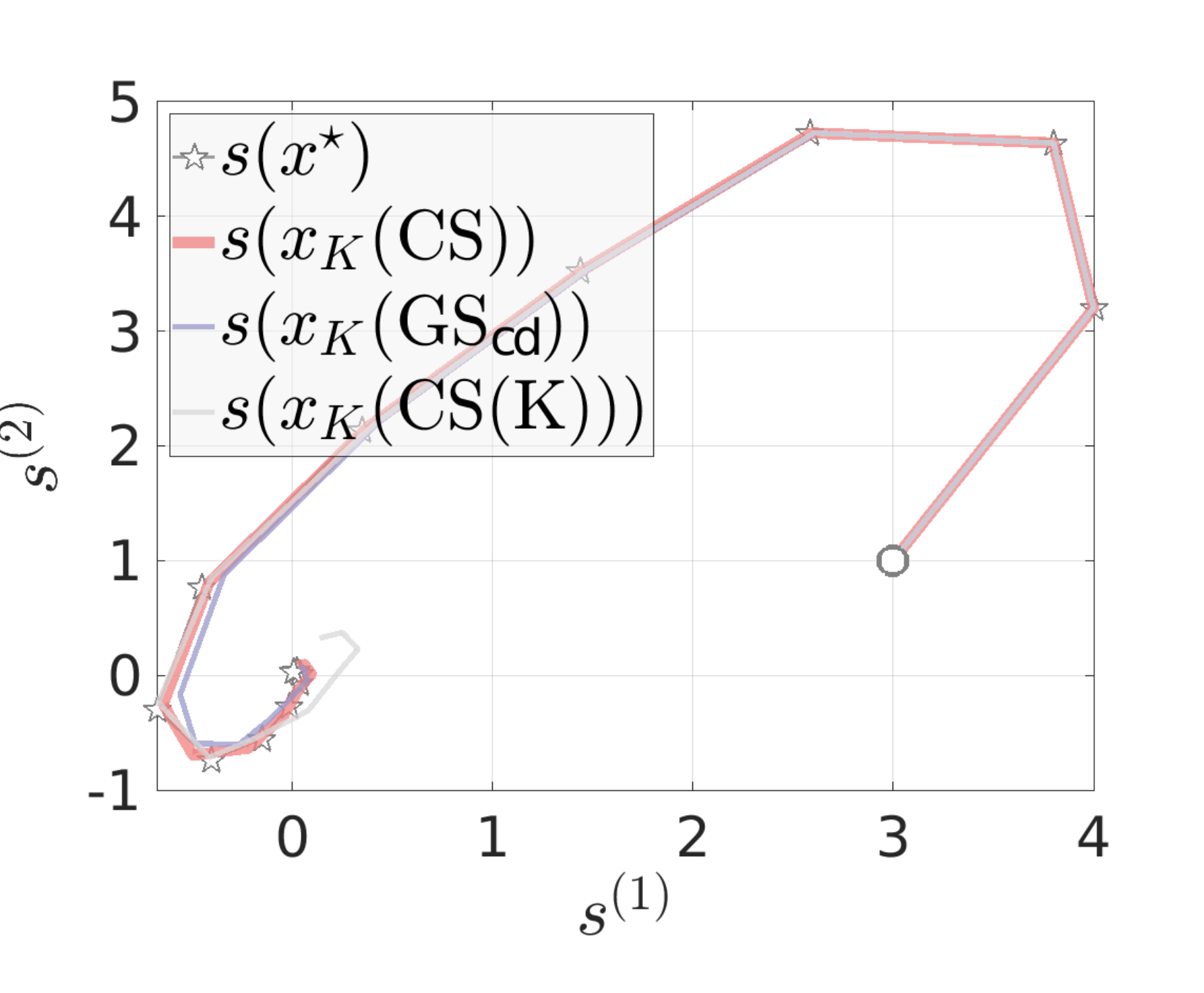}
        \caption{State trajectories generated by the control policies.}
        \label{fig:MPCx}
    \end{subfigure}\,
    \begin{subfigure}[b]{0.24\textwidth}
        \includegraphics[width=\textwidth]{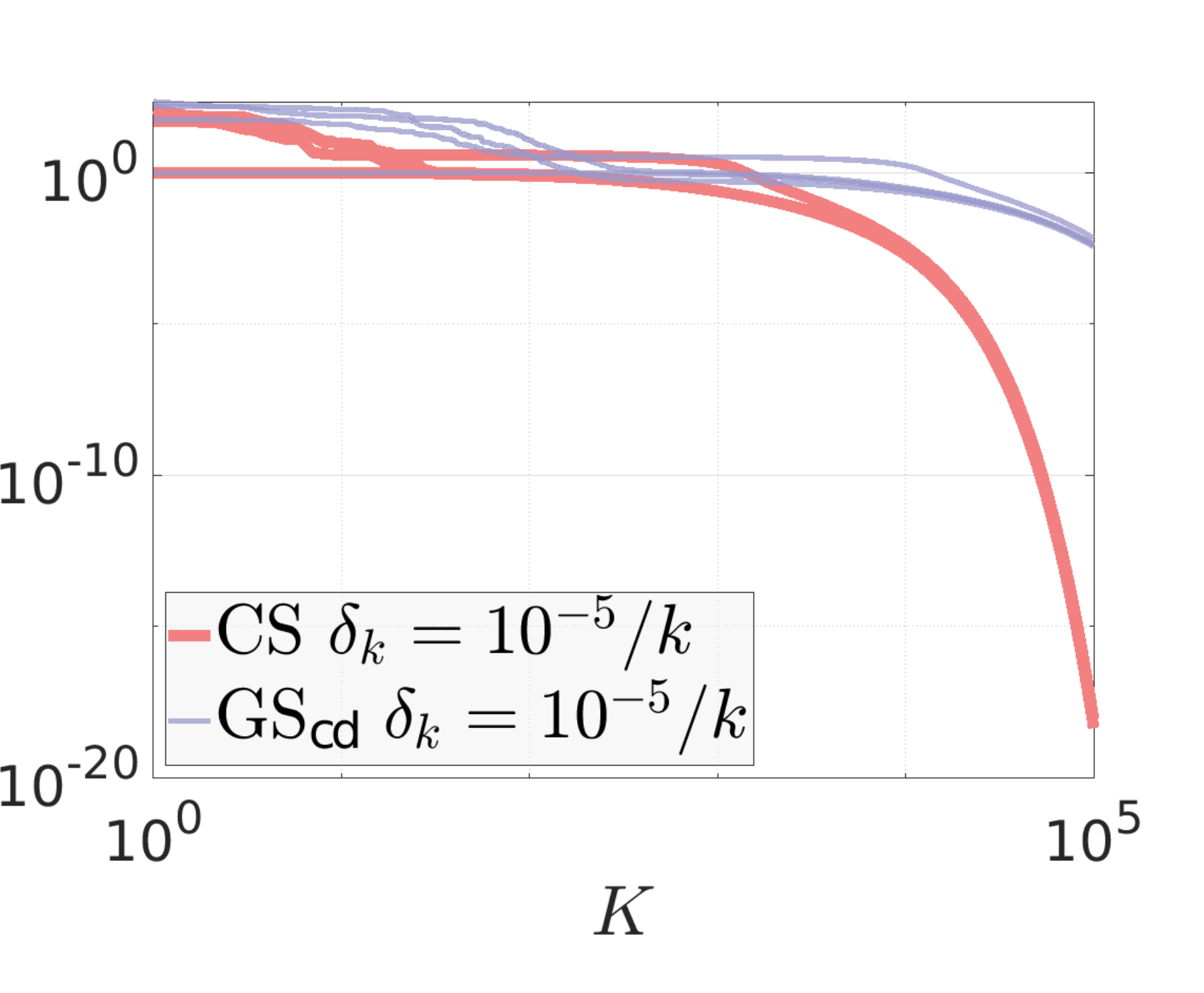}
        \caption{Suboptimality gap $f({x}_K)-f^{\star}$ for Example~\ref{ex:rosenbrock}.}
        \label{fig:ncvx:cost}
    \end{subfigure}\,
    \begin{subfigure}[b]{0.24\textwidth}
        \includegraphics[width=\textwidth]{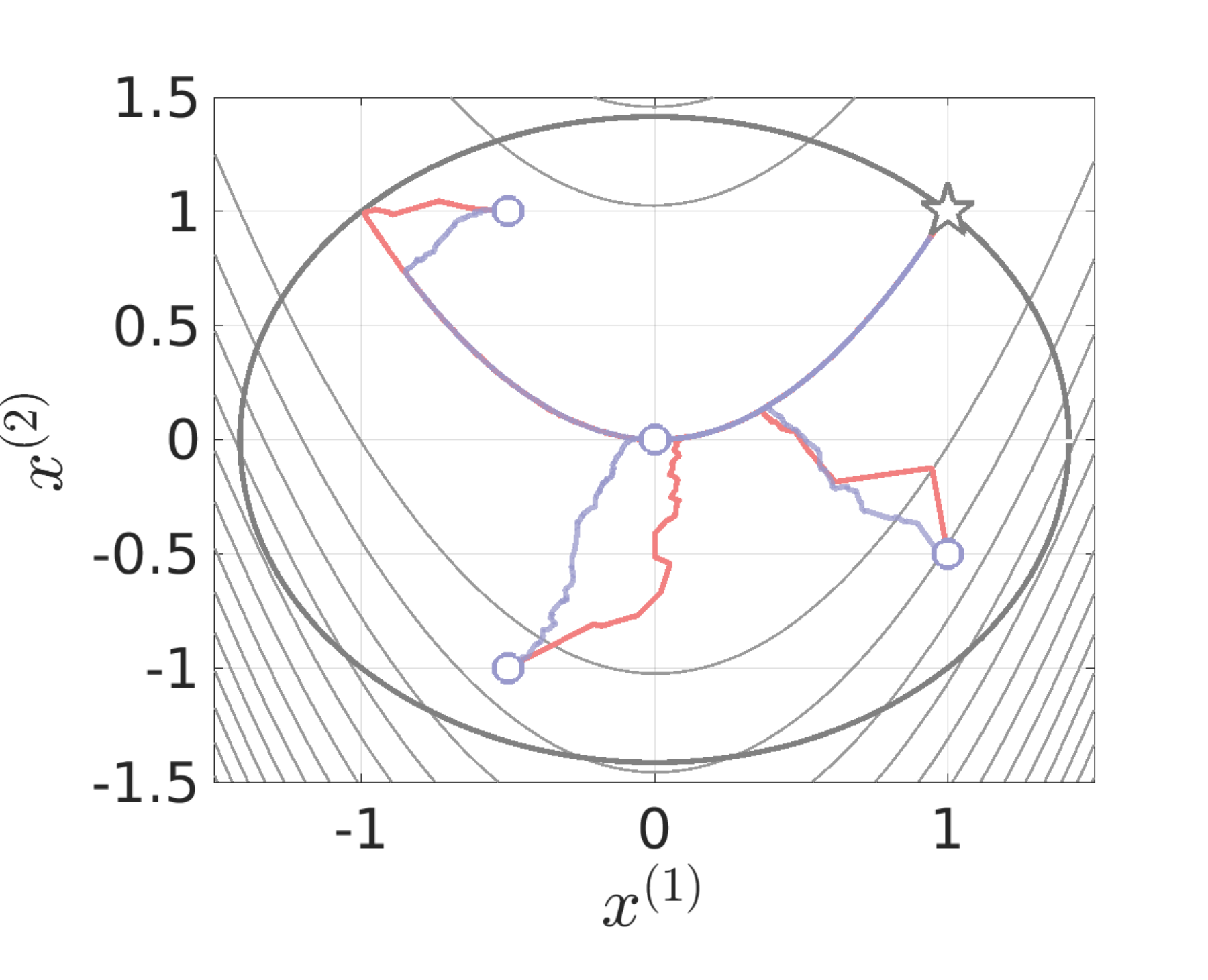}
        \caption{Paths of iterates starting at four different initial points.}
        \label{fig:ncvx:path}
    \end{subfigure}
    \caption[]{Comparison of the single-point complex-step estimator $g_{\mathsf{cs}}$ {($\mathrm{CS}$ for $\mu_k=\mu$ and $\mathrm{CS}(K)$ for $\mu_k=1/k$)} against the multi-point Gaussian smoothing estimator~$g_{\mathsf{cd}}$ ($\mathrm{GS}_{\mathsf{cd}}$) on the optimization problems of Examples~\ref{ex:MPC} and~\ref{ex:rosenbrock}.}
    \label{fig:constrained:ncvx}
\end{figure*}

The next example revolves around a constrained optimization problem grounded in optimal control. We remark that the (unconstrained) infinite-horizon version of this problem could be addressed with the policy iteration scheme proposed in~\cite{ref:Fazel_18,ref:malik2019derivative}. 
\begin{example}[Policy iteration]
\label{ex:MPC}
\upshape{
We now address the MPC problem
\begin{equation}
\label{equ:MPC}
    \begin{aligned}
    \minimize_{x=\{x_t\}_{t=0}^{T-1}\subseteq \mathbb{R}^{n_x}} \quad & \textstyle\sum^{T-1}_{t=0}\langle Qs_t,s_t \rangle + \langle Rx_t, x_t\rangle + \langle Qs_T, s_T\rangle \\
    \subjectto \hspace{6mm}& s_{t+1} = As_t + Bx_t \quad\forall t=0,\ldots, T-1 \\
    & \|x_t\|_{\infty} \leq 1\quad\forall t=0,\ldots, T-1
    \end{aligned}
\end{equation}
with planning horizon~$T\in\mathbb{Z}_{\geq 0}$ and initial state~$s_0\in\mathbb R^{n_s}$. Note that the dynamic constraints in~\eqref{equ:MPC} can be used to express the state trajectory~$s=\{s_t\}_{t=0}^{T}$ as an affine function of the $n$-dimensional input trajectory~$x=\{x_t\}_{t=0}^{T-1}$ with~$n=Tn_x$. We can thus eliminate~$s$ and express the objective function of~\eqref{equ:MPC} as a quadratic function~$f(x)$ of the inputs~$x$ alone. Similarly, we can identify the feasible set of~\eqref{equ:MPC} with the compact hypercube~$\mathcal X=\{x\in\mathbb R^n:\|x\|_\infty\leq 1\}$. Hence, the MPC problem~\eqref{equ:MPC} constitutes an instance of~\eqref{equ:opt:main}. We further assume that the cost matrices $Q\succeq 0$ and $R\succ 0$ are known, that the system matrices~$A$ and $B$ are unknown, and that the costs of a given input trajectory~$x$ can be evaluated by simulation. This implies that~$f$ is unknown but admits a zeroth-order oracle. Throughout this experiment we set $(A,B,Q,R)$ to the standard two-dimensional MPC instance\footnote{\url{https://yalmip.github.io/example/standardmpc/}} in Yalmip~\cite{ref:Lofberg2004}, and we set~$T=15$ and~$s_0=(3,1)$. We emphasize that the optimal solution of~\eqref{equ:opt:main} may reside on the boundary of~$\mathcal X$, and thus the theoretical guarantees of Sections~\ref{sec:convex} and~\ref{sec:strongly:convex} do not apply. Nevertheless, we will show that Algorithm~\ref{alg:convex_unconstrained} performs better when the complex-step estimator~$g_{\mathsf{cs}}$ is used instead of the central-difference estimator~$g_{\mathsf{cd}}$. We initialize the algorithm at the origin and upper bound the Lipschitz modulus of~$\nabla f$ by~${L_1}=4\cdot 10^4$. This crude bound is merely based on the operator norms of~$A$ and~$B$. Figure~\ref{fig:MPC} visualizes the suboptimality gap of~$x_K$ as a function of~$K$. The oscillations in the suboptimality gap corresponding to~$g_{\mathsf{cs}}$ emerge because the optimizer~$x^\star$ of~\eqref{equ:MPC} resides on the boundary of~$\mathcal X$. {Note that Theorems~\ref{thm:unconstrained} and~\ref{thm:strong:convex} do not apply even though~$f$ is strongly convex. The reason is again that~$\nabla f(x^
\star)\neq 0$. Convergence results for optimization problems with boundary solutions have recently been obtained in~\cite{ref:JongeneelIZO2021} by allowing the stepsize~$\mu_k$ to decay with~$k$. These results even hold if we have only access to noisy function evaluations. We thus solve problem~\eqref{equ:MPC} once again with Algorithm~\ref{alg:convex_unconstrained} and the complex-step estimator $g_{\mathsf{cs}}$ but set $\mu_k=1/k$ instead of $\mu_k=1/(2nL_1)$, both for $\delta_k=\delta/k$ with $\delta=10^{-5}$. In this case, \cite{ref:JongeneelIZO2021} guarantees the suboptimality gap to decay as $O(1/K)$. Figure~\ref{fig:MPC} 
empirically validates this theoretical result. Note also that initial convergence is slower under a harmonically decaying stepsize. Figure~\ref{fig:MPCx} shows the state trajectories corresponding to differently computed inputs~${x}_K$ for~$K=5\cdot 10^4$.}}
\end{example}

\subsection{Non-convex optimization}
We finally apply our method to a classical non-convex test problem. 
\begin{example}[Rosenbrock function]
\label{ex:rosenbrock}
\upshape{
Set~$\mathcal X=\sqrt{2}\mathbb{B}^2$, and let~$f$ be the Rosenbrock function defined through~$f(x)= (1-x^{(1)})^2 + 100[x^{(2)}-(x^{(1)})^2]^2$. Then, problem~\eqref{equ:opt:main} is uniquely solved by~$x^{\star}=(1,1)$, which coincides with the global minimizer of~$f$ over~$\mathbb{R}^2$. We compare the complex-step estimator~$g_{\mathsf{cs}}$ against~$g_{\mathsf{cd}}$ but remark that the convergence behavior of Algorithm~\ref{alg:convex_unconstrained} does not change noticeably when~$g_{\mathsf{fd}}$  is replaced with~$g_{\mathsf{cd}}$. We also set~$x_1$ to one of four different points in~$\mathcal X$ as visualized in Figure~\ref{fig:ncvx:path}. Figures~\ref{fig:ncvx:cost} and~\ref{fig:ncvx:path} show the convergence of the suboptimality gap of~$x_K$ and the paths of iterates generated by Algorithm~\ref{alg:convex_unconstrained}, respectively. Again, the complex-step estimator~$g_{\mathsf{cs}}$ leads to significantly faster convergence. An additional acceleration can be achieved by decreasing~$\delta$ below~$10^{-5}$. In this case, however, the Gaussian smoothing method eventually breaks down. {We point out that, compared to other derivative-free approaches such as the Nelder-Mead algorithm, the convergence is slow, and future work should aim at improving our understanding of the relative merits of these methods, \textit{e.g.}, fast empirical convergence (Nelder-Mead algorithm) versus slower but guaranteed convergence (Algorithm~\ref{alg:convex_unconstrained}).}
}
\end{example}

{
\subsection{Outlook}
\label{sec:outlook}
To close the paper, we discuss potential applications of our methods in the context of simulation-based optimization, where evaluating the objective function~$f$ requires the solution of an ordinary differential equation (ODE) or a partial differential equation (PDE). This section is illustrative only, and additional work is required to derive rigorous convergence guarantees. We start with {an optimization problem involving an ODE, which---due to its chaotic nature---often serves as a benchmark problem in the dynamical systems literature; see, \textit{e.g.}, \cite{ref:kaheman2022automatic}.}

\begin{example}[Lorenz system]
\label{ex:chaos}
The ODE
 \begin{equation}
 \label{equ:Lorenz}
     \frac{\mathrm{d}}{\mathrm{d}t}
     \begin{pmatrix}
    {{\ell}}_1(t) \\ 
    \ell_2(t)\\
    \ell_3(t)
    \end{pmatrix} = 
    \begin{pmatrix}
    \sigma(\ell_2(t)-\ell_1(t))\\
     \ell_1(t)(r-\ell_3(t)) - \ell_2(t)\\
     \ell_1(t)\ell_2(t) - b \ell_3(t) 
    \end{pmatrix}
 \end{equation}
is commonly known as the Lorenz system~\cite[Ch.~9]{ref:strogatz1994nonlinear}. It was developed as a stylized model of atmospheric convection, with $\ell_1$, $\ell_2$ and $\ell_3$ representing the rate of convection, the horizontal temperature variation and the vertical temperature variation, respectively. However, the Lorenz system also arises in the study of chemical reactions, population dynamics or electric circuits etc. In the following we denote by $\varphi^t(x)$ the time-$t$ state of a Lorenz system with initial state~$\ell(0)=x$. Given a potentially noisy measurement~$p$ of the state at time~$t\geq 0$, a problem of practical interest is to estimate the initial state~$x$ that led to~$p$. If~$x$ is known to belong to a closed set~$\mathcal X\subseteq \mathbb R^3$, then it can conveniently be estimated by solving an instance of problem~\eqref{equ:opt:main} with objective function $f(x)=\|p-\varphi^t(x)\|_2^2$ and feasible set~$\mathcal X$. We expect this problem to be challenging because the Lorenz system is known to be chaotic. Thus, slight changes in the initial state have a dramatic impact on the future trajectory. Moreover, the objective function is not available in closed form but must be evaluated with a numerical ODE solver. As all commonly used ODE solvers map the initial state~$x$ to an approximation of~$\varphi^t(x)$ by recursively applying analytic (in fact, polynomial) transformations, the resulting instance of problem~\eqref{equ:opt:main} can be addressed with Algorithm~\ref{alg:convex_unconstrained}. We remark that most out-of-the-box ODE solvers accept complex-valued initial conditions. Here we use MATLAB's \texttt{ode45} routine.  
 
In the following we set the problem parameters to $\sigma=10$, $r=28$ and $b=8/3$. In addition, we define $t=2$ and $\mathcal X=\{x\in\mathbb R^3:\|x-\ell(0)\|_2\leq 2\}$, and we sample~$p$ from the normal distribution $\mathcal N(\varphi^2(\ell(0)), \epsilon^2 I_3)$, where $\ell(0)=(10,10,10)$ and $\epsilon=10^{-3}$. Finally, we sample the initial iterate~$x_1$ from the uniform distribution on the boundary of~$\mathcal X$, use $L_1=1{,}000$ as a conservative estimate for the Lipschitz modulus of $\nabla f$ and set the smoothing parameter to $\delta=10^{-10}$.
By Theorem~\ref{thm:ncvx}, Algorithm~\ref{alg:convex_unconstrained} converges to a stationary point~$x^\star$ of the objective function with $\nabla f(x^\star)=0$ provided that $\|\nabla f(x_1)\|_2$ is sufficiently small. See the discussion below Example~\ref{ex:MPC} for the case $\nabla f(x^\star)\neq 0$. Figure~\ref{fig:L1} shows the decay of $f(x_K)$ with the total number~$K$ of iterations for 10 independent simulation runs. Figure~\ref{fig:L2} visualizes the corresponding state trajectories\footnote{A video of the state evolution is available from~\url{https://wjongeneel.nl/Lorenz.gif}.} $\{\varphi^t(x_K)\}_{t\in[0,2]}$ for~$K=10^5$. Only 6 of the 10 trajectories are shown for better visibility.

\begin{figure*}[t!]
    \centering 
    \begin{subfigure}[b]{0.3\textwidth}
        \includegraphics[width=\textwidth]{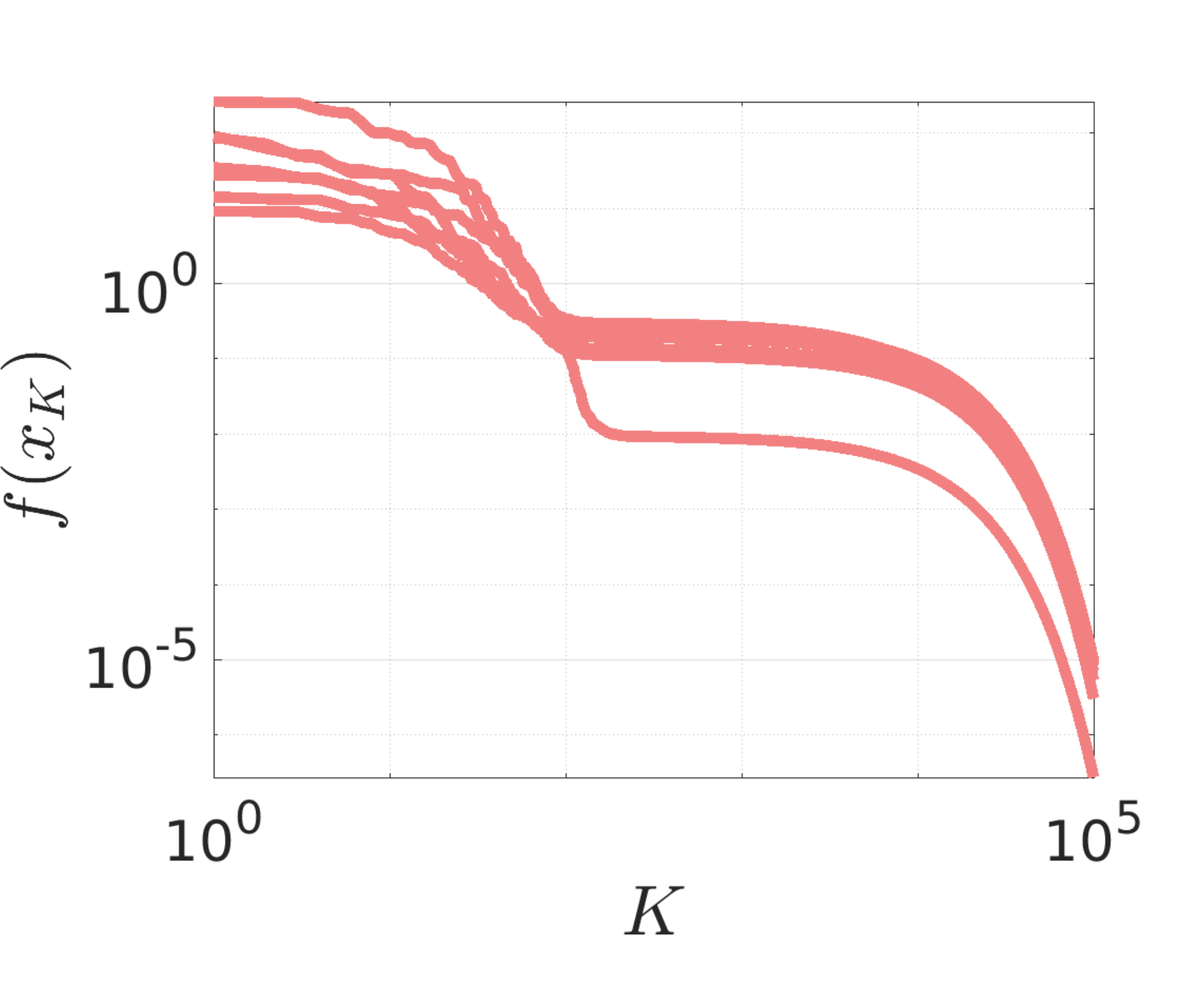}
        \caption{Decay of $f(x_K)$ with the total number~$K$ of iterations for 10 independent simulation runs.\\}
        \label{fig:L1}
    \end{subfigure}\quad 
    \begin{subfigure}[b]{0.3\textwidth}
        \includegraphics[width=\textwidth]{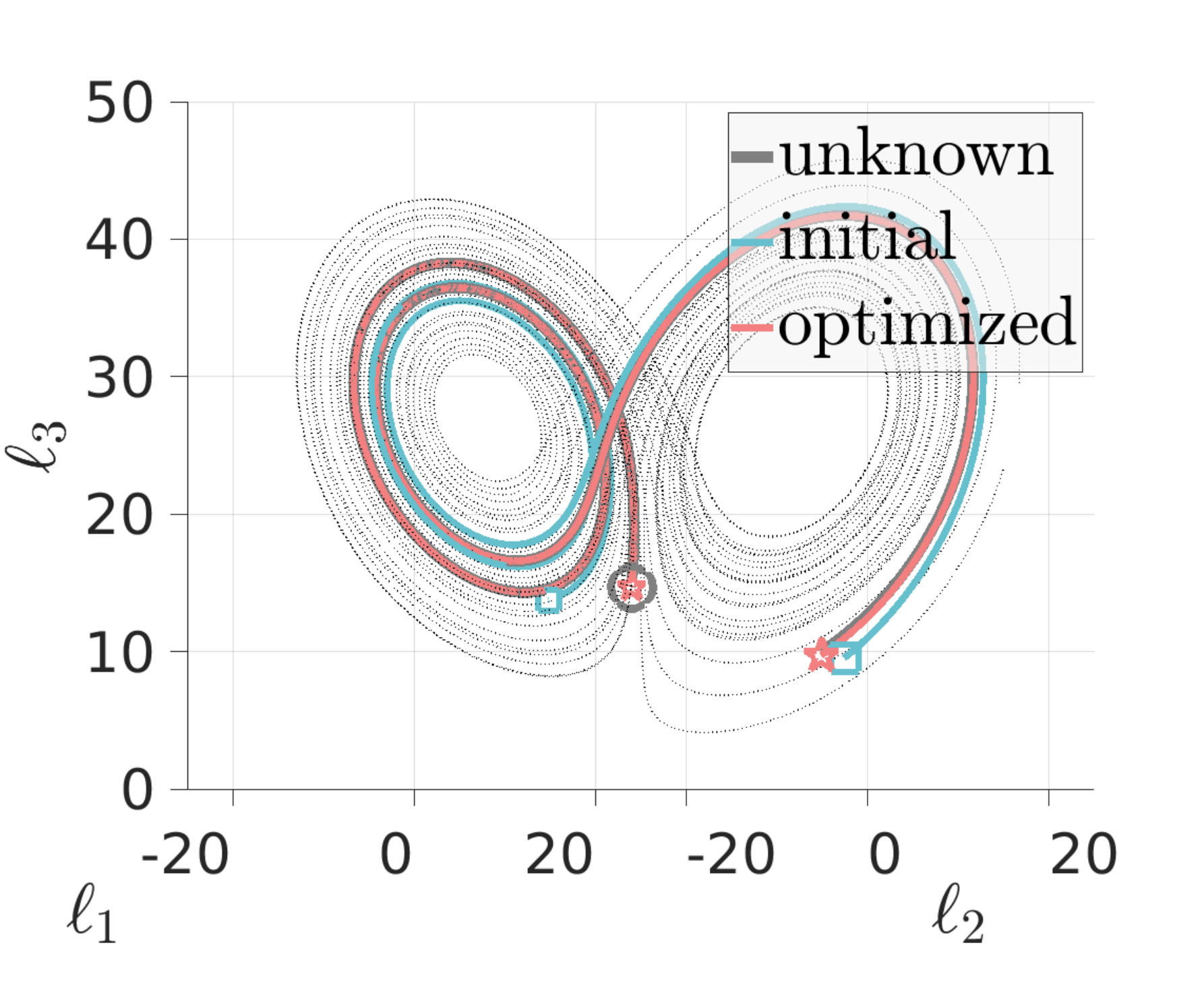}
        \caption{Trajectories of~\eqref{equ:Lorenz} starting from $\ell(0)$ (unknown), $x_0$ (initial) and $x_K$ for $K=10^5$ (optimized).}
        \label{fig:L2}
    \end{subfigure}
    \caption[]{{
    Estimating the initial state $\ell(0)$ of a Lorenz system~\eqref{equ:Lorenz} from a noisy measurement~$p$ of the state $\ell(2)=\varphi^2(\ell(0))$ (grey circle in~\ref{fig:L2}) at time~2. Even though the initial estimate $x_0$ is close to the optimized estimate~$x_K$, $\varphi^2(x_0)$ is far from $\varphi^2(\ell(0))$.  
    }}
    \label{fig:outlook}
\end{figure*}
\end{example}

{To close this section, we highlight that the complex-step method offers distinct benefits in the context of simulation-based optimization, where the objective function~$f$ can only be evaluated within prescribed error tolerances. In Example~\ref{ex:chaos}, for instance, the evaluation of~$f$ is corrupted by ODE integration errors. Unfortunately, such errors can have a detrimental impact on classical finite-difference-based optimization schemes. Indeed, the central-difference estimator $\frac{1}{2\delta_k
}(f(x_k+\delta_k y_k)-f(x_k-\delta_k y_k))y_k$ for $\nabla f(x_k)$ is useless for optimization unless the numerical errors in the evaluation of~$f$ are significantly smaller than~$\delta_k$. As $\delta_k$ must decay to $0$ as $k$ grows, so must the numerical tolerances. Otherwise, the ODE integration errors would dominate, which could be seen as another manifestation of catastrophic cancellation. Inexact evaluations of~$f$ can conveniently be modeled as outputs of a \textit{noisy} zeroth-order oracle. While this paper was under review, it has been shown that convergence guarantees for Algorithm~\ref{alg:convex_unconstrained} can be obtained even if the complex zeroth-order oracle is affected by independently and identically distributed noise and even if the sequence of smoothing parameters $\{\delta_k\}_{k\in \mathbb{Z}_{>0}}$ is chosen \textit{independently} of the noise statistics~\cite{ref:JongeneelIZO2021}. This provides strong evidence that the complex-step approach may be able to overcome the practical obstructions outlined above that plague classical finite-difference schemes in simulation-based optimization. As integration errors are arguably not purely random and serially independent, however, further research is needed.}

We highlight that complex-step derivatives are routinely used in PDE-constrained optimization. For example, they are used in a recent airfoil optimization package\footnote{\url{https://mdolab-cmplxfoil.readthedocs-hosted.com}.} developed in 2021. The underlying algorithm relies on \textit{sequential quadratic programming}~\cite[Ch.~18]{noc} and assumes that the complex-step derivative equals the gradient. In contrast, our analysis provides a rigorous treatment of approximation errors.

Additional applications of the complex-step derivative are discussed in~\cite[\S~3.2]{ref:martins2013review}.
}

\section{Conclusions and future work}
\label{sec:conclusions}
The cancellation effects that plague all multi-point gradient estimators tend to have a detrimental effect on the numerical stability and the convergence behavior of zeroth-order algorithms. These numerical problems can sometimes be mitigated by replacing the terminal iterate~$x_K$ with the averaged iterate~$\bar{x}_{K}=\tfrac{1}{K}\sum^{K}_{k=1}x_k$, at the cost of slower convergence. The single-point complex-step gradient estimator thus provides an attractive alternative to the classical gradient estimators because it leads to provably fast and numerically stable algorithms. As pointed out in~\cite{ref:Higham}, smoothness is not a necessary condition for the applicability of the complex-step approximation, which suggests that the analyticity assumption used in this paper can perhaps be relaxed. Other promising research directions would be to extend our convergence guarantees to the class of weakly convex functions and to investigate multi-batch as well as online settings. 

\appendix
\section{Appendix}

The proofs of our convergence results rely on the following lemma borrowed from~\cite{schmidt2011convergence}.

\begin{lemma}[{\cite[Lem.~1]{schmidt2011convergence}}]\label{lem:inequality_recursion}
If $\{t_k\}_{k\in\mathbb{Z}_{> 0}}$ and $\{\nu_k\}_{k\in\mathbb{Z}_{> 0}}$ are two sequence of non-negative real numbers, while $\{T_K\}_{K\in \mathbb{Z}_{> 0}}$ is a non-decreasing sequence of real numbers with $T_1\geq t_1^2$ such that $t_K^2 \le T_K + \textstyle\sum_{k=1}^{K} \nu_k t_k$ $\forall k\in\mathbb{Z}_{> 0}$, then we have
\[ t_K \le \textstyle\frac{1}{2} \sum_{k = 1}^{K} \nu_k + \left( T_K + (\frac{1}{2} \sum_{k = 1}^{K} \nu_k )^2 \right)^{\frac{1}{2}} \quad\forall K\in\mathbb{Z}_{> 0}.\]
\end{lemma}

In addition, several proofs in the main text make use of the inequalities
\begin{equation}
\label{equ:zeta}
    \textstyle\sum^J_{j=1}{j^{-2}}\leq \zeta(2)= \frac{1}{6}\pi^2\quad\text{and}\quad \sum^J_{j=1}{j^{-4}}\leq\zeta(4)= \frac{1}{90}\pi^4\quad \forall J\in \mathbb{Z}_{>0}, 
\end{equation}
which are obtained by truncating the series that defines the Riemann zeta function. 

\bigskip

\textbf{Acknowledgements}
WJ and DK are supported by the Swiss National Science Foundation under the NCCR Automation, grant agreement 51NF40\_180545. MCY is supported by the Hong Kong Research Grants Council under the grant 25302420. WJ wishes to  thank Prof. Arkadi Nemirovski for the supplied reference material, Prof. Timm Faulwasser for the pointer to the imaginary trick and Roland Schwan for fruitful discussions.


\pagestyle{basicstyle}
\subsection*{Bibliography}
\printbibliography[heading=none]


\end{document}